\documentclass[11pt,a4paper]{article}
\usepackage{graphicx}
\graphicspath{{figures/}} 
\usepackage[english]{babel}
\usepackage[export]{adjustbox}
\usepackage{amsfonts}
\usepackage{amsmath}
\usepackage{amssymb}
\usepackage{bbm}
\usepackage{epsfig}
\usepackage{color}
\usepackage{enumerate}
\usepackage{amsthm} 
\usepackage{soul}
\usepackage{enumitem}
\newtheorem{theorem}{Theorem}
\newtheorem{lemma}[theorem]{Lemma}
\newtheorem{corollary}[theorem]{Corollary}
\newtheorem{proposition}[theorem]{Proposition}
\newtheorem{remark}[theorem]{Remark}
\newtheorem{assumption}[theorem]{Assumption}
\newtheorem{definition}[theorem]{Definition}
\newtheorem{example}{Example}
\usepackage{caption}
\usepackage{float}
\usepackage{subcaption}
\usepackage{setspace}
\usepackage{amssymb, amsfonts, amsthm,amsmath,mathrsfs}
\usepackage{ifthen,bookmark}
\usepackage{graphicx}
\usepackage{epsfig}
\usepackage{lineno,hyperref}
\usepackage{booktabs}
\usepackage{multirow}
\usepackage{bm}
\usepackage[font={bf,footnotesize},textfont=md]{caption}
\usepackage{indentfirst}
\usepackage{array}
\usepackage{arydshln}
\usepackage{graphicx}
\usepackage{listings}
\usepackage{mathtools}
\hypersetup{hidelinks}
\usepackage{float}
\usepackage[labelfont=bf]{caption}

\definecolor{mygreen}{RGB}{0,120,80}

\newcommand{\bN}{\mathbb{N}}
\newcommand{\bR}{\mathbb{R}}
\newcommand{\bC}{\mathbb{C}}

\renewcommand{\Re}{\operatorname*{Re}}
\renewcommand{\Im}{\operatorname*{Im}}
\newcommand{\real}{\operatorname*{Re}}

\newcommand{\tol}{\mathrm{tol}}
\newcommand{\bone}{\mathbbm{1}}
\newcommand{\opi}{\mathrm{i}}

\newcommand{\e}{\mathrm{e}}

\def \Sv{{\mathbf{S}}}
\def \bv{{\mathbf{b}}}
\def \cv{{\mathbf{c}}}
\def \Wv{{\mathbf{W}}}
\def \Pv{{\mathbf{P}}}
\def \Rv{{\mathbf{R}}}
\def \Vv{{\mathbf{V}}}

\newcommand{\Iv}{{\mathbf{I}}}

\newcommand{\vv}{{\mathbf{v}}}
\newcommand{\Uv}{{\mathbf{U}}}

\newcommand{\Av}{{\mathbf{A}}}

\newcommand{\fv}{{\mathbf{f}}}
\newcommand{\ev}{{\mathbf{e}}}

\newcommand{\Ev}{{\mathbf{E}}}
\newcommand{\Yv}{{\mathbf{Y}}}
\newcommand{\Dv}{{\mathbf{D}}}

\newcommand{\lv}{{\mathbf{l}}}
\def \Qv{{\mathbf{Q}}}

\usepackage{enumerate}

\begin{document}
	
	\title{Runge--Kutta generalized Convolution Quadrature for sectorial problems}
	\author{J. Guo \thanks{ School of Mathematics and Statistics, Guangdong University of Technology, Guangdong, Guangzhou 510006, China. Email: {\tt jingguo@gdut.edu.cn}} \and M. Lopez-Fernandez\thanks{Department of Mathematical Analysis, Statistics and O.R., and Applied Mathematics.
			Fa\-cul\-ty of Sciences. University of Malaga.
			Bulevar Louis Pasteur, 31
			29010  Malaga,  Spain. Email: {\tt maria.lopezf@uma.es}}}
	\maketitle
	\begin{abstract}
		We study the application of the generalized convolution quadrature (gCQ) based on Runge--Kutta methods to approximate the solution of an important class of sectorial problems. The gCQ  generalizes Lubich's original convolution quadrature (CQ) to variable steps. High-order versions of the gCQ  have been developed in the last decade, relying on certain Runge--Kutta methods. The Runge--Kutta based gCQ  has been studied so far in a rather general setting, which includes applications to boundary integral formulations of wave problems. The  {more general} stability and convergence results for these new methods are suboptimal compared to those known for the uniform-step CQ, both in terms of convergence order and regularity requirements of the data.  Here  {we restrict the class of problems under study and} focus on a special  {type} of sectorial problems.  {We prove that for this important class of applications} it is possible to achieve the same order of convergence as for the original CQ, under the same regularity hypotheses on the data, and for very general time meshes. In the particular case of data with some known algebraic type of singularity, we also show how to choose an optimally graded time mesh to achieve convergence with maximal order, overcoming the well-known order reduction of the original CQ in these situations. 
		An important advantage of the gCQ  method is that it allows for a fast and memory-efficient implementation. We describe how the fast and oblivious Runge--Kutta based gCQ  can be implemented and illustrate our theoretical results with several numerical experiments. The codes implementing the examples in this article are available in \cite{CodesGuoLo25Zenodo}. To the best of our knowledge, comparable methods for the class of problems considered here do not currently offer the same combination of stability, high-order convergence, variable time steps, and reduced memory requirements and computational complexity.
	\end{abstract}
	
	

	\section{Introduction}
	
	In this work, we consider the high-order numerical approximation of Volterra-type convolution operators of the abstract form
	\begin{equation}\label{conv}
		u(t) = \int_0^tk(t-s) f(s)\, ds, \quad t > 0,
	\end{equation}
	on general non-uniform temporal meshes
	\begin{equation}\label{mesh}
		\Delta := \{0 = t_0 < t_1 < \cdots < t_N = T\},
	\end{equation}
	with step sizes defined as 
	\[\tau_n := t_n - t_{n-1},\quad n = 1, \ldots, N.\]
	The maximum and minimum step sizes are denoted by
	\begin{equation}\label{tau_maxmin}  
		\tau_{\max} := \max_{1 \leq j \leq N} \tau_j, \quad \tau_{\min} := \min_{1 \leq j \leq N} \tau_j.  
	\end{equation} 
	The  numerical discretization of \eqref{conv} has been extensively studied due to its  appearance in a wide range of models, including fractional differential equations \cite{BanLo19,JingRebecca_fi,Lu88II,Podlubny,McSloTho2006}, viscoelasticity \cite{Gorenflo2020} and time-domain boundary integral equations \cite{BanLuMe,Lu94,LuScha2002}. 	The convolution kernel \( k(t) \) can take various forms depending on the problem setting. Typical examples include scalar-valued kernels like \( k(t) = t^{\alpha} \) with \( \alpha > -1 \), which arise in fractional integral operators \cite{BanLo19,SamkoKilbasMarichev,Lu86}; operator-valued kernels such as \( k(t) = e^{\mathcal{A}t} \), representing semigroups generated by linear operators \cite{LoLuPaScha,CuLuPa}; and Green's functions in boundary integral equations \cite{BanLo20,SautSch}, including distributional kernels such as the Dirac delta that appears in wave propagation problems \cite[Chapter 2]{BanSay22}.
	
	Let $K(z)$ denote the Laplace transform of  $k(t)$. We assume that, {for the values of $z\in \bC$ where $K(z)$ is defined, it} acts as a bounded linear operator between normed vector spaces $\mathcal{X}$ and $\mathcal{Y}$, with norm
	\[
	\|K(z)\| := \sup_{\substack{u\in\mathcal{X}, \|u\|_{\mathcal{X}}=1}} \|K(z)u\|_{\mathcal{Y}}.
	\]
	The Laplace transform \(K(\cdot)\), which is referred to as the \emph{transfer operator} in many applications, plays a fundamental role in  Convolution Quadrature (CQ) methods, introduced in  Lubich's pioneering work \cite{Lu88I,Lu88II}. {The general theory of CQ methods and their gene\-ra\-li\-za\-tions to variable steps in \cite{LoSau13} and \cite{LoSau16} assumes that
		\begin{enumerate}
			\item There exists some abscissa $\sigma_{K} \in \bR$ such that $K(z)$ is well defined and the mapping
			\begin{equation}\label{def_Kop}
				K: \bC_{\sigma_{K}} \to \mathcal{L}(\mathcal{X},\mathcal{Y}) 
			\end{equation}
			is analytic, with the notation  $\bC_{\sigma}:=\{z\in \bC: \real z > \sigma\}$.
			\item For each $\sigma>\sigma_{K}$, there exists a constant $M_{\sigma}$ and some $\mu \in \bR$ such that
			\begin{equation}\label{gengrowthK}
				\|K(z)\|_{\mathcal{X}\to\mathcal{Y}} \leq M_{\sigma} (1 + |z|)^{\mu}, \qquad \forall z \in \bC_{\sigma}.
			\end{equation}
		\end{enumerate}
		This general framework covers important applications such as the boundary integral formulation of different types of hyperbolic equations \cite{Lu94, BanSay22}.}
	
	In this work, we reduce the class of problems and focus on a special type of {\em sectorial} problems, for which the Laplace transform of the convolution kernel admits a real inversion formula, see \cite{Henrici_II}. More precisely, we assume the following conditions on $K$:
	\begin{assumption}\label{assumptionK}
		The transfer operator $K$ has a branch cut along the negative real axis and satisfies:
		\begin{enumerate}  
			\item $K$ is holomorphic in the sector $|\arg(z)|<\pi$.
			\item There exist constants $M>0$ and $\alpha \in (0,1)$ such that
			\begin{equation}\label{Kz}
				\|K(z)\|_{\mathcal{X}\to\mathcal{Y}} \leq M |z|^{-\alpha}, \qquad |\arg(z)|<\pi.
			\end{equation}
			\item The operator $K$ extends continuously
				to the cut from both the upper and lower half-planes, except possibly at $z = 0$.
		\end{enumerate}  
	\end{assumption}

	When unambiguous, norm subscripts are omitted. 
	The conditions on \( K(z) \) cover a large and important range of examples, such as \( K(z) = z^{-\alpha} \), which corresponds to fractional integrals \cite{BanLo19,JingRebecca_fi,Lu86,Podlubny}; \( K(z) = (z^\alpha + r)^{-1} \) with \( r > 0 \), which is the Laplace transform of kernels usually used in damped wave models, see for instance \cite{BaBanPtas24}; also \( K(z) = (z^{\alpha} I + \mathcal{A})^{-1} \),  where \( \mathcal{A} \) is a symmetric positive definite operator or matrix \cite{GuoLo,JinLiZh,LiMa}.
	
	{A closely  related problem to the one we described above is the resolution of the integral equation \eqref{conv} for a given $u$, thus computing $f$ under the integral. Indeed, as proven in \cite[Theorem 28]{LoSau16}, the composition rule of the continuous convolution is preserved under discretization by the Runge--Kutta Generalized Convolution Quadrature method. Thus the results in this paper apply to the resolution of the integral equation by Runge--Kutta gCQ, as long as the inverse transfer operator $K^{-1}$ is well defined and belongs to the class in Assumption~\ref{assumptionK}. This is for instance the case when solving linear subdiffusion equations, see \cite[Section 4]{GuoLo}, which is one of our main applications.}
	
	Under Assumption~\ref{assumptionK}, the convolution kernel $k$ admits a real integral representation, cf. \cite{BanLo19, GuoLo}:	
	\begin{equation}\label{Gx}  
		k(t) = \int_0^\infty e^{-xt} G(x)\, dx, \quad 
		G(x) = \frac{K(x e^{-i\pi}) - K(x e^{i\pi})}{2\pi i},
	\end{equation}
	where $K(x e^{\pi i})$ and $K(x e^{-\pi i})$ denote the boundary values on the upper and lower edges of the branch cut, respectively, following the notation in \cite{Henrici_II}. The operator $G$ satisfies	the bound
	\begin{equation}\label{Gx_bnd}  
		\|G(x)\| \leq \frac{M}{\pi} x^{-\alpha}, \quad x > 0, \quad 0 < \alpha < 1.  
	\end{equation}  
	For the particular example of the fractional integral, with \( k(t) = t^{\alpha - 1}/\Gamma(\alpha) \), the corresponding function \( G(x) \) is given by
	\begin{equation}\label{Gx_fracint}
		G(x) =  {\frac{x^{-\alpha}  e^{i\pi \alpha} - x^{-\alpha} e^{-i\pi \alpha} }{2\pi i}}=\frac{\sin(\pi \alpha)}{\pi} x^{-\alpha}.
	\end{equation}
	For $K(z) = (z^{\alpha} I + \mathcal{A})^{-1}$, explicit formulas for $G$ are derived in \cite[Section 4]{GuoLo}. Other examples are considered in detail in Section~\ref{sec:num_test}.
	
	The representation in \eqref{Gx} has been successfully employed  in \cite{BanLo19}  for the efficient nu\-me\-ri\-cal solution of subdiffusion equations on uniform temporal grids, allowing for a special fast and oblivious quadrature which significantly improves the
	computational efficiency, the memory requirements and the coding.  In \cite{BanMak2022}, the Euler based gCQ weights for the fractional derivative are derived using \eqref{Gx_fracint}, and {\em a posteriori} error estimates for the application of this method to the subdiffusion equation are established. The study in \cite{GuoLo} is also based on \eqref{Gx}, where under Assumption~\ref{assumptionK}, optimal {\em a priori} error bounds are proven for the Euler based gCQ  method, together with a generalization of the fast and oblivious algorithm in \cite{BanLo19} to the case of general nonuniform meshes. So far all available studies based on \eqref{Gx} for the gCQ  are limited to the first order version of it. High order versions of the gCQ  based on Runge--Kutta methods have been proposed and analyzed in \cite{LoSau16} for more general problems than the ones considered in the present manuscript, for which high order but non optimal error estimates are proven, under very strong regularity assumptions on the data $f$.  {In this paper we generalize the ideas and theoretical developments in \cite{GuoLo} in order to derive sharper error estimates for the Runge--Kutta based gCQ for problems satisfying Assumption~\ref{assumptionK}. An essential tool in this error analysis is the use of the Mittag-Leffler function, see Appendix~\ref{sec:regularity}. To our knowledge, this theoretical approach was first introduced in \cite{GuoLo}.} 
	
	Original CQ methods, defined on uniform time meshes, require quite some regularity of the input data \( f(t) \) in order to achieve their maximal order. In particular, the extension of $f$ to $t<0$ by $0$ is typically required to be several times differentiable, which leads to the request that a certain number of the moments of $f$ at $t=0$ vanish \cite{LuOs,BanLu11}. In many practical applications this is not the case and a well-known order reduction takes place \cite{Lu88I}, even if correction terms can be added to improve the performance at times away from 0, see for instance \cite{BanMak2022,CuLuPa,JinLiZh,Lu2004}. Moreover, singularities of the data $f$ at other time points cannot be properly resolved by the original CQ formulas.  {These}  limitations have motivated the development and analysis of the so-called generalized Convolution Quadrature (gCQ) methods \cite{BanFer23,LoSau13,LoSau15apnum,LoSau16, GuoLo}. In particular, under Assumption~\ref{assumptionK}, a sharp error analysis is performed  in \cite{GuoLo}, for data $f$ satisfying the minimal regularity requirement 
	\begin{equation}\label{dataform}
		f(t) = t^\beta g(t), \quad \beta > -1,
	\end{equation}
	where \( g \) is smooth. In this work, the study in \cite{GuoLo} is generalized to high-order gCQ  based on Runge--Kutta methods. We actually improve further the regularity requirements on the data $f$ in \cite{GuoLo} by using distributional derivatives rather than the mean value theorem in the estimation of the local error remainders, see Remark~\ref{remark:improvegCQ1}. In this way, for problems satisfying Assumption~\ref{assumptionK}, we manage to achieve full order of convergence for the Runge--Kutta gCQ under the very same assumptions as the original Lubich's CQ with uniform step size, see \cite{LuOs}.  The results in \cite{LuOs} are asymptotically sharp and optimal with respect to both the smoothness of the data and the number of vanishing moments imposed at $t=0$.  Notice though that the class of problems we consider here is smaller, since we restrict to Assumption~\ref{assumptionK}, and that the state of the art theory for more general problems, including applications to hyperbolic time-domain integral equations, is still \cite{LoSau16}.
	
	{For the particular application to subdiffusion equations, we mention \cite{LiMa,LiMaRa}, where  high-order approximations of semilinear problems with a source term \( f \) that is nonsmooth at $t=0$ are studied. Both works are based on quadrature approximations of the contour integral representation of the solution on prescribed graded time grids. We notice that subdiffusion equations are only one of the applications of the theory developed here and that our convergence result Theorem~\ref{thm:conv_gCQrkGenM} is established for arbitrary time meshes, either prescribed {\em a priori} or dynamically generated by means of some step control mechanism. The particular case of graded meshes, for the whole class of problems satisfying Assumption 
		~\ref{assumptionK}, is carefully analyzed in Theorem~\ref{thm:conv_gCQrkGrad}. It is also worth noticing that, at least for linear subdiffusion equations, the preservation of the composition rule by the Runge--Kutta gCQ method allows only one matrix inversion per time step, as described in Subsection~\ref{subsec:subdiff}. This property avoids the evaluation of the multiple resolvents of the operator required by the contour integral representations used in the aforementioned works. The analysis of our method for semilinear subdiffusion equations is outside the scope of the present paper.
	}
	
	In what follows we consider \( s \)-stage Runge--Kutta schemes with coefficient matrix \( \Av \), weights \( \bv \),  and abscissas \( \cv \), satisfying the following conditions.
	\begin{assumption}\label{ass_RK}
		The Runge--Kutta method is A-stable, that is,
		\begin{equation}\label{Rz_Astab}
			|R(z)| \le 1 \quad \text{for } \Re z \le 0,
		\end{equation} 
		{where $R$ denotes the stability function of the method, given by \cite{HairWan10}
			\begin{equation}\label{Rz}
				R(z) := 1 + z \,\bv^\top (\Iv - z \Av)^{-1} \bone,
				\quad \text{with } \bone = (1,\dots,1,1)^\top \in \bR^s,
		\end{equation}}
		which   is further  assumed to be  a Pad\'e approximation to the exponential function. The method has classical order $p \ge 1$, stage order $q \le p$, and is stiffly accurate, meaning
		\begin{equation}\label{stiff_RK}
			\bv=\Av^\top\ev_s, \quad \text{with }\, \ev_s=(0,\cdots,0,1)^\top\in \bR^s.
		\end{equation}
		{
			The entries of the vector of abscissas $\cv = (c_1,\dots,c_s)^\top$ satisfy
			\begin{equation}\label{ci}   
				0\le c_1\le \cdots\le c_s\le 1.
		\end{equation}}
	\end{assumption}
	{\begin{remark}\label{remark:cbet}
			If $c_1=0$ in \eqref{ci}, then we must restrict to $\beta\ge0$ in \eqref{dataform}, since the numerical scheme in this case requires the value of $f$ at $0$, see \eqref{u_nweights}.
		\end{remark}
	}
	
	Notice that \eqref{Gx}  yields the alternative formulation of \eqref{conv}:
	\begin{equation*}\label{conv_real}
		u(t) = \int_0^\infty G(x) y(x,t)\, dx, 
	\end{equation*}
	where $y(x,t)$ satisfies
	\begin{equation}\label{ode}  
		\partial_t y(x,t) = -x y(x,t) + f(t), \quad y(x,0) = 0.
	\end{equation}  
	The gCQ  approximation \( u_n \) to $u(t)$ at $t_n$, for \( 1 \leq n \leq N \), is obtained by applying the $s$-stage Runge--Kutta   method to approximate $y(x,t_n)$, this is 
	\begin{equation*}
		u_n = \int_0^\infty G(x) y_n(x)\, dx, 
	\end{equation*}
	with $y_n(x) \approx y(x,t_n)$. This definition of $u_n$ is equivalent to the definition in \cite{LoSau16} for general problems \eqref{conv}, where a complex contour integral representation is needed for both $u(t)$ and its numerical approximation. {The available, general convergence theory for the Runge--Kutta gCQ in \cite{LoSau16} applies for problems satisfying Assumption~\ref{ass_RK}, with $\sigma_{K}=0$ in \eqref{def_Kop} and $\mu=-\alpha$ in \eqref{gengrowthK}. The precise convergence result from \cite{LoSau16} in this particular situation is given below.}   
	\begin{theorem}[Theorem 16 in \cite{LoSau16}]\label{th:convergence_regular}  
		Assume that	the Runge--Kutta method satisfies Assumption~\ref{ass_RK}, $f \in C^{p+2}([0,T])$ and $f^{(\ell)}(0) = 0$ for all $0 \leq \ell \leq q$.  Then, it holds
		\begin{equation}\nonumber
			\left\| u(t_n) - u_n \right\|_\mathcal{Y} \le C e^{\sigma T} \tau_{\max}^{q} \left\| f \right\|_{C^{p+2}([0,T],\mathcal{X})}, \quad n\ge 1.
		\end{equation}
	\end{theorem}  
	Here we analyze the error of the Runge--Kutta based gCQ  method applied to \eqref{conv}, with $f(t)$ like in \eqref{dataform}. More precisely, we assume that $f$ admits the fractional power expansion
	\begin{equation}\label{gen_f}
		f(t)= \sum_{\ell=0}^{m} \frac{f^{(\ell+\beta)}(0)}{\Gamma(\ell+\beta+1)} t^{\ell+\beta} + \frac{1}{\Gamma(m+\beta+1)} \left( t^{m+\beta} \ast f^{(m+\beta+1)}\right)(t),
	\end{equation}
	with the Riemann--Liouville derivative of order $\mu > 0$ defined as
	\[
	f^{(\mu)}(t) = \frac{d^\nu}{dt^\nu} \left(\int_0^t\frac{(t - s)^{\nu - \mu - 1}}{\Gamma(\nu - \mu)} f(s)\, ds\right),
	\]
	with $\nu := \lceil \mu \rceil$. The following theorem is the main result of this paper and establishes the convergence error bounds of the Runge--Kutta based gCQ  for problems satisfying Assumption~\ref{assumptionK}. 
	\begin{theorem}\label{thm:conv_gCQrkGenM}
		For \( f(t) \) admitting the expansion \eqref{gen_f} with \(\beta > -1\) and \( m = \lceil p - \beta -1\rceil \), and under Assumptions~\ref{assumptionK} and \ref{ass_RK}, the following error estimates hold true.
		\begin{enumerate}[label=(\roman*)]
			\item \label{itm:thm_genmesh} If  \( f^{(\ell+\beta)}(0) = 0 \) for all \( \ell = 0, 1, \dots, \lceil p - \alpha - \beta - 1 \rceil \), then on the general mesh \eqref{mesh}, the error satisfies
			\begin{align*}
				\left\Vert u(t_n) - u_n \right\Vert_{\mathcal{Y}} 
				&\le  {C} M  {
					\left(1+|\log (\tau_{\min})|\right)} \tau_{\max}^{\min\{p,q+1+\alpha\}} 
				\left\| f^{(m+\beta+1)} \right\|_{C^0([0,t_n], \mathcal{X})},\, n\ge 1.
			\end{align*}
			
			\item \label{itm:thm_quasimesh} If the mesh \eqref{mesh} is  {locally} quasi-uniform, i.e.,  
			\begin{equation}\label{quasimesh}
				c_{\Delta} := \frac{1}{2} \max_{2 \le i \le N} \left( \frac{\tau_i}{\tau_{i-1}} + \frac{\tau_{i-1}}{\tau_i} \right)
			\end{equation}  
			is a bounded constant, and  
			\begin{equation*}
				f^{(\ell + \beta)}(0) = 0 \quad \text{for all } \ell = 0, 1, \dots, \min\left\{\lceil p - \alpha - \beta - 1 \rceil, \lceil q - \beta \rceil\right\},
			\end{equation*}  
			then for all \( n \ge 1 \), the error satisfies
			\begin{align*}
				&\left\Vert u(t_n) - u_n \right\Vert_{\mathcal{Y}} \\[.5em]
				& \hspace{4em} \le  {C} M c_\Delta^{\max\{0,\, p- q-\alpha\}} 
				{
					\left(1+|\log (\tau_{\min})|\right)} \tau_{\max}^{\min\{p,q+1+\alpha\}} 
				\left\| f^{(m+\beta+1)} \right\|_{C^0([0,t_n], \mathcal{X})}.
			\end{align*}
		\end{enumerate}
		Here, \( {C}\) is a positive constant independent of the mesh  {which grows at most algebraically with respect to the final time $T$ and depends on $\alpha$ and $\beta$}.  
		
		{In particular, if we set \( \beta = 0 \), we see that the full convergence order \( \min\{p, q + 1 + \alpha\} \) is achieved when \( f \in C^{p}([0,T]; \mathcal{X}) \) and satisfies the vanishing moment conditions \( f^{(\ell)}(0) = 0 \) for all \( \ell = 0, 1, \dots, p-1 \), for the general mesh \eqref{mesh}, and  \( f^{(\ell)}(0) = 0 \) for all \( \ell = 0, 1, \dots, \min\{p-1,q \} \), for the locally quasi-uniform one.}
	\end{theorem}
	For graded meshes, which are  {locally} quasi-uniform, the error can be estimated more precisely, without assuming any vanishing moments of the data $f$ at the origin.
	\begin{theorem}\label{thm:conv_gCQrkGrad}
		Suppose $f(t)$ admits an expansion of the form \eqref{gen_f} with \( \beta > -1 \) and \( m = \lceil p - \beta -1\rceil \). Under Assumptions~\ref{assumptionK} and~\ref{ass_RK}, the error of the Runge--Kutta based gCQ  method for \eqref{conv} on the graded mesh
		\begin{equation}\label{gmesh}
			t_n= (n\tau)^{\gamma}, \qquad \tau=\frac{T^{1/\gamma}}{N},\quad n=1,\dots,N,\quad \gamma\ge 1,
		\end{equation}
		satisfies
		\begin{align*}
			\nonumber
			\left\Vert u(t_n) - u_n \right\Vert_{\mathcal{Y}} &\le C M \sum_{\ell = 0}^{m}  {T^{\alpha+\ell+\beta}}\log(N) \, N^{-\min\{ \gamma(\alpha + \ell + \beta),\, p,\, q + 1 + \alpha \}} \, \Vert f^{(\ell + \beta)}(0)\Vert_{\mathcal{X}} \\
			&\quad + CM {T^{\alpha+m+\beta+1}} \log(N) \, N^{-\min\{p,\, q + 1 + \alpha \}} \, \left\Vert f^{(m +  \beta+1)} \right\Vert_{C^0([0, t_n], \mathcal{X})},
		\end{align*}
		where $C$ is  a  positive constant independent of the mesh  {and the final time $T$.}   
		
	\end{theorem}
	
	
	{
		\begin{remark}
			The explicit form of the error constant $C$ in Theorem~\ref{thm:conv_gCQrkGenM} can be obtained from \eqref{Const_tnbet}, see Proposition~\ref{prop:conv_gCQrk}. It follows that all error bounds in Theorems~\ref{thm:conv_gCQrkGenM} and~\ref{thm:conv_gCQrkGrad} grow at most algebraically with respect to $T$, whereas in the state-of-the-art result of Theorem~\ref{th:convergence_regular} the growth is exponential.
		\end{remark}
	}
	{
		\begin{remark}
			When $\left\Vert f^{(\beta)}(0)\right\Vert_{\mathcal{X}}\neq 0$, the estimate in Theorem~\ref{thm:conv_gCQrkGrad} shows that, on the graded mesh \eqref{gmesh}, the error decreases for $\alpha+\beta>0$ with convergence rate $\min\{\gamma(\alpha+\beta),\, p,\, q+1+\alpha\}$, while for $\alpha+\beta\le 0$ the error grows  as $O\left(\log(N)\, N^{-\gamma(\alpha+\beta)}\right)$.
		\end{remark}
	}
	{
		\begin{remark} The results above can be easily generalized to the case of $f$ having an algebraic singularity at some $\sigma>0$ that can be isolated in the form $H(t-\sigma)(t-\sigma)^{\beta}$, for $H$ the Heaviside function. In this situation, appropriate meshes which are graded around $\sigma$ can still be devised in order to achieve the maximal order of convergence, see Example~\ref{ex_fracODE} in Subsection~\ref{subsec:expFDE}.
		\end{remark}
	}
	The paper is organized as follows.  
	Section~\ref{sec:rkgCQ} presents the framework of the Runge--Kutta  based gCQ  scheme.  
	Section~\ref{sec:stability} establishes the stability analysis, while Section~\ref{sec:err} derives error estimates for general and  {locally} quasi-uniform meshes, along with a gra\-ding stra\-tegy for graded meshes.  
	Implementation details, including a fast algorithm for computing convolution integrals and its application to subdiffusion equations, are discussed in Section~\ref{sec:algfi}.  
	The numerical experiments in Section~\ref{sec:num_test} evaluate the performance of the proposed algorithm for convolution integrals, subdiffusion equations, and nonlinear wave equations with damping terms of the form \eqref{conv}.  
	In all cases, the Runge--Kutta based gCQ  scheme achieves significantly reduced maximum absolute errors compared to the uniform step size approximation provided by original Lubich's CQ, cf.~\cite{JinLiZh} for the subdiffusion equation and \cite{BaBanPtas24} for the damped Westervelt equation.  
	In particular, for the nonlinear Westervelt equation considered in~\cite{BanFer}, employing the two-stage Radau--IIA gCQ  method on a suitable graded mesh increases the convergence order from $1+\alpha$ to $3$. Also, the achieved errors are perfectly consistent with our theoretical results and much better than predicted by the available theory to date in \cite{LoSau16}.

	\section{Runge--Kutta  based generalized convolution quadrature}\label{sec:rkgCQ}
	In this section, the formulation of Runge--Kutta based gCQ  on the real axis is developed. The $s$-stage Runge--Kutta method, applied to \eqref{ode} on the time mesh \eqref{mesh}, defines the stage vectors $\Yv_n(x)$ and solutions $y_n(x)$ via
	\begin{subequations}\label{RK_ode_vec}
		\begin{align}
			\Yv_{n}(x) &= y_{n-1}(x)\bone + \tau_n \Av \left(-x\Yv_{n}(x)+\fv_n\right), && n \geq 1, \\
			y_{n}(x) &= y_{n-1}(x) + \tau_n \bv^\top \left(-x\Yv_{n}(x)+\fv_n\right), && n \geq 1,
		\end{align}
	\end{subequations}
	where $\Yv_n = (Y_{ni})_{i=1}^s$,  $\fv_n = (f(t_{n-1}+c_i\tau_n))_{i=1}^s$. Rewriting \eqref{RK_ode_vec} as recurrences for $\Yv_n$ and $y_n$ yields
	\begin{equation*}
		\Yv_{n}(x) = \Rv(-\tau_n x) \ev_s^\top \Yv_{n-1}(x) + \tau_n \Av (\Iv + \tau_n x \Av)^{-1} \fv_n, \quad n \geq 1,
	\end{equation*}
	and
	\begin{equation}\label{RK_sol1}
		y_{n}(x) = R(-\tau_n x) y_{n-1}(x) + \tau_n \bv^\top (\Iv + \tau_n x \Av)^{-1} \fv_n, \quad n \geq 1,
	\end{equation}
	{ where
		\begin{equation}\label{RzVec}
			\Rv(z):= (\Iv - z \Av)^{-1} \bone,
		\end{equation} and $R(z)$ is given in  \eqref{Rz}.}
	Solving the recursion relations yields
	\begin{equation}\label{Yvn_sol}
		\Yv_{n}(x) = \sum_{j=1}^n \left( \prod_{l=j+1}^{n} \Rv(-\tau_l x) \ev_s^\top \right) \tau_j \Av (\Iv + \tau_j x \Av)^{-1} \fv_j, \quad n \geq 1,
	\end{equation}
	\begin{equation}\label{RK_sol}
		y_{n}(x) = \sum_{j=1}^n \left( \prod_{l=j+1}^{n} R(-\tau_l x) \right) \tau_j \bv^\top (\Iv + \tau_j x \Av)^{-1} \fv_j, \quad n \geq 1,
	\end{equation}
	where 
	\[
	\prod_{l=j+1}^{n} \Rv(-\tau_l x) \ev_s^\top = \Rv(-\tau_n x) \ev_s^\top \cdots \Rv(-\tau_{j+1} x) \ev_s^\top,
	\]
	with initial conditions $\Yv_0(x) = \mathbf{0}$ and $y_0(x) = 0$. The  Runge--Kutta based gCQ  approximation to \eqref{conv} is then defined as follows:
	\begin{definition}
		For $n\ge 1$, the Runge--Kutta based gCQ  approximation to \eqref{conv} on the time mesh \eqref{mesh} is defined by
		\begin{equation}\label{gCQ}
			u_n = \int_{0}^{\infty}G(x) y_n(x) \, dx,
		\end{equation}
		with $y_n(x)$ given in \eqref{RK_sol}.
	\end{definition}

	Let $\Uv(t_n):= \left(u(t_{n-1}+c_i\tau_n)\right)_{i=1}^s$. The Runge--Kutta based gCQ  approximation to $\Uv(t_n)$ is given by
	\begin{equation}\label{Un_stage}
		\Uv_n := \int_{0}^\infty G(x) \Yv_n(x)\,dx,\quad 1\le n\le N,
	\end{equation}
	where $\Yv_n(x)$ is given in \eqref{Yvn_sol}. Define the gCQ   weights
	\begin{equation}\label{gCQ_w}
		\Wv_{n,j}:= \tau_j  \int_{0}^\infty G(x) \left(\prod_{l=j+1}^{n} \left(\Rv(-\tau_l x)\ev_s^\top\right)\right) \left  (\Av(\Iv+\tau_jx\Av)^{-1} \right) \,dx,
	\end{equation}
	then the gCQ  scheme for the internal stages can be expressed as
	\begin{equation}\label{u_nweights}
		\Uv_n=\sum\limits_{j=1}^n \Wv_{n,j} \fv_j,\quad 1\le n\le N,
	\end{equation}
	where $u_n$ is the last component of $\Uv_n$. The weights \eqref{gCQ_w} also admit an alternative representation as contour integrals \cite{LoSau16}:
	\begin{equation}\label{gCQ_wContour}
		\Wv_{n,j}= \frac{\tau_j }{2\pi i} \int_{\Gamma} K(z) \left(\prod_{l=j+1}^{n} \left(\Rv(\tau_l z)\ev_s^\top\right)\right) \left  (\Av(\Iv-\tau_jz\Av)^{-1} \right) \,dz.
	\end{equation}

	\section{Stability}\label{sec:stability}
	The Runge--Kutta method satisfying Assumption~\ref{ass_RK} has the following properties:
	\begin{enumerate}
		\item By \eqref{Rz_Astab}, the eigenvalues of $\Av^{-1}$ lie in the region $\Re z > 0$, implying that $\Av$ is invertible and its eigenvalues are also in $\Re z > 0$.
		\item From \cite[Theorem 4.12]{HairWan10}, the matrix $\Av$ has $s$ distinct eigenvalues. Hence, $\Av$ is diagonalizable, i.e., there exists an invertible matrix $\Pv=(p_{ij})_{i,j=1}^s$ consisting of all eigenvectors of $\Av$ such that
		\begin{equation}\label{A_eigdecp}
				\Av=\Pv\Lambda \Pv^{-1},
		\end{equation}
		where
		\begin{equation}\label{lambda}
			\Lambda=\mathrm{diag}([\lambda_1,\cdots,\lambda_s]),\quad \lambda_{\min}=\min_{1\le i\le s} \Re{\lambda_i},
		\end{equation}
		and $\Av\Pv(:,i)=\lambda_i\Pv(:,i)$ for $1\le i\le s$.
		\item The function $R(z)$ satisfies
		\begin{align}\label{Rz_prop}
			R(\infty)=0, \qquad R(0)=R'(0)=1.
		\end{align}
	\end{enumerate}
	
	\begin{lemma}\cite{LoLuPaScha,Maria_thesis}\label{lem:Rz_bnd}
		Assume that the $s$-stage Runge--Kutta method satisfies Assumption \ref{ass_RK},  then there exists a constant $\hat{b}>0$ such that
		\begin{equation}\label{Rz_bnd}
			\left\vert R( - x) \right\vert\le
			\frac{1}{1+\hat{b}x},\quad \text{for all} \quad x > 0.
		\end{equation}
	\end{lemma}
	
	\begin{lemma}\label{lem:gCQ w_bnd}
		Under  Assumptions~\ref{assumptionK} and \ref{ass_RK}, the gCQ  weights given in \eqref{gCQ_w} satisfy
		\begin{align}
			\Vert \Wv_{n,j}\Vert\le CM\tau_j (t_n-t_{j-1})^{\alpha-1},\quad n\ge j\ge 1.
		\end{align}
	\end{lemma}
	\begin{proof}
		Applying the bound \eqref{Gx_bnd} to \eqref{gCQ_w} yields
		\begin{align*}
			&\Vert \Wv_{n,j}\Vert\le   \frac{M\tau_j }{\pi}  \int_{0}^\infty x^{-\alpha} \left\Vert \Av(\Iv+\tau_jx\Av)^{-1}\right\Vert\left\Vert \Rv(-\tau_n x)\ev_s^\top \right\Vert \left\vert\prod_{l=j+1}^{n-1} R(-\tau_l x)\right\vert\,dx,
		\end{align*}
		where $\Rv(\cdot)$ and $R(\cdot)$ are the stability functions defined in \eqref{RzVec} and \eqref{Rz}, respectively.
		For any \(x>0\), the eigenvalue decomposition \eqref{A_eigdecp} gives
		\begin{equation*}
			\left\Vert \Av(\Iv+\tau_jx\Av)^{-1}\right\Vert
			\le
			\frac{\|\Pv^{-1}\|\,\|\Pv\|\max\limits_{1\le i\le s}|\lambda_i|}
			{1+\min\limits_{1\le i\le s}\{\Re\lambda_i\}\tau_jx},
		\end{equation*}
		and
		\begin{align*}
			\left\Vert\Rv(-\tau_nx)\ev_s^\top\right\Vert
			&=\left\Vert(\Iv+\tau_nx\Av)^{-1}\bone\ev_s^\top\right\Vert
			\le
			\frac{\|\Pv^{-1}\|\,\|\Pv\|\sqrt{s}}
			{1+\min\limits_{1\le i\le s}\{\Re\lambda_i\}\tau_nx}.
		\end{align*}
		Hence,
		\begin{equation*}
			\left\Vert\Av(\Iv+\tau_jx\Av)^{-1}\right\Vert
			\left\Vert\Rv(-\tau_nx)\ev_s^\top\right\Vert
			\le
			\frac{\|\Pv^{-1}\|^2\|\Pv\|^2\sqrt{s}\max\limits_{1\le i\le s}|\lambda_i|}
			{1+\min\limits_{1\le i\le s}\{\Re\lambda_i\}(\tau_n+\tau_j)x}.
		\end{equation*}
		In addition, Lemma~\ref{lem:Rz_bnd} provides
		\begin{equation*}
			\left|\prod_{l=j+1}^{n-1}R(-\tau_lx)\right|
			\le \frac{1}{1+\hat b(t_{n-1}-t_j)x},
			\qquad n\ge j+1.
		\end{equation*}
		The case \(n=j\), for which the product is empty, follows directly from the preceding estimate. Let
		\[
			\hat v=\min\left\{\hat b,\min\limits_{1\le i\le s}\{\Re\lambda_i\}\right\}.
		\]
		Combining these estimates and using \eqref{int_Betafun}, we obtain
		\begin{align*}
			\Vert\Wv_{n,j}\Vert
			&\le \frac{M\tau_j}{\pi}
			\int_0^\infty x^{-\alpha}
			\frac{\|\Pv^{-1}\|^2\|\Pv\|^2\sqrt{s}\max\limits_{1\le i\le s}|\lambda_i|}
			{1+\hat v(t_n-t_{j-1})x}\,dx\\[.3em]
			&=\frac{M\tau_j}{\pi}\|\Pv^{-1}\|^2\|\Pv\|^2\sqrt{s}B(\alpha,1-\alpha)
			\max\limits_{1\le i\le s}|\lambda_i|
			\bigl(\hat v(t_n-t_{j-1})\bigr)^{\alpha-1}.
		\end{align*}
		{where $B(\cdot,\cdot)$ denotes the Beta function \cite[(1.28)]{Podlubny}.} The proof is thus complete.
	\end{proof}
	The stability of scheme \eqref{gCQ} follows directly from the preceding weight estimates, as formalized in the following theorem.
	\begin{theorem}\label{thm:stab}
		Under Assumptions~\ref{assumptionK} and \ref{ass_RK}, the Runge--Kutta based gCQ scheme \eqref{gCQ} satisfies
		\begin{align*}
			\|u_n\|_{\mathcal{Y}} \le CM \sum_{j=1}^n \tau_j \left(t_n - t_{j-1}\right)^{\alpha-1} \left\| f_j \right\|_{C^0\left ([t_{j-1},t_j];\mathcal{X}\right)}.
		\end{align*}
	\end{theorem}
	
	\section{Proof of Theorems~\ref{thm:conv_gCQrkGenM} and \ref{thm:conv_gCQrkGrad}}\label{sec:err} 
Following Lubich's operational framework \cite{Lu88I},  we adopt the following notation for the error analysis:
\begin{align*}  
	&K(\partial_t)f := k * f, \quad \text{ the continuous convolution defined in \eqref{conv}}, \\  
	&	K(\partial_t^\Delta)f,\quad  \text{the Runge--Kutta based gCQ  approximation of \eqref{conv} on the time mesh  } \Delta .  
\end{align*}   
\subsection{Error representation on the real axis}
Let $\Yv(x, t_n):=(y(x,t_{n-1}+c_i\tau_n))_{i=1}^s$. Substituting the exact solution into the Runge--Kutta scheme \eqref{RK_ode_vec} gives
\begin{subequations}\label{RK_ode_defec}
	\begin{align}
		\label{Y_tn}
		\Yv(x, t_n) &= y(x, t_{n-1})\bone + \tau_n \Av \left(-x\Yv(x, t_{n})+\fv_n\right) + \Dv_n(x), &&\text{for } n \ge 1, \\
		\label{y_tn}
		y(x, t_n) &= y(x, t_{n-1}) + \tau_n \bv^\top \left(-x\Yv(x, t_{n})+\fv_n\right) + d_n(x), &&\text{for } n \ge 1,
	\end{align}
\end{subequations}
where \(\Dv_n = (d_{ni})_{i=1}^s\) and \(d_n\) denote the defects.  
Under Assumption~\ref{ass_RK}, these defects admit the representations \cite[proof of Theorem 3.3]{LuOs}:
\begin{align}
	\label{dni}
	d_{ni}(x) &= \sum\limits_{\ell=q+1}^p \delta_{i}^{(\ell)} \tau_n^\ell \frac{\partial^{\ell}y}{\partial t^{\ell}}(x,t_{n}) + \tau_n^p \int_{t_{n-1}}^{t_n} \kappa_i\left(\frac{t-t_{n-1}}{\tau_n}\right) \frac{\partial^{p+1}y}{\partial t^{p+1}}\,dt, \\
	\label{dn}
	d_n(x) &= \tau_n^p \int_{t_{n-1}}^{t_n} \kappa\left(\frac{t-t_{n-1}}{\tau_n}\right) \frac{\partial^{p+1}y}{\partial t^{p+1}}\,dt,
\end{align}
where $\kappa_i$ and $\kappa$ are bounded Peano kernels, and
{	\begin{equation}\label{delt_ell}
		\delta^{(\ell)} := (\delta_i^{(\ell)})_{i=1}^s,
		\qquad
		\delta_i^{(\ell)}
		= \frac{1}{(\ell-1)!}
		\left(
		\frac{c_i^\ell}{\ell}
		- \sum_{j=1}^s a_{ij} c_j^{\ell-1}
		\right),
		\quad q+1 \le \ell \le p .
\end{equation}}
Denote $\Ev_n(x):=\Yv(x, t_n) - \Yv_n(x)$ and $e_n(x):=y(x, t_n) - y_n(x)$. 
Subtracting \eqref{RK_ode_vec} from \eqref{RK_ode_defec} gives
\begin{equation}
	\label{RK_err}
	\begin{alignedat}{2}
		\Ev_n(x) &= e_{n-1}(x)\bone - \tau_n x \Av \Ev_n(x) + \Dv_n(x), &&\quad \text{for } n \ge 1, \\
		e_n(x) &= e_{n-1}(x) - \tau_n x \bv^\top \Ev_n(x) + d_n(x), &&\quad \text{for } n \ge 1.
	\end{alignedat}
\end{equation}
Solving the recursion \eqref{RK_err}, we obtain
\begin{equation}\label{en}
	e_n(x) = \sum\limits_{j=1}^n \left(d_j(x)-\tau_j x \bv^\top(\Iv+\tau_j x \Av)^{-1} \Dv_j(x) \right) \prod_{l=j+1}^{n} R(-\tau_l x),
\end{equation}
with $e_0(x) = 0$ and $R(z)$ defined in \eqref{Rz}. Then, the error representation of the  Runge--Kutta  based gCQ  scheme \eqref{gCQ} can be expressed as
\begin{equation}
	\label{gCQ_error_decomposition}
	u(t_n) - u_n = \int_{0}^{\infty} G(x) e_n(x) \, dx :=\mathcal{E}_n^1 + \mathcal{E}_n^2,
\end{equation}
where
\begin{align}
	\label{Err_1}
	\mathcal{E}_n^1 &= \int_0^{\infty} G(x) \left(   d_1(x)-\tau_1 x \bv^\top(\Iv + \tau_1 x \Av)^{-1} \Dv_1(x) \right) \left(\prod_{j=2}^{n} R(-\tau_j x)\right)\,dx,\\[.5em]
	\label{Err_2}
	\mathcal{E}_n^2 &= \sum_{j=2}^{n} \int_0^{\infty} G(x)\left( d_j(x)-\tau_j x \bv^\top(\Iv + \tau_j x \Av)^{-1} \Dv_j(x)\right) \left(\prod_{l=j+1}^{n} R(-\tau_l x)\right)\,dx.
\end{align}


\subsection{Convergence analysis for   power-type integrands }
To establish error estimates for general integrands~\( f \) satisfying~\eqref{gen_f}, we begin by analyzing the  case \( f(t) = t^\beta \mathbf{v} \) with \( \beta > -1 \). The error analysis proceeds by examining each component of the error decomposition~\eqref{gCQ_error_decomposition} on the general mesh~\eqref{mesh}.

\subsubsection{Error estimate for the first time interval}
We now establish the error bound for    $\mathcal{E}_n^1 $ in \eqref{Err_1}.
\begin{lemma}\label{lem:Err_1}
	Given  \(f(t) = t^\beta \vv\) with \(\beta > -1\), under Assumptions~\ref{assumptionK} and \ref{ass_RK}, the approximation error \(\mathcal{E}_n^1\) in \eqref{Err_1} on the mesh \eqref{mesh} satisfies
	\begin{equation}\label{T1_bnd}
		\left\Vert \mathcal{E}_n^1 \right\Vert_{\mathcal{Y}} \leq C M \tau_1^{\alpha + \beta},
	\end{equation}
	where the positive constant \(C\) is independent of the mesh  {and the final time $T$.}
\end{lemma}

\begin{proof}
	{Taking $n=1$ in \eqref{RK_ode_defec}, we have
		\begin{align*}
			\Yv(x, t_1) &=   \tau_1 \Av \left(-x\Yv(x, t_{1})+\fv_1\right) + \Dv_1(x),\\
			y(x, t_1) &=  \tau_1 \bv^\top \left(-x\Yv(x, t_{1})+\fv_1\right) + d_1(x).
		\end{align*}
		Solving the first equation for $\Yv(x, t_1)$ gives
		\[ \Yv(x, t_1) = (\Iv + \tau_1 x \Av)^{-1}(\tau_1 \Av \fv_1 + \Dv_1(x)).\]
		Substituting the above expression  into the equation for $y(x, t_1)$, we obtain
		\begin{align*}
			y(x, t_1) &= \tau_1 \bv^\top \left(-x(\Iv + \tau_1 x \Av)^{-1}(\tau_1 \Av \fv_1 + \Dv_1(x)) + \fv_1\right) + d_1(x)\\
			&=-\tau_1 x \bv^\top(\Iv + \tau_1 x \Av)^{-1} \Dv_1(x)+\tau_1 \bv^\top\left(-\tau_1 x \Av(\Iv + \tau_1 x \Av)^{-1}\fv_1+\fv_1\right)+d_1(x).
		\end{align*}
		Consequently, the defect-related term in \eqref{Err_1} can be written as
		\begin{align*}
			d_1(x)-\tau_1 x \bv^\top(\Iv + \tau_1 x \Av)^{-1} \Dv_1(x)
			&=y(x, t_1)-\tau_1 \bv^\top\left(-\tau_1 x \Av(\Iv + \tau_1 x \Av)^{-1}\fv_1+\fv_1\right)\\
			&=y(x, t_1)-\tau_1 \bv^\top(\Iv + \tau_1 x \Av)^{-1}\fv_1. 
		\end{align*}
		Combining the above expression with \eqref{Rz_Astab} leads to
		\begin{align*}
			\left\Vert \mathcal{E}_n^1 \right\Vert_{\mathcal{Y}} 
			&\leq C M \int_{0}^\infty x^{-\alpha} \left\| y(x, t_1)-\tau_1 \bv^\top(\Iv + \tau_1 x \Av)^{-1}\fv_1\right\|_{\mathcal{X}}  \, dx\\
			&\le C M \int_{0}^\infty x^{-\alpha} \left(\left\| y(x, t_1)\right\|_{\mathcal{X}}+\left\|\tau_1 \bv^\top(\Iv + \tau_1 x \Av)^{-1}\fv_1\right\|_{\mathcal{X}} \right) \, dx.
		\end{align*}
		For \(\fv_1 =\left ((c_i\tau_1)^\beta\vv\right)_{i=1}^s\) with \(c_1\) in \eqref{ci} as specified in Remark~\ref{remark:cbet} (i.e., \(c_1 \neq 0\) for \(-1 < \beta < 0\)), applying \eqref{ydiff_bnd} with \(m=0\) and \eqref{A_eigdecp} yields
		\begin{align*}
			\left\Vert \mathcal{E}_n^1 \right\Vert_{\mathcal{Y}} 
			&\leq C M \tau_1^{\beta+1} \int_0^\infty x^{-\alpha} \frac{1}{1 +  \tau_1 x} \, dx 
			= C M B(\alpha, 1 - \alpha) \tau_1^{\alpha + \beta},
		\end{align*}
		where the integral is evaluated using \eqref{int_Betafun}. This completes the proof.}
\end{proof}

\subsubsection{Error estimates from the second time interval}
This subsection focuses on establishing an error bound for the accumulated error \( \mathcal{E}_n^2 \) defined in \eqref{Err_2}, for $n\ge 2$. The analysis begins with the following lemma.
\begin{lemma}\label{lem:Matstab_estimate}
	Under Assumption~\ref{ass_RK}, for \(\tau > 0\), \(x > 0\), and \(\delta^{(\ell)}\) defined in \eqref{delt_ell} with \(q + 1 \leq \ell \leq p\), it holds
	\begin{equation}\label{stab_combined}
		\left| \tau x\bv^\top(\Iv + \tau x\Av)^{-1} \delta^{(\ell)}\right| \leq 
		\begin{cases}
			C(\tau x)^{p-\ell+1}, & 0 < x < (\tau\|\Av\|)^{-1}, \\
			C, & x \geq (\tau\|\Av\|)^{-1},
		\end{cases}
	\end{equation}
	where \(C > 0\) is a constant independent of \(\tau\) and \(x\), and \(\|\mathbf{A}\|\) denotes the induced matrix norm.
\end{lemma}
\begin{proof}
	We first establish a uniform bound valid for all $x > 0$, then refine it for small $x$. From the stiff accuracy condition \eqref{stiff_RK} and the eigendecomposition \eqref{A_eigdecp}, we have
	\begin{align}
		\label{mat_bnd}
		\left\|\tau x\bv^\top(\Iv + \tau x\Av)^{-1}\right\| 
		= \left\|\ev_s^\top - \ev_s^\top(\Iv + \tau x\Av)^{-1}\right\| \leq 1 + \frac{\|\Pv\|\,\|\Pv^{-1}\|}{1 + \tau x\lambda_{\min}} \leq 1+\|\Pv\|\,\|\Pv^{-1}\|,
	\end{align}
	where $\lambda_{\min}$ is given in \eqref{lambda}. Thus, for \emph{any} $x > 0$ 
	\[
	\left| \tau x\bv^\top(\Iv + \tau x\Av)^{-1} \delta^{(\ell)}\right| 
		\leq \bigl(1+\|\Pv\|\,\|\Pv^{-1}\|\bigr)\left\Vert \delta^{(\ell)}\right\Vert \leq C.
	\]
	When \(0 < x < (\tau\|\Av\|)^{-1}\), the Neumann series expansion gives
	\begin{equation*}
		(\Iv + \tau x\Av)^{-1} = \sum_{m=0}^\infty (-\tau x\Av)^m.
	\end{equation*}
	Together with the order conditions
	\begin{equation*}
		\bv^\top\Av^{m+1}\cv^{\ell-1} = \frac{1}{\ell}\bv^\top\Av^m\cv^\ell \quad \text{for } m+1+\ell \leq p,
	\end{equation*}
	{where	$\cv^\ell=(c_i^\ell)_{i=1}^s$}, and the definition of \(\delta^{(\ell)}\) in \eqref{delt_ell}, it follows for each \(q+1 \leq \ell \leq p\) that
	\begin{equation*}
		\left| \tau x\bv^\top(\Iv + \tau x\Av)^{-1} \delta^{(\ell)}\right|
		= \left| \tau x\bv^\top \sum_{m=p-\ell}^\infty (-\tau x\Av)^m \delta^{(\ell)}\right|
		\leq C(\tau x)^{p-\ell+1}.
	\end{equation*}
	This completes the proof.
\end{proof}
\begin{lemma}\label{lem:integral_bound}
	For the defects \(d_j(x)\) in \eqref{RK_ode_defec} with \(f(t) = t^\beta \vv\), where \(\beta > -1\), the following bound holds on the general mesh \eqref{mesh}:
	\begin{equation}\label{intdj_bnd}
		\sum_{j=2}^n \int_0^\infty x^{-\alpha}		\left\| d_j(x)  \right\|_{\mathcal{X}}\, dx \leq C
		\begin{cases}
			{\log\left(\frac{t_n}{\tau_1}\right)\max\limits_{2\le j\le n}\left\{\tau_j^p t_j^{\alpha+\beta - p}\right\}}, & \text{if } \beta \ge  p - \alpha, \\[.8em]
			\sum\limits_{j=2}^n \tau_j^{p+1} \xi_j^{\alpha + \beta - p - 1}, & \text{if } \beta < p - \alpha,
		\end{cases}
	\end{equation}
	where $2\le n\le  N$,    $\xi_j\in(t_{j-1},t_j)$ and  $C$ is a positive constant independent of the mesh.
\end{lemma}

\begin{proof}
	Using  \eqref{dn}, we have
	\begin{equation*}
		\sum_{j=2}^n \int_0^\infty x^{-\alpha}		\left\| d_j(x)  \right\|_{\mathcal{X}}\, dx \leq C \int_0^\infty x^{-\alpha}  \left( \sum_{j=2}^n \tau_j^p\int_{t_{j-1}}^{t_j} \left\| \frac{\partial^{p+1} y}{\partial t^{p+1}} \right\|_{\mathcal{X}} \, dt \right)\, dx.
	\end{equation*}
	From \eqref{ydiff_bnd}, it follows that
	\begin{equation}\label{intdj_allj}
		\sum_{j=2}^n \tau_j^p \int_{t_{j-1}}^{t_j} \left\| \frac{\partial^{p+1} y(t)}{\partial t^{p+1}} \right\|_{\mathcal{X}} \, dt
		\le C \Gamma(\beta + 1) \sum_{j=2}^n \tau_j^p \int_{t_{j-1}}^{t_j} \frac{t^{\beta - p}}{1 + x t} \, dt, \quad \text{for } x > 0.
	\end{equation}
	Hence,
	\begin{equation*}
		\sum_{j=2}^n\int_0^\infty  x^{-\alpha}  	\left\|d_j(x) \right\|_{\mathcal{X}} \, dx\leq C\Gamma(\beta+1) \sum\limits_{j=2}^n \tau_j^p \int_{t_{j-1}}^{t_j} t^{\beta - p} \left( \int_0^\infty \frac{x^{-\alpha}}{1 + x t} \, dx \right) dt.
	\end{equation*}
	By applying the definition of the Beta function given in \eqref{int_Betafun}, we derive
	\begin{equation*}
		\sum_{j=2}^n\int_0^\infty  x^{-\alpha}  	\left\|d_j(x) \right\|_{\mathcal{X}} \, dx \leq C\Gamma(\beta+1)B(\alpha,1-\alpha) \sum\limits_{j=2}^n \tau_j^p \int_{t_{j-1}}^{t_j} t^{\alpha+\beta - p-1} \, dt.
	\end{equation*}
	The final bound \eqref{intdj_bnd} follows from the estimate  {
		\begin{equation*}
			\sum_{j=2}^n \tau_j^p \int_{t_{j-1}}^{t_j} t^{\alpha+\beta - p-1} \, dt
			\le
			\left(\sum_{j=2}^n \int_{t_{j-1}}^{t_j} t^{-1} \, dt\right)
			\max_{2\le j\le n}\left\{\tau_j^p t_j^{\alpha+\beta - p}\right\}
		\end{equation*}
		for the case $\beta \ge p-\alpha$, and from the integral mean value theorem for $-1 < \beta < p-\alpha$, } thereby completing the proof.
\end{proof}

\begin{lemma}\label{lem:Pj_bnd}
	Suppose the \(s\)-stage Runge--Kutta method satisfies Assumption~\ref{ass_RK}. Define
	\begin{equation}
		\label{Pj}
		P_j := \sum_{\ell = q+1}^{p} \int_0^\infty x^{-\alpha} \left\| \tau_j x \bv^\top (\Iv + \tau_j x \Av)^{-1} \delta^{(\ell)} \tau_j^\ell \frac{\partial^\ell y}{\partial t^\ell}(x, t_j) \right\|_{\mathcal{X}} \prod_{l=j+1}^n \left| R(-\tau_l x) \right| \, dx,
	\end{equation}
	for \(2 \leq j \leq n\), where \(y(x,t)\) is the solution to \eqref{ode} with \(f(t) = t^\beta \vv\), \(\beta > -1\). Then, on the general time mesh \eqref{mesh},   it holds
	\begin{equation}
		\label{Pj_bnd}
		\sum_{j=2}^n P_j \le C  {\log\left(2+\frac{t_n}{\tau_n}\right)}
		\begin{dcases}
			\tau_{\max}^{\min\{p,q+1+\alpha\}} t_n^{\alpha+\beta-\min\{p,q+1+\alpha\}}, & \text{if } \beta \ge p-\alpha, \\[.5em]
			\max_{2 \leq j \leq n} \left\{  {\tau_{j}^{\min\{p,q+1+\alpha\}} t_j^{\alpha+\beta-\min\{p,q+1+\alpha\}}} \right\}, & \text{if } \beta < p-\alpha,
		\end{dcases}
	\end{equation}
	where $C$ is a positive constant independent of the mesh.
\end{lemma}
\begin{proof}
	By Lemma~\ref{lem:Rz_bnd} and the estimate in \eqref{ydiff_bnd}, we have
	\begin{align*}
		P_j &\le C \sum_{\ell = q+1}^p \tau_j^\ell t_j^{\beta - \ell + 1} \int_0^\infty x^{-\alpha} \left| \tau_j x \bv^\top (\Iv + \tau_j x \Av)^{-1} \delta^{(\ell)} \right| \frac{1}{(1 + t_j x)(1 + \hat{b}(t_n - t_j)x)}\, dx.
	\end{align*}
	For the case \(j=n  \), Lemma~\ref{lem:Matstab_estimate} yields
	\begin{align*}
		P_n \le& C \sum_{\ell = q+1}^{p-1}  \left(\tau_n^{p+1}t_n^{\beta - \ell + 1}\int_0^{1 / (\tau_n \Vert \Av \Vert)} \frac{x^{-\alpha+p-\ell+1}}{1 + t_n x}\, dx+ \tau_n^\ell t_n^{\beta - \ell + 1}\int_{1 / (\tau_n \Vert \Av \Vert)}^\infty \frac{x^{-\alpha}}{1 + t_n x}\, dx\right) \\[.5em]
		&+C\tau_n^p t_n^{\beta-p+1}\int_0^\infty\frac{x^{-\alpha}}{1 + t_n x}\, dx,
	\end{align*}
	which gives
	\begin{equation}\label{Pn_esti}
		P_n \le C  \sum_{\ell = q+1}^{p-1} \tau_n^{\alpha + \ell} t_n^{\beta - \ell}+C\tau_n^p t_n^{\alpha+\beta-p}.		\end{equation}
	Using 
	\begin{equation}\label{ineq:stepsize}
		\tau_j / t_j < 1 \quad \text{for all } j \ge 2,
	\end{equation}
	we further deduce that
	\[
	P_n \le C \left(\tau_n^{q+1+\alpha} t_n^{\beta - q - 1} + \tau_n^p t_n^{\alpha + \beta - p}\right).
	\]
	Now consider the case \( 2 \le j \le n - 1\). In this case, we have
	\[P_j \le C \sum_{\ell = q+1}^p \tau_j^\ell(\hat{b}(t_n - t_j))^{-1}  t_j^{\beta - \ell +1} \int_0^\infty x^{-\alpha-1} \left| \tau_j x \bv^\top (\Iv + \tau_j x \Av)^{-1} \delta^{(\ell)} \right| \frac{1}{1 + t_j x}\, dx.\]
	Applying Lemma~\ref{lem:Matstab_estimate}  gives
	\begin{align} \label{Tn3_m1p1}
		P_j \le
		C \sum_{\ell = q+1}^p \left( |I^{(\ell)}_{\mathrm{near}}| + \tau_j^{\alpha + \ell + 1} (\hat{b}(t_n - t_j))^{-1} t_j^{\beta - \ell} \right),
	\end{align}
	where
	\[
	I^{(\ell)}_{\mathrm{near}} := \tau_j^{p+1} (\hat{b}(t_n - t_j))^{-1} t_j^{\beta - \ell + 1} \int_0^{1 / (\tau_j \Vert \Av \Vert)} x^{-\alpha + p - \ell} \frac{1}{1 + t_j x}\, dx.
	\]
	Evaluating the integral explicitly yields
	\begin{align} \label{Inear_combined}
		|I^{(\ell)}_{\mathrm{near}}| \le
		\begin{cases}
			C \tau_j^{\alpha + \ell + 1} (\hat{b}(t_n - t_j))^{-1} t_j^{\beta - \ell}, & q+1 \le \ell \le p - 1, \\[.7em]
			B(\alpha, 1 - \alpha) \tau_j^{p + 1} (\hat{b}(t_n - t_j))^{-1} t_j^{\alpha + \beta - p}, & \ell = p.
		\end{cases}
	\end{align}
	Substituting \eqref{Inear_combined} into \eqref{Tn3_m1p1}, we obtain
	\begin{align}
		\label{Pj_esti}
		P_j &\le C \sum_{\ell = q+1}^{p-1} \tau_j^{\alpha + \ell + 1} (\hat{b}(t_n - t_j))^{-1} t_j^{\beta - \ell}
		+ C \tau_j^{p + 1} (\hat{b}(t_n - t_j))^{-1} t_j^{\alpha + \beta - p} \\[.5em]
		\nonumber
		&\le C (\hat{b}(t_n - t_j))^{-1} \left( \tau_j^{q+2 + \alpha} t_j^{\beta - q - 1} + \tau_j^{p+1} t_j^{\alpha + \beta - p} \right),
	\end{align}
	where \eqref{ineq:stepsize} is applied in the last step. Combining all cases yields
	\[
	\sum_{j=2}^n P_j \le C \left(1 + \sum_{j=2}^{n-1} \tau_j (\hat{b}(t_n - t_j))^{-1} \right) \max_{2 \le j \le n} \left\{ \tau_j^{q+1+\alpha} t_j^{\beta - q - 1} + \tau_j^p t_j^{\alpha + \beta - p} \right\}.
	\]
	Using the inequality
	\begin{equation}\label{log_bnd}
		1+ \sum_{j=2}^{n-1} \tau_j (\hat{b}(t_n - t_j))^{-1} \le 1+\hat{b}^{-1}+\hat{b}^{-1}  \int_0^{t_{n-1}} (t_n - s)^{-1} \, ds \le C  {\log\left(2+\frac{t_n}{\tau_n}\right)},
	\end{equation}
	we have
	\[
	\sum_{j=2}^n P_j \le C {\log\left(2+\frac{t_n}{\tau_n}\right)} \max_{2 \le j \le n} \left\{ \tau_j^{q+1+\alpha} t_j^{\beta - q - 1} + \tau_j^p t_j^{\alpha + \beta - p} \right\}.
	\]
	Furthermore, application of \eqref{ineq:stepsize} gives,
	{%
		for any \(\nu_1,\nu_2\in\mathbb{R}\) and \(m\ge 0\),
		\[
		\tau_j^{\nu_1+m}t_j^{\nu_2-m} \le \tau_j^{\nu_1}t_j^{\nu_2}.
		\]
		Setting 
		\[
		\nu_1 = \min\{p,q+1+\alpha\},\quad \nu_2 = \alpha+\beta-\min\{p,q+1+\alpha\}
		\]
		in the above inequality, and choosing
		\[
		m = q+1+\alpha-\min\{p,q+1+\alpha\},
		\]
		for the term \(\tau_j^{q+1+\alpha}t_j^{\beta-q-1}\), and
		\[
		m = p-\min\{p,q+1+\alpha\},
		\]
		for \(\tau_j^p t_j^{\alpha+\beta-p}\), respectively,} we obtain
	\begin{equation}\label{max_meshbnd}
		\tau_j^{q+1+\alpha}t_j^{\beta-q-1} + \tau_j^p t_j^{\alpha+\beta-p} \le 2 \tau_j^{\min\{p,q+1+\alpha\}} t_j^{\alpha+\beta-\min\{p,q+1+\alpha\}}, \quad \text{for all } 2\leq j\leq n.
	\end{equation}
	Thus we can get \eqref{Pj_bnd}, which completes the proof.
\end{proof}

\begin{lemma} \label{lem:Err_2}
	Let the right-hand side term in \eqref{ode} be \( f(t) = t^\beta \vv \) with \(\beta > -1\). Under Assumptions~\ref{assumptionK} and \ref{ass_RK}, for the error term \(\mathcal{E}_n^2\) defined in \eqref{Err_2} on the general time mesh \eqref{mesh}, it holds that
	\begin{enumerate}
		\item[(i)] If \(\beta \ge p - \alpha\), then
		\begin{equation}\label{Err_2_bnd_ge_p_minus_alpha}
			\left\Vert \mathcal{E}_n^2 \right\Vert_{\mathcal{Y}} 
			\le CM \, {\widehat{C}(t_n,\alpha,\beta)}\left(1+\left\vert\log\left(\tau_{\min}\right)\right\vert\right)
			\tau_{\max}^{\min\{p,q+1+\alpha\}}.
		\end{equation}
		
		\item[(ii)] If \( {-1} < \beta < p - \alpha\), then
		\begin{multline} \label{Err_2_bnd_lt_p_minus_alpha}
			\left\Vert \mathcal{E}_n^2 \right\Vert_{\mathcal{Y}} 
			\le CM \Bigg[ \sum_{j=2}^n \tau_j^{p+1} \xi_j^{\alpha + \beta - p - 1} \\
			+  {\log\!\left(2+\frac{t_n}{\tau_{\min}}\right) }
			\max_{2 \leq j \leq n}\left\{ \tau_{j}^{\min\{p,q+1+\alpha\}} t_j^{\alpha+\beta-\min\{p,q+1+\alpha\}} \right\} \Bigg].
		\end{multline}
	\end{enumerate}
	Here, $n\ge 2$,
	\(\xi_j \in (t_{j-1}, t_j)\). The positive constant \( C \) is independent of the mesh \eqref{mesh} and final time $T$.  { $ \widehat{C}(\cdot,\cdot,\cdot)$ is defined by
		\begin{equation}\label{Const_tnbet}
			\widehat{C}(t_n,\alpha,\beta):=\left(1+\left\vert\log\left(t_{n}\right)\right\vert\right) \max\left\{1,t_n^{\alpha+\beta-\min\{p, q + 1 + \alpha\}}\right\}.
	\end{equation}}
\end{lemma}
\begin{proof}
	From the definition of \( \mathcal{E}_n^2 \) in \eqref{Err_2}, together with \eqref{Gx} and \eqref{Rz_bnd}, it follows that
	\begin{equation}
		\label{Err2_split}
		\begin{split}
			\left\| \mathcal{E}_n^2 \right\|_{\mathcal{Y}} &\le CM \sum_{j=2}^n \left( \int_0^\infty x^{-\alpha} \| d_j(x) \|_{\mathcal{X}} \left|\prod_{l=j+1}^{n} R(-\tau_l x)\right| \, dx + P_j \right)\\
			&\le CM \sum_{j=2}^n \left( \int_0^\infty x^{-\alpha} \| d_j(x) \|_{\mathcal{X}} \, dx + P_j \right),
		\end{split}
	\end{equation}
		where \( d_j \) is the defect in \eqref{dn} and $P_j$ is defined in \eqref{Pj}. Applying Lemmas~\ref{lem:integral_bound} and \ref{lem:Pj_bnd},  {together with 
		\begin{equation}
			\label{logtntau_bnd}
			\begin{split}
				\log\left(2+\frac{t_n}{\tau_{\min}}\right) \le  \log\left(3\frac{t_n}{\tau_{\min}}\right)&\le \left(1+|\log (t_n)|\right)+\left(1+|\log(\tau_{\min})|\right)\\[.5em]
				& \le  2\left(1+|\log (t_n)|\right)\left(1+|\log(\tau_{\min})|\right), 
			\end{split}
		\end{equation}
		and \eqref{max_meshbnd}}
	to \eqref{Err2_split} yields  \eqref{Err_2_bnd_ge_p_minus_alpha} and \eqref{Err_2_bnd_lt_p_minus_alpha}, completing the proof.
\end{proof}

Combining Lemmas~\ref{lem:Err_1} and \ref{lem:Err_2}, the convergence result for $f(t)=t^{\beta}\vv$ is stated below.
\begin{proposition}\label{prop:gCQ_conv}
	Assume that \( K(z) \) satisfies Assumption~\ref{assumptionK} and the Runge--Kutta method satisfies Assumption~\ref{ass_RK}. Given the general time mesh \eqref{mesh}, the error of the Runge--Kutta based gCQ  approximation \eqref{gCQ} applied to \eqref{conv} with \( f(t) = t^\beta \vv \) \(( \beta > -1 )\) satisfies
	\begin{enumerate}
		\item[(i)] If \(\beta \ge p - \alpha\), then
		\begin{equation}\label{err_tbeta_ge_p_minus_alpha}
			\left\Vert u(t_n)-u_n \right\Vert_{\mathcal{Y}} 
			\le CM \, {\widehat{C}(t_n,\alpha,\beta)}\left(1+\left\vert\log\left(\tau_{\min}\right)\right\vert\right)
			\tau_{\max}^{\min\{p,q+1+\alpha\}}.
		\end{equation}
		
		\item[(ii)] If \( {-1} < \beta < p - \alpha\), then
		\begin{multline} \label{err_tbeta_lt_p_minus_alpha}
			\left\Vert u(t_n)-u_n \right\Vert_{\mathcal{Y}} 
			\le CM \Bigg[ \sum_{j=2}^n \tau_j^{p+1} \xi_j^{\alpha + \beta - p - 1} \\
			+  {\log\!\left(2+\frac{t_n}{\tau_{\min}}\right) }
			\max_{1 \leq j \leq n}\left\{ \tau_{j}^{\min\{p,q+1+\alpha\}} t_j^{\alpha+\beta-\min\{p,q+1+\alpha\}} \right\} \Bigg].
		\end{multline}
	\end{enumerate}
	where  $n\ge 1$, \( \xi_j \in (t_{j-1}, t_j) \),  $C$ is a positive constant independent of the time mesh \eqref{mesh}  {and \( \widehat{C}(t_n,\alpha,\beta)\)  is given by \eqref{Const_tnbet}}.
\end{proposition}
Below we derive specific error estimates for two  {special} cases: the case of  {locally} quasi-uniform meshes like in \eqref{quasimesh} and, among them, the case of graded meshes \eqref{gmesh}. For general  {locally} quasi-uniform meshes we can relax the requirement about the number of vanishing moments at the origin, whereas for graded meshes we do not need to assume any of them, since we can always achieve full order of convergence by choosing properly the grading parameter.

\begin{corollary}\label{cor:conv_tbeta}
	Under the assumptions in Proposition~\ref{prop:gCQ_conv}, the following error estimates hold for \( f(t) = t^{\beta} \vv \).
	\begin{enumerate}[label=(\roman*)]
		\item  \label{itm:tbet_quasi} If the mesh is  {locally} quasi-uniform as in \eqref{quasimesh}, we can relax the condition on $\beta$ in \eqref{err_tbeta_ge_p_minus_alpha} to  \( \beta \ge \min\{p-\alpha, q+1 \} \), and the error is bounded by
		\begin{equation}\label{fullrate_quasimesh}
			\left \Vert u(t_n) - u_n \right\Vert_{\mathcal{Y}} \le C M c_\Delta^{\max\{0,\, p-q- \alpha\}} { \widehat{C}(t_n,\alpha,\beta)\left(1+\left\vert\log\left(\tau_{\min}\right)\right\vert\right)}  \tau_{\max}^{\min\{p,q+1+\alpha \}} , \, n \ge 1,
		\end{equation}
		{with $ \widehat{C}(\cdot,\cdot,\cdot)$ given in \eqref{Const_tnbet}.}
		\item  \label{itm:tbet_grad} For any $\beta >-1$, the error on the graded mesh \eqref{gmesh} can be estimated by
		\begin{equation*}
			\left \Vert u(t_n) - u_n \right\Vert_{\mathcal{Y}} \!\leq  \!\! CM  {T^{ \alpha+\beta}} \log(N)\!\!
			\begin{cases}
				\!	N^{-\gamma(\alpha + \beta)}, & \gamma(\alpha + \beta) \!< \min\{p, q + 1 + \alpha\}, \\[.5em]
				\!	N^{-\min\{p, q + 1 + \alpha\}}, & \gamma(\alpha + \beta) \!\ge \min\{p, q + 1 + \alpha\},
			\end{cases}
		\end{equation*}
		where \( 1 \le n \le N \), \( C \) is independent of \( N \).
	\end{enumerate}
\end{corollary}
\begin{proof}
	When \( \beta \ge p - \alpha \), the result in Cases~\ref{itm:tbet_quasi}  follows directly from Proposition~\ref{prop:gCQ_conv}.

	For Case~\ref{itm:tbet_quasi}, it remains to derive the error bound under the condition
	\begin{equation}\label{bet_cond}
		\beta \ge q+1 \quad \text{with} \quad q+1 = \min\{p - \alpha, q + 1\}.
	\end{equation}
	We begin by analyzing the summation in \eqref{err_tbeta_lt_p_minus_alpha}. The  {locally} quasi-uniformity property of the mesh \eqref{quasimesh} implies the lower bound
	\[
	t_{j-1} = \sum_{\ell = 1}^{j-1} \tau_\ell > \sum_{\ell = 1}^{j-1} \frac{\tau_\ell}{2} + \frac{\tau_j}{4c_\Delta} \ge \frac{t_j}{4c_\Delta}.
	\]
	Applying this bound to the summation in \eqref{err_tbeta_lt_p_minus_alpha} yields
	\begin{align*}
		\sum_{j=2}^n \tau_j^{p+1} \xi_j^{\alpha + \beta - p - 1}
		&\le (4c_\Delta)^{\max\{0,\, p+1-\alpha-\beta\}} \sum_{j=2}^n \tau_j^{p+1} t_j^{\alpha + \beta - p - 1} \\
		&\le (4c_\Delta)^{\max\{0,\, p+1-\alpha-\beta\}} \left( \sum_{j=2}^n \tau_j t_j^{-1} \right) \max_{2 \le j \le n} \left\{ \tau_j^p t_j^{\alpha + \beta - p} \right\}.
	\end{align*}
	Using the integral estimate
	\[
	\sum_{j=2}^n \tau_j t_j^{-1} < 1 + \int_{t_1}^{t_n} t^{-1} \, dt = 1 + \log(t_n / t_1)
	\]
	along with \eqref{max_meshbnd},  {\eqref{logtntau_bnd}} and \eqref{bet_cond}, we obtain
	\begin{equation}\label{sum_quasimeshbnd}
		\sum_{j=2}^n \tau_j^{p+1} \xi_j^{\alpha + \beta - p - 1} 
		\le C\, c_\Delta^{\max\{0,\, p-q - \alpha\}}  {\left(1+\left\vert\log\left(t_n\right)\right\vert\right)\left(1+\left\vert\log\left(\tau_{\min}\right)\right\vert\right)}\tau_{\max}^{q + 1 + \alpha}\, t_n^{\beta - q - 1}.
	\end{equation}
	Substituting the bound above into \eqref{err_tbeta_lt_p_minus_alpha} and then applying  {\eqref{logtntau_bnd} }completes the proof of Case~\ref{itm:tbet_quasi}.

	For case~\ref{itm:tbet_grad}, the bound in Proposition~\ref{prop:gCQ_conv} on the graded mesh \eqref{gmesh} satisfies
	\begin{align*}
		&\left\| u(t_n) - u_n \right\|_{\mathcal{Y}} 
		\leq \\[.5em]       
		&\qquad C M \left( \sum_{j=2}^n \tau_j^{p+1} \xi_j^{\alpha+\beta - p - 1} + \log(N) \max_{1 \leq j \leq n} 
		\left\{ \tau_{j}^{\min\{p,q+1+\alpha\}} t_j^{\alpha+\beta-\min\{p,q+1+\alpha\}}\right\} \right).
	\end{align*}
	Applying \eqref{grad_sumprodbnd} gives
	\begin{align*}
		\sum_{j=2}^n \tau_j^{p+1} \xi_j^{\alpha+\beta - p - 1} 
		\leq C {T^{ \alpha+\beta}}
		\begin{cases}
			N^{-\gamma(\alpha+\beta)}, & \gamma(\alpha+\beta) < p, \\[0.5em]
			N^{-p} \log(n), & \gamma(\alpha+\beta) = p, \\[0.5em]
			N^{-p}, & \gamma(\alpha+\beta) > p.
		\end{cases}
	\end{align*}
	Next, using \eqref{grad_prodbnd}, we  can derive
	\begin{align*}
		\nonumber
		& \tau_{j}^{\min\{p,q+1+\alpha\}} t_j^{\alpha+\beta-\min\{p,q+1+\alpha\}} \\[.5em]
		&\le C   {T^{ \alpha+\beta}}
		\begin{cases}
			N^{-\gamma(\alpha + \beta)}, & \text{if } \gamma(\alpha + \beta) < \min\{p,\, q+1+\alpha\}, \\[0.5em]
			N^{-\min\{p,\, q+1+\alpha\}}, & \text{if } \gamma(\alpha + \beta) \ge \min\{p,\, q+1+\alpha\}.
		\end{cases}
	\end{align*}
	Combining these estimates proves \ref{itm:tbet_grad}.
\end{proof}

\subsection{Convergence analysis for general integrands}
\begin{lemma}
	
	\label{lem:int_dj_peano}
	Let \( y(x, t) \) be the solution to \eqref{ode} with right-hand side given by
	\begin{equation}\label{f_heavis}
		f(t) = H(t - \eta)(t - \eta)^{\beta} \mathbf{v}, \quad \eta \geq 0,
	\end{equation}
	where \( H(\cdot) \) is the Heaviside step function and \( \mathbf{v} \in \mathcal{X} \). Let \( d_j(x) \) denote the associated defect in \eqref{RK_ode_defec}. Under Assumptions~\ref{assumptionK} with \( \beta \geq p - 1 \geq 0 \), for \( 2 \leq j \leq n-1 \), it holds that
	\begin{equation*}
		\| d_j(x) \|_{\mathcal{X}} \left| \prod_{l=j+1}^{n} R(-\tau_l x) \right| \leq C\tau_j^{p} { t_j^{\beta-p+1}}
		\begin{cases}
			\displaystyle\int_{t_{j-1}}^{t_j} \frac{1} {(t - \eta)(1 + \hat{b}(t_n - t_j) x)} \, dt, & \eta \in [0, t_{j-1}),
			\\[0.7em]
			(1 + \hat{b}(t_n - t_j) x)^{-1},    & \eta \in [t_{j-1}, t_j], \\[0.5em]
			0, & \eta \in (t_j, \infty).
		\end{cases}
	\end{equation*}
	Moreover, for \( j = n \), we have
	\begin{equation*}
		\| d_n(x) \|_{\mathcal{X}}   \leq C     \tau_n^{p}  { t_n^{\beta-p+1}}
		\begin{cases}
			\displaystyle\int_{t_{n-1}}^{t_n} \frac{1}{(t - \eta)(1 + (t - \eta) x)} \, dt, & \eta \in [0, t_{n-1}), \\[0.7em]
			\left\Vert H(t_{n-1}+\cv \tau_n - \eta)(1 + (t_{n-1}+\cv \tau_n - \eta) x)^{-1}\right\Vert_\infty,    & \eta \in [t_{n-1}, t_n], \\[0.5em]
			0, & \eta \in (t_n, \infty), 
		\end{cases}
	\end{equation*}
	where  $C$ and $\hat{b}$ are positive constants independent of the mesh. 
\end{lemma}
\begin{proof}
	From the representation in \eqref{RK_ode_defec}, the defect $d_j(x)$ can be expressed as
	\begin{equation*}\label{dj_peano}
		d_j(x) = 
		\begin{cases}    \tau_j^p \int_{t_{j-1}}^{t_j} \kappa\left(\dfrac{t - t_{j-1}}{\tau_j}\right)
			\dfrac{\partial^{p+1} y}{\partial t^{p+1}}(x,t) dt, 
			& \eta \in [0, t_{j-1}), \\[0.5em]
			y(x,t_j) - y(x,t_{j-1}) - \tau_j \mathbf{b}^\top \mathbf{Y}_t(x, t_j), 
			& \eta \in [t_{j-1}, t_j] , \\[0.5em]
			0, & \eta \in(t_j,\infty).
		\end{cases} 
	\end{equation*}
	Applying \eqref{Rz_bnd} and \eqref{ydiff_bndpeano}, we obtain for \( \eta \in [0, t_{j-1}) \):
	\begin{equation*}\label{dj_nearbnd}
		\|d_j(x)\|_{\mathcal{X}} \left|\prod_{l=j+1}^n R(-\tau_l x)\right|
		\leq C \tau_j^p\int_{t_{j-1}}^{t_j} \dfrac{(t - \eta)^{\beta-p}}{(1+(t-\eta)x)(1 + \hat{b}(t_n - t_j)x)} dt,   
	\end{equation*}
	and for \( \eta \in [t_{j-1}, t_j] \):
	\begin{equation*}\label{dj_midbnd}
		\|d_j(x)\|_{\mathcal{X}} \left|\prod_{l=j+1}^n R(-\tau_l x)\right|
		\leq C \tau_j^{\beta+1}\left\Vert  \dfrac{H(t_{j-1}+\cv \tau_j- \eta)}{(1+(t_{j-1}+\cv \tau_j - \eta)x)(1 + \hat{b}(t_n - t_j)x)}\right\Vert_\infty, 
	\end{equation*}
		where \( \hat{b} > 0 \) is the constant from \eqref{Rz_bnd}. The final bounds then follow by using $\beta\ge p-1$ and
	\[
	(t-\eta)^{\beta-p} \leq  t_j^{\beta-p+1}/(t-\eta), \quad \text{for }  \, t \in (t_{j-1}, t_j), \ \eta \in (0, t_{j-1}), \ 2 \leq j \leq n.
	\]
\end{proof}

\begin{lemma}\label{lem:E1Pj_peano}
	Let \( y(x, t) \) be the solution to \eqref{ode} with \( f \) given by \eqref{f_heavis}, where \( \beta \ge p-1 \ge 0 \) and \( 0 < \alpha < 1 \). Under Assumptions~\ref{assumptionK} and \ref{ass_RK}, the following estimates hold:
	\begin{equation*}
		\left\| \mathcal{E}_n^1 \right\|_{\mathcal{Y}} 
		\leq C M \tau_1 H(\tau_1 - \eta)(\tau_1 - \eta)^{\alpha + \beta - 1},
	\end{equation*}
	and
	\begin{align*}
		\sum_{j=2}^n P_j 
		&\leq  \!C {\log\left(\,2+\frac{t_n}{\tau_n}\,\right)} \!\max_{2 \le j \le n}\! \left\{\! H(t_j - \eta)\! {\left(\!\sum_{\ell=q+1}^{p-1} \tau_j^{\alpha+\ell}(t_j-\eta)^{\beta-\ell}+\tau_j^p(t_j - \eta)^{\alpha+\beta-p}\!\right)}\!\!\right\},
	\end{align*}
	where \( \mathcal{E}_n^1 \) denotes the approximation error defined in \eqref{Err_1}, and \( P_j \) is given by \eqref{Pj}.
\end{lemma}
\begin{proof}
	Following the proof of Lemma~\ref{lem:Err_1}, we obtain
	{\begin{align*}
			\left\Vert \mathcal{E}_n^1 \right\Vert_{\mathcal{Y}} 
			&\le C M \int_{0}^\infty x^{-\alpha} \left(\left\| y(x, t_1)\right\|_{\mathcal{X}}+\left\|\tau_1 \bv^\top(\Iv + \tau_1 x \Av)^{-1}\fv_1\right\|_{\mathcal{X}} \right) \, dx.
		\end{align*}
		Applying \eqref{ydiff_bndpeano} with \(m=0\) to \(y(x, t_1)\), together with the expression for \(\fv_1\), and proceeding as in the proof of Lemma~\ref{lem:Err_1} yields the stated bound.
	}

	For \( P_j \) in \eqref{Pj},    adapting the proofs of \eqref{Pn_esti} and \eqref{Pj_esti} and applying \eqref{ydiff_bndpeano} yields
	\begin{align*}  
		P_n &\leq C H(t_n-\eta) \left( \sum_{\ell=q+1}^{p-1} \tau_n^{\alpha+\ell}(t_n-\eta)^{\beta-\ell} + \tau_n^p(t_n-\eta)^{\alpha+\beta-p} \right), \\[0.5em]
		P_j &\leq C H(t_j-\eta) \sum_{\ell=q+1}^{p-1} \tau_j^{\alpha+\ell+1}(\hat{b}(t_n-t_j))^{-1}(t_j-\eta)^{\beta-\ell} \\[0.5em]
		&\quad + C H(t_j-\eta)\tau_j^{p+1}(\hat{b}(t_n-t_j))^{-1}(t_j-\eta)^{\alpha+\beta-p}, \quad 2 \leq j \leq n-1.
	\end{align*}
	{Applying  \eqref{log_bnd} yields the final bound.}
\end{proof}

\begin{proposition}\label{propo:err_remainder}
	Let the data $f$ in \eqref{conv} be given by
	\begin{equation}\label{int_remainder}
		f(t) = \frac{1}{\Gamma(m+\beta+1)} \left(t^{m+\beta} \ast f^{(m+\beta+1)}\right)(t),
	\end{equation}
	where \(\beta > -1\) and \(m = \lceil p - \beta - 1\rceil\). Under Assumptions~\ref{assumptionK} and~\ref{ass_RK}, the error of the Runge-Kutta based gCQ  approximation \eqref{gCQ} on the general mesh \eqref{mesh} satisfies
	\begin{equation*}
		\begin{split}
			&\left\| u(t_n) - u_n \right\|_{\mathcal{Y}} \leq \\[.5em]
			&\hfill\qquad\qquad\quad   CM  { \widehat{C}(t_n,\alpha,m+\beta+1)}  {\left(1+\left\vert\log\left(\tau_{\min}\right)\right\vert\right)}\tau_{\max}^{\min\{p,q+1+\alpha\}} \left\| f^{(m+\beta+1)} \right\|_{C^0([0,t_n],\mathcal{X})},    
		\end{split}
	\end{equation*}
	where $1 \leq n \leq N$, \(C\) is a constant independent of the temporal discretization \eqref{mesh} and the final time $T$,  { \( \widehat{C}(\cdot,\cdot,\cdot)\) is defined in \eqref{Const_tnbet}
		.}
\end{proposition}
\begin{proof}
	Denote \((t-\eta)_+^{m+\beta}:=H(t-\eta) (t-\eta)^{m+\beta}.\)
	Following the argument in \cite[Proposition~13]{GuoLo}, the continuous convolution and its Runge--Kutta based gCQ  approximation, when applied to \( f(t) \) in \eqref{int_remainder}, can be represented as
	\[
	\left[K(\partial_t)f\right](t_n) = \int_0^{t_n} \frac{f^{(m+\beta+1)}(\eta)}{\Gamma(m+\beta+1)} \left[ K(\partial_t) (\cdot - \eta)^{m+\beta}_+\right](t_n) \, d\eta,
	\]
	and
	\[
	\left[K(\partial_t^\Delta) f\right]_n = \int_0^{t_n} \frac{f^{(m+\beta+1)}(\eta)}{\Gamma(m+\beta+1)} \left[ K(\partial_t^\Delta) (\cdot - \eta)_+^{m+\beta} \right]_n\, d\eta.
	\]
	The approximation error then satisfies
	\begin{align*}
		&\left\| \left[K(\partial_t)f\right](t_n) - \left[K(\partial_t^\Delta) f\right]_n\right\|_{\mathcal{Y}} \\
		&= \left\| \int_{0}^{t_n} \frac{f^{(m+\beta+1)}(\eta)}{\Gamma(m+\beta+1)} \left( \left[K(\partial_t)(\cdot-\eta)_+^{m+\beta} \right] (t_n) - \left[ K(\partial_t^\Delta)(\cdot-\eta)_+^{m+\beta} \right]_n\right) d\eta \right\|_{\mathcal{Y}} \\
		&\leq \int_0^{t_n} \left\|\left[K(\partial_t)(\cdot-\eta)_+^{m+\beta} \right] (t_n) - \left[ K(\partial_t^\Delta)(\cdot-\eta)_+^{m+\beta} \right]_n \right\|_{\mathcal{Y}} d\eta \cdot \left\| f^{(m+\beta+1)} \right\|_{C^0([0,t_n],\mathcal{X})}.
	\end{align*}
	Let \( y(x,t) \) be the solution to \eqref{ode} with right hand side \( (t-\eta)_+^{m+\beta} \), and let \( e_n(x) \) denote the associated Runge-Kutta approximation error from \eqref{en}. Combining the error representation in \eqref{gCQ_error_decomposition} with the bound from \eqref{Err2_split}, we obtain 
	\begin{align}\nonumber
		\int_0^{t_n} &\left\|\left[K(\partial_t)(\cdot-\eta)_+^{m+\beta} \right] (t_n) - \left[ K(\partial_t^\Delta)(\cdot-\eta)_+^{m+\beta} \right]_n \right\|_{\mathcal{Y}}\, d\eta \\
		\label{djsum_nminus_split}
		&\leq \int_0^{t_n}\left[ \left\| \mathcal{E}_n^1 \right\|_{\mathcal{Y}}  + CM \sum_{j=2}^n \left( \int_0^\infty x^{-\alpha} \| d_j(x) \|_{\mathcal{X}}\left|\prod_{l=j+1}^{n} R(-\tau_l x)\right| dx + P_j \right)\right]\, d\eta.
	\end{align}
	For \( 2 \le j \le n - 1 \), applying Lemma~\ref{lem:int_dj_peano} under the condition \( m + \beta \ge p - 1 \), we obtain
	\begin{align*}
		&\int_0^{t_n}\left( \int_0^\infty x^{-\alpha} \| d_j(x) \|_{\mathcal{X}} \left|\prod_{l=j+1}^{n} R(-\tau_l x)\right| \,dx \right) \, d\eta \\
		\nonumber
		&\le C  \tau_j^{p}  {t_j^{m+\beta-p+1}}\int_0^\infty \frac{x^{-\alpha}}{1 + \hat{b}(t_n - t_j)x} \, dx 
		\left(  \int_0^{t_{j-1}} \int_{t_{j-1}}^{t_j} \frac{1}{t - \eta} \, dt \, d\eta+\tau_j  \right)\\
		&\overset{\eqref{int_Betafun}}{=} C \tau_j^{p}  {t_j^{m+\beta-p+1}}(\hat{b}(t_n - t_j))^{\alpha - 1}
		\left( \int_0^{t_{j-1}} \int_{t_{j-1}}^{t_j} \frac{1}{t - \eta} \, dt \, d\eta+\tau_j  \right).
	\end{align*}
	For the integral term, we have
	\begin{align*}
		\int_0^{t_{j-1}} \int_{t_{j-1}}^{t_j} \frac{1}{t - \eta} \, dt \, d\eta & =
		t_j \log(t_j) - \tau_j \log(\tau_j) - t_{j-1} \log(t_{j-1})\\
		&\le  {\tau_j \log(t_j)-\tau_j \log(\tau_j) = \log\left(\frac{t_j}{\tau_{j}}\right)}\tau_j.
	\end{align*}
	Thus, we can bound the summation w.r.t. $j$ in \eqref{djsum_nminus_split} by
	\begin{equation}\label{djsum_nmius1_bnd}
		\begin{split}
			&\sum\limits_{j=2}^{n-1}\int_0^{t_n}\left( \int_0^\infty x^{-\alpha} \| d_j(x) \|_{\mathcal{X}} \left|\prod_{l=j+1}^{n} R(-\tau_l x)\right| \,dx \right) \, d\eta \\
			&\le C  {\log\left(\frac{t_n}{\tau_{\min}}\right)}\sum\limits_{j=2}^{n-1} \tau_j^{p+1}   {t_j^{m+\beta-p+1}}(\hat{b}(t_n - t_j))^{\alpha - 1}\\
			&\le C  {\log\left(\frac{t_n}{\tau_{\min}}\right)t_n^{\alpha}\max\limits_{2\le j\le n-1}\{\tau_{j}^p t_j^{m+\beta-p+1}\}}
			\\[.5em]
			& {\le C\log\left(\frac{t_n}{\tau_{\min}}\right)\tau_{\max}^{\min\{p,q+1+\alpha\}}\max\left\{1,t_n^{\alpha+m+\beta+1-\min\{p,q+1+\alpha\}}\right\}. }
		\end{split}
	\end{equation}
	Now, for the term corresponding to \( j = n \), applying Lemma~\ref{lem:int_dj_peano} and \eqref{int_Betafun}, we get
	\begin{equation}\label{dj_n_bnd}
		\begin{split}
			& \int_0^{t_n} \int_0^\infty x^{-\alpha} \| d_n(x) \|_{\mathcal{X}} \, dx d\eta 
			\le \tau_n^p {t_n^{m+\beta-p+1}} \int_0^{t_{n-1}} \int_{t_{n-1}}^{t_n} (t - \eta)^{\alpha - 2} \, dt \, d\eta \\[.5em]
			&+ C\tau_n^{p}  {t_n^{m+\beta-p+1}}\int_{t_{n-1}}^{t_n}  \left\|H(t_{n-1} + \cv \tau_n - \eta) (t_{n-1} + \cv \tau_n - \eta)^{\alpha - 1} \right\|_\infty \, d\eta \\[.5em]
			&\le C \tau_n^p  {t_n^{m+\beta-p+1}}\left(\int_0^{t_{n-1}} \left[ (t_n - \eta)^{\alpha - 1} - (t_{n-1} - \eta)^{\alpha - 1} \right] d\eta + \tau_n^{\alpha}\right)\\[.5em]
			&\le C\tau_n^p { t_n^{\alpha+m+\beta-p+1}\le C\tau_n^{\min\{p,q+1+\alpha\}}t_n^{\alpha+m+\beta+1-\min\{p,q+1+\alpha\}}.}     
		\end{split}
	\end{equation}
	From \eqref{djsum_nmius1_bnd} and \eqref{dj_n_bnd}, we  {obtain} that
	\begin{align*}
		& \sum\limits_{j=2}^n\int_0^{t_n}\left( \int_0^\infty x^{-\alpha} \| d_j(x) \|_{\mathcal{X}} \left|\prod_{l=j+1}^{n} R(-\tau_l x)\right| \,dx \right) \, d\eta\\
		&\le C  {\log\left(2+\frac{t_n}{\tau_{\min}}\right)\tau_{\max}^{\min\{p,q+1+\alpha\}}\max\left\{1,t_n^{\alpha+m+\beta+1-\min\{p,q+1+\alpha\}}\right\}.}  
	\end{align*}
	Applying Lemma~\ref{lem:E1Pj_peano} together with the bound above to \eqref{djsum_nminus_split}  {and using \eqref{logtntau_bnd}} yields
	\begin{equation*}
		\begin{split}
			&\left\| u(t_n) - u_n \right\|_{\mathcal{Y}} \leq\\[.5em]
			&\hfill\qquad\qquad\quad  CM  { \widehat{C}(t_n,\alpha,m+\beta+1)}  {\left(1+\left\vert\log\left(\tau_{\min}\right)\right\vert\right)}\tau_{\max}^{\min\{p,q+1+\alpha\}} \left\| f^{(m+\beta+1)} \right\|_{C^0([0,t_n],\mathcal{X})},    
		\end{split}
	\end{equation*}
	with  {$ \widehat{C}(\cdot,\cdot,\cdot)$ given in \eqref{Const_tnbet}.}

\end{proof}
\begin{remark}\label{remark:improvegCQ1}
	By applying the  bound \eqref{int_remainder} with \( p = 1 \) to the analysis of the Euler based gCQ  method, we can relax the regularity condition in \cite[Proposition 13]{GuoLo} 
	to data \( f \) in $C^1([0,T],\mathcal{X})$, rather than in $C^3([0,T],\mathcal{X})$. 
\end{remark}
By applying Proposition~\ref{prop:gCQ_conv} and Proposition~\ref{propo:err_remainder}, we can establish the following convergence results for general integrands \(f(t)\) of the form \eqref{gen_f}.

\begin{proposition}\label{prop:conv_gCQrk}
	Suppose that the transfer operator \(K(z)\) satisfies Assumption~\ref{assumptionK} and the Runge--Kutta method satisfies Assumption~\ref{ass_RK}. Let the integrand \(f(t)\) in  the convolution  \eqref{conv} take the form of  \eqref{gen_f} with \(\beta > - 1\) and \(m = \lceil p - \beta-1\rceil\). Denote 
	\begin{equation*}
		\mu:=\min\{p,q + 1+\alpha\} \text{ and } \nu:=\lceil p-\alpha-\beta - 1\rceil. 
	\end{equation*}
	Then, on the general mesh \eqref{mesh}, it holds
	\begin{align*}
		&\|u(t_n)-u_n\|_{\mathcal Y}
		\le
		\\
		&CM\left[
		\sum_{\ell=0}^{\nu}
		\left(
		\sum_{j=2}^{n}
		\tau_j^{p+1}\xi_j^{\alpha+\beta+\ell-p-1}
		+  {\log\!\left(2+\frac{t_n}{\tau_{\min}}\right)}
		\max_{1\le j\le n}
		\left\{
		\tau_j^\mu t_j^{\alpha+\beta+\ell-\mu}
		\right\}
		\right)
		\|f^{(\ell+\beta)}(0)\|_{\mathcal X}
		\right.
		\\
		&\!\left.
		+  { \widehat{C}(t_n,\alpha,m+\beta+1)\!
			\left(1+|\log (\tau_{\min})|\right)}
		\tau_{\max}^{\mu}\!
		\left(\!
		\sum_{\ell=\nu+1}^{m}
		\|f^{(\ell+\beta)}(0)\|_{\mathcal X}
		+
		\|f^{(m+\beta+1)}\|_{C^0([0,t_n],\mathcal X)}
		\!\right)\!
		\right].
	\end{align*}
	Here, \(n\ge 1\), \(\xi_j\in(t_{j-1},t_j)\), the constant \(C\) is independent of the mesh and the final time \(T\), and \( \widehat{C}(\cdot,\cdot,\cdot)\) is given in \eqref{Const_tnbet}.
\end{proposition}
\begin{remark}
	Proposition~\ref{prop:conv_gCQrk} directly implies the convergence result for general meshes, as stated in Case~\ref{itm:thm_genmesh} of Theorem~\ref{thm:conv_gCQrkGenM}. For  {locally} quasi-uniform meshes, the convergence result follows from Case~\ref{itm:tbet_quasi} of Corollary~\ref{cor:conv_tbeta}, which is summarized as \ref{itm:thm_quasimesh} in Theorem~\ref{thm:conv_gCQrkGenM}. In the case of graded meshes, the corresponding convergence result is obtained via Case~\ref{itm:tbet_grad} of Corollary~\ref{cor:conv_tbeta}, and is formally stated in Theorem~\ref{thm:conv_gCQrkGrad}.
\end{remark}

\section{Fast and oblivious Runge--Kutta gCQ}\label{sec:algfi}
\subsection{Algorithm for the convolution}\label{sec:algconv}
For efficient implementation, we extend the algorithm proposed in \cite{BanLo19} to the general mesh \eqref{mesh}. Denoting
\[\fv_n=\left(f(t_{n-1}+c_i\tau_n)\right)_{i=1}^s,\quad 1\le n\le N.\]
Following  \cite{GuoLo} and recalling the notation in \eqref{Un_stage}, we split the gCQ  approximation into two parts:
\begin{equation}\label{gCQshort_split}
	\Uv_n = \left[K\left( \partial^{\Delta}_t \right)\fv \right]_n  =\sum\limits_{j=\max\left(1, n-n_0+1\right)}^{n}\Wv_{n,j} \fv_j+ \sum\limits_{j=1}^{n-n_0}\Wv_{n,j} \fv_j = \Sv_n^{loc}+\Sv_n^{\mathrm{his}},
\end{equation}
Here, $n_0$ is a fixed moderate constant; for example, we use $n_0=5$, although $n_0=1$ is also admissible. Using the definition of $\Wv_{n,j}$ in \eqref{gCQ_w}, we have
\begin{align}
	\nonumber
	\Sv_n^{\mathrm{loc}}&=\frac{1}{2\pi i}\int_{\mathcal{C}}  K(z) \sum\limits_{j=\max\left(1, n-n_0+1\right)}^{n}\left(\prod_{l=j+1}^{n} \left(\Rv(\tau_lz)\ev_s^\top\right)\right) \left  (\tau_j\Av(\Iv-\tau_jz\Av)^{-1}\fv_j \right) \,dz\\
	\label{Snloc}
	&:=\frac{1}{2\pi i}\int_{\mathcal{C}} K(z)\Qv_n^{\mathrm{loc}}(z) \,dz,
\end{align}
where
\begin{equation}\label{recursion_loc}
	\Qv^{\mathrm{loc}}_{n}(z)=\Rv(\tau_nz) \ev_s^\top\Qv^{\mathrm{loc}}_{n-1}(z)+\tau_n\Av(\Iv-\tau_n\Av z)^{-1}\fv_n,\qquad \Qv^{\mathrm{loc}}_{(n-n_0)_+}(z)=\mathbf{0}.
\end{equation}
To approximate the local term, we select $\mathcal{C}$ as a circle in the right half of the complex plane, enclosing the poles of the integrand, and parameterized by using elliptic functions. The same idea appears in \cite{GuoLo} for Euler based gCQ, but we introduce the adjustments in $\mathcal{C}$ indicated for Runge--Kutta methods in \cite{LoSau16}. The radius of $\mathcal{C}$ is then given by
\begin{equation}\label{defMloc}
	M^{\mathrm{loc}}=5\max\left\{\left\vert {\rm{spec}}\left(\Av^{-1}\right)\right\vert\right\}/\min\limits_{\max\left(1, n-n_0+1\right)\le j\le \max(n,n_0)}\{\tau_j\},
\end{equation}
and centered at  {$M^{\mathrm{loc}}+m^{\mathrm{loc}}/10$}, where
\begin{equation}\label{defmloc}
	m^{\mathrm{loc}}=\min\left\{\Re\left  ({\rm{spec}}\left (\Av^{-1}\right)\right )\right\}/\max\limits_{\max\left(1, n-n_0+1\right)\le j\le \max(n,n_0)}\{\tau_j\}.
\end{equation}
{For high-order Runge--Kutta methods, the factor \(5\) in $M^{\mathrm{loc}}$ is chosen to ensure that the contour encloses all poles of \(\Qv_n^{\mathrm{loc}}(z)\) in \eqref{Snloc}, while preventing it from invading too far into the region where \(|R(\tau_{\min} z)| > 1\), as it is done in \cite{LoSau16}. Our contour is further shifted by \(m^{\mathrm{loc}}/10\) to guarantee that \(K(z)\) in \eqref{Snloc} remains bounded along \(\mathcal{C}\). The choice of these precise factors is heuristic and turned out to be stable and to work very well in all our numerical examples.} 
Subsequently, the quadrature rule for the local term can be formulated as follows
\begin{equation}\label{Sloc_quad}
	\Sv^{\mathrm{loc}}_n \approx \sum_{l=1}^{N_{Q}^{\mathrm{loc}}} w_l z_l^{-\alpha} \Qv_n^{\mathrm{loc}}(z_l).
\end{equation}
We notice that the circle $\mathcal{C}$ is fixed for the first $n_0$ steps, until we compute $\Uv_{n_0}$, and then it changes at every time step $n>n_0$.

For the history part in \eqref{gCQshort_split} we notice that
\begin{align}
	\nonumber
	\Sv_n^{\mathrm{his}}&=\int_{0}^\infty G(x) \sum\limits_{j=1}^{n-n_0}\left(\prod_{l=j+1}^{n} \left(\Rv(-\tau_lx)\ev_s^\top\right)\right) \left  (\tau_j\Av(\Iv+\tau_j x\Av)^{-1} \fv_j\right) \,dx\\
	\label{S_hisint}
	&:= \int_{0}^\infty G(x)\left(\prod_{l=n-n_0+1}^{n} \left(\Rv(-\tau_lx)\ev_s^\top\right)\right)\Qv^{\mathrm{his}}_{n-n_0}(x) \,dx,
\end{align}
where
\[\Qv^{\mathrm{his}}_{n-n_0}(x)= \sum\limits_{j=1}^{n-n_0}\left(\prod_{l=j+1}^{n-n_0} \left(\Rv(-\tau_lx)\ev_s^\top\right)\right) \left  (\tau_j\Av(\Iv+\tau_j x\Av)^{-1}\fv_j \right),\]
satisfying
\begin{equation}\label{recursion_his}
	\Qv^{\mathrm{his}}_{n-n_0}(x)=\Rv(-\tau_{n-n_0}x)\ev_s^\top\Qv^{\mathrm{his}}_{n-n_0-1}(x)+\tau_{n-n_0}\Av (\Iv+\tau_{n-n_0}\Av x)^{-1}\fv_{n-n_0},\quad \Qv^{\mathrm{his}}_{0}(x)=\mathbf{0}.
\end{equation}
{As discussed below, in the numerical results reported in Section~\ref{sec:num_test} we apply different quadrature rules according to the form of \( G(x) \).} In any case, this yields an approximation of the form
\begin{align}
	\label{Shis_quad}
	\Sv^{\mathrm{his}}_n \approx \sum\limits_{l=1}^{N^{\mathrm{his}}_Q} \varpi_l G(x_l)\left(\prod_{j=n-n_0+1}^{n} \left(\Rv(-\tau_jx_l)\ev_s^\top\right)\right) \Qv_{n-n_0}^{\mathrm{his}}(x_l).
\end{align}
Thus, we obtain the approximation
\begin{equation}\label{gCQ_appInt}
	\left[K\left( \partial^{\Delta}_t\right)\fv\right]_n  \approx \sum_{l=1}^{N_{Q}^{\mathrm{loc}}} w_l K(z_l) \Qv_n^{\mathrm{loc}}(z_l)+\sum\limits_{l=1}^{N^{\mathrm{his}}_Q} \varpi_lG(x_l)\left(\prod\limits_{j=n-n_0+1}^{n}\left(\Rv(-\tau_jx_l)\ev_s^\top\right)\right)\Qv_{n-n_0}(x_l).
\end{equation}
In our numerical experiments, we set
\begin{equation*}
	n_0=\min\left(5,N\right).
\end{equation*}
Applying the formula from \cite[Corollary 16]{LoSau15apnum}, we choose
\begin{equation}\label{NQloc}
	N^{\mathrm{loc}}_Q = \left\lceil n_0^{\max\left\{1,\, \frac{1}{2} \log(M^{\mathrm{loc}})/\log(m^{\mathrm{loc}}) \right\}} \log(n_0)\left( \log(n_0) + \log(1/\mathrm{tol}) \right) \right\rceil
\end{equation}
in the quadrature rule \eqref{Sloc_quad}, where  {\(0<\mathrm{tol}<1\)} denotes the prescribed tolerance for the quadrature used in both \eqref{Sloc_quad} and \eqref{Shis_quad}.  {Notice that, for $n>n_0$, $N^{\mathrm{loc}}_Q$ does not depend on $N$, $n$, $\tau_{\min}$ nor $\tau_{\max}$, but only on $n_0$ and the local minimal and maximal steps in the interval $[t_{n-n_0+1}, t_n]$.} 

The total computational cost of the algorithm is \(O(n_0 N_Q^{\mathrm{loc}} {N} + N_Q^{\mathrm{his}}(N - n_0) )\), while the memory requirement is \( {O(N_Q^{\mathrm{his}})}\).  When \( G(x)  \) takes the form  \eqref{Gx_fracint}, our choice of $N_Q^\mathrm{his}$ relies on a generalization of the algorithm in \cite{BanLo19} to variable steps, which can be achieved by following the proofs in this paper and replacing where appropriate $\tau$ by $\tau_{\max}$ or $\tau_{\min}$, cf. \cite{GuoLo}. We have not derived so far rigorous complexity estimates for this more general setting but only tested it numerically. For other types of kernels, we apply the trapezoidal rule instead, see Example~\ref{ex_genker}. For the kernel $k_a$ in Example~\ref{ex_genker} the splitting into a local and a history term can actually be reduced to the case $n_0=1$, in such a way that $N_Q^{\mathrm{loc}}=0$. A more optimized quadrature for each specific kernel might be developed, but we do not pursue this issue further in the present paper. 

\begin{figure}[H]
	\centering
	\begin{subfigure}[t]{0.48\textwidth}
		\centering
		\includegraphics[width=0.9\textwidth]{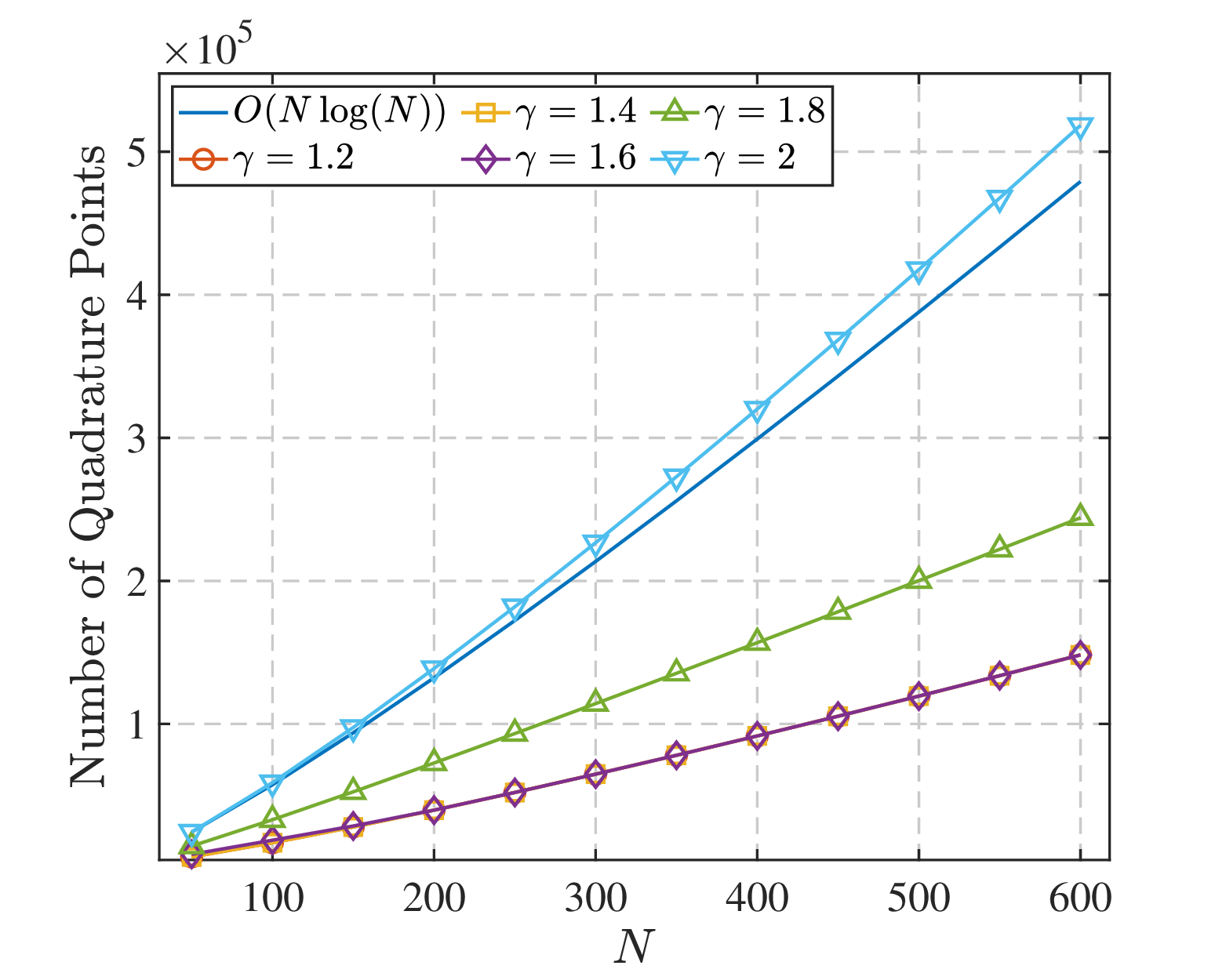}
		\caption{Reference implementation using the quadrature rule from \cite{LoSau15apnum}.}
		\label{subfig:ellip_quad}
	\end{subfigure}
	\hfill
	\begin{subfigure}[t]{0.48\textwidth}
		\centering
		\includegraphics[width=0.9\textwidth]{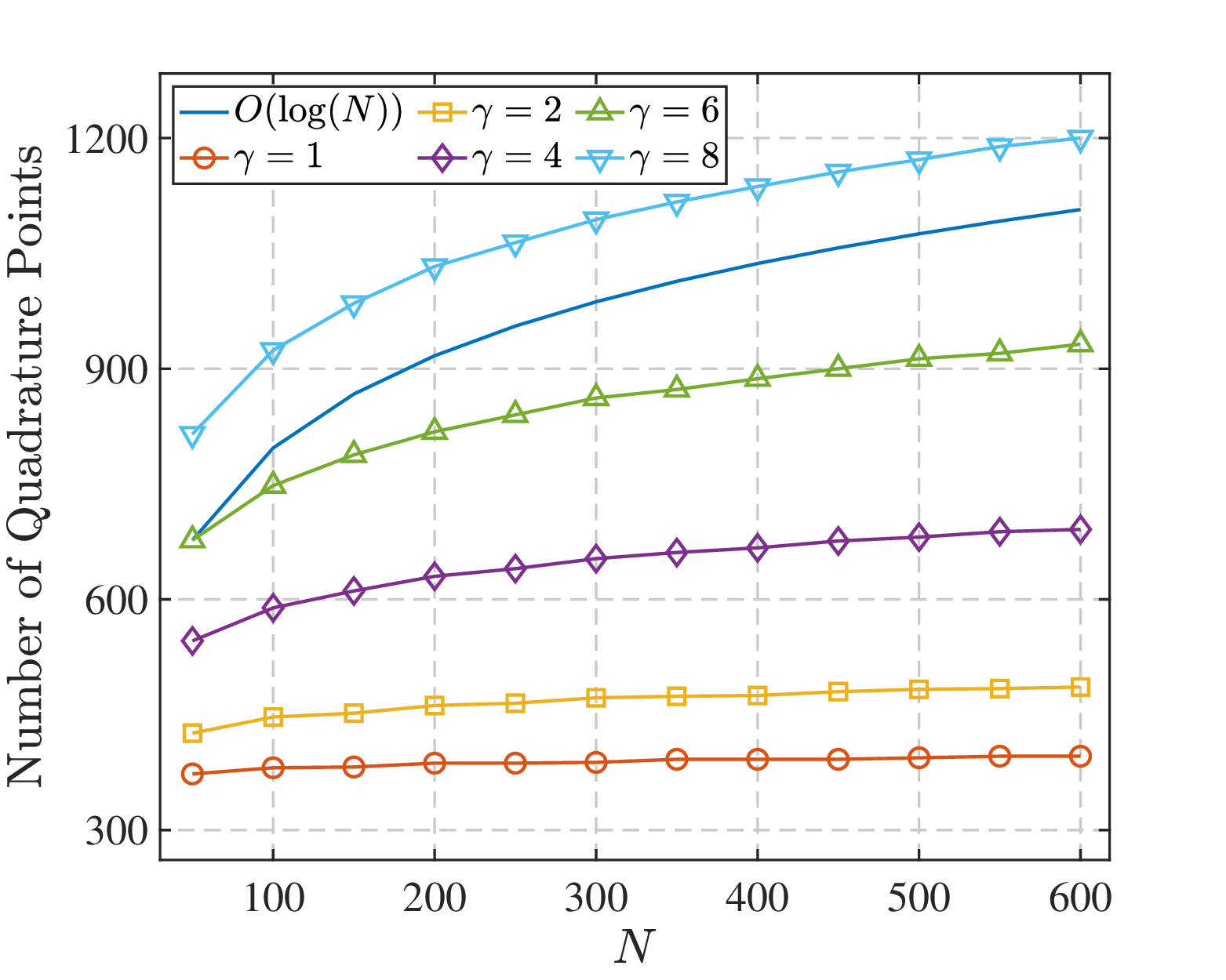}
		\caption{Fast implementation introduced in Section~\ref{sec:algconv}.}
		\label{subfig:fast_quad}
	\end{subfigure}
	\caption{Comparison of the total number of quadrature points used in the gCQ scheme \eqref{gCQshort_split} for Example~\ref{ex_fracint} on the graded mesh \eqref{gmesh} with \( \alpha = 0.5 \) and \( \tol = 10^{-14} \).}
	\label{fig:num_quadcomp}
\end{figure}

\begin{figure}[H]
	\centering
	\begin{subfigure}[t]{0.32\textwidth}
		\centering
		\includegraphics[width=1\textwidth]{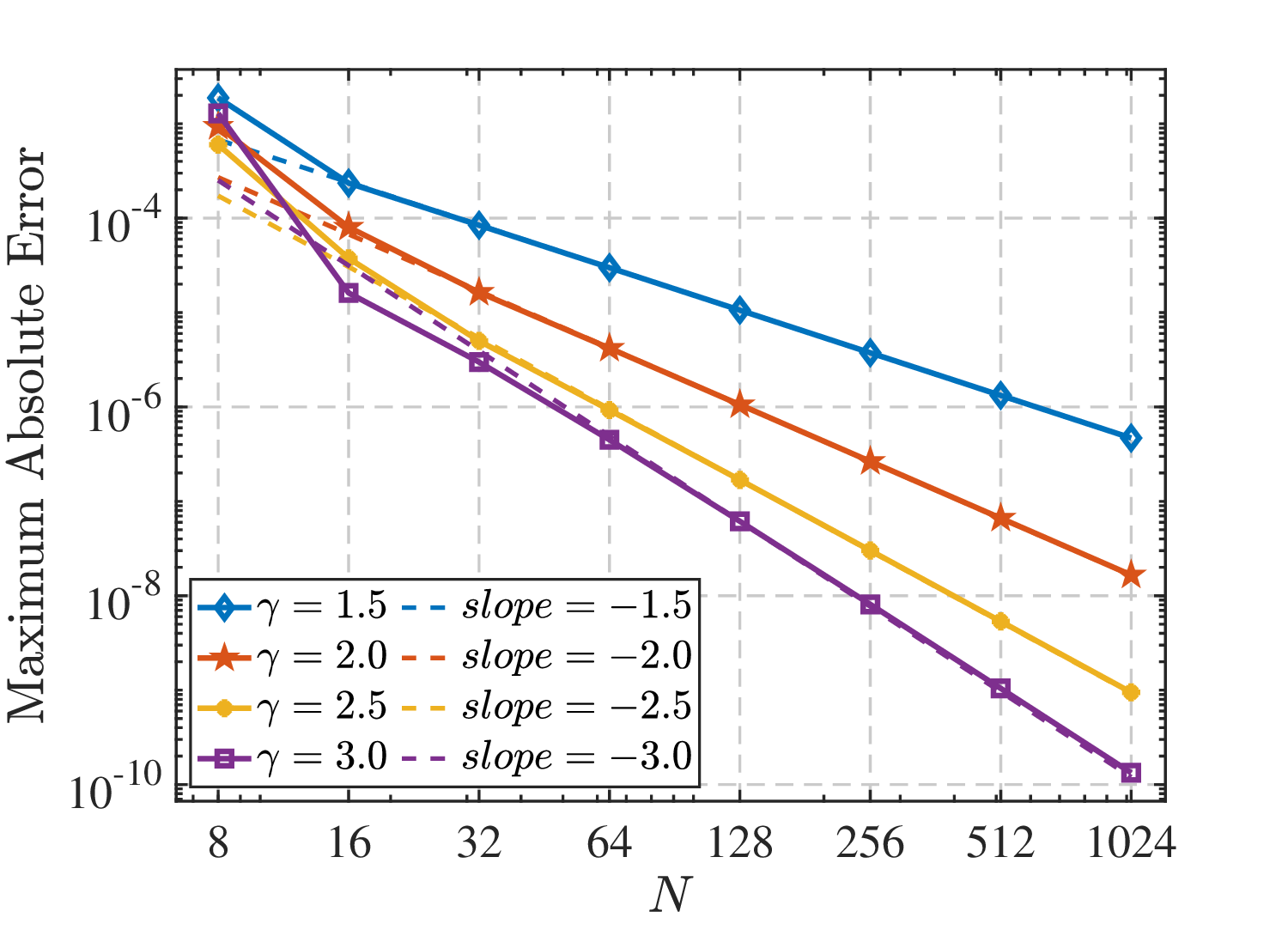}
		\caption{Maximum absolute error.}
		\label{fig:rate_diffgam_tol}
	\end{subfigure}
	\hfill
	\begin{subfigure}[t]{0.32\textwidth}
		\centering
		\includegraphics[width=1\textwidth]{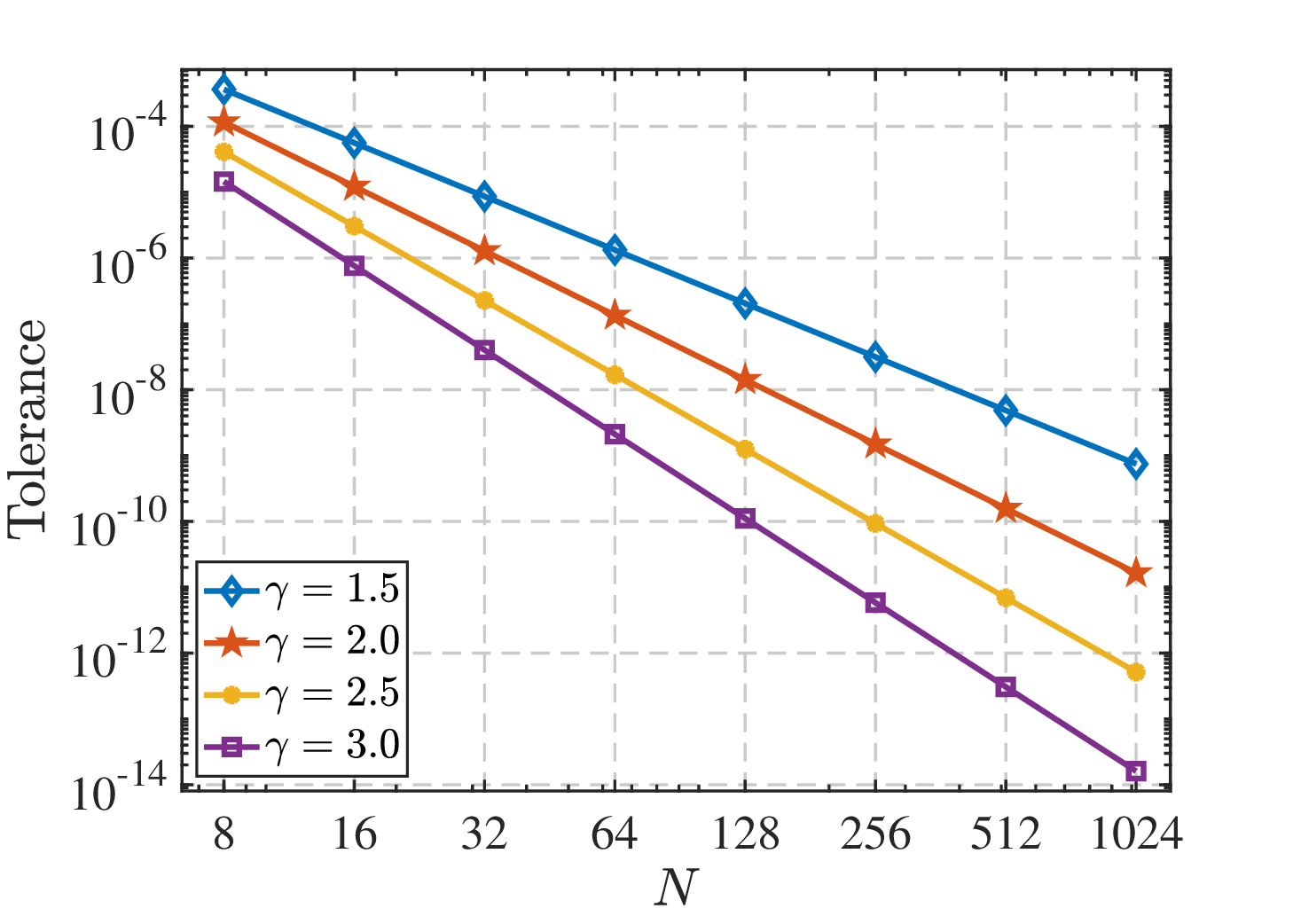}
		\caption{Tolerance values.}
		\label{fig:tol_diffgam}
	\end{subfigure}
	\hfill
	\begin{subfigure}[t]{0.32\textwidth}
		\centering
		\includegraphics[width=1\textwidth]{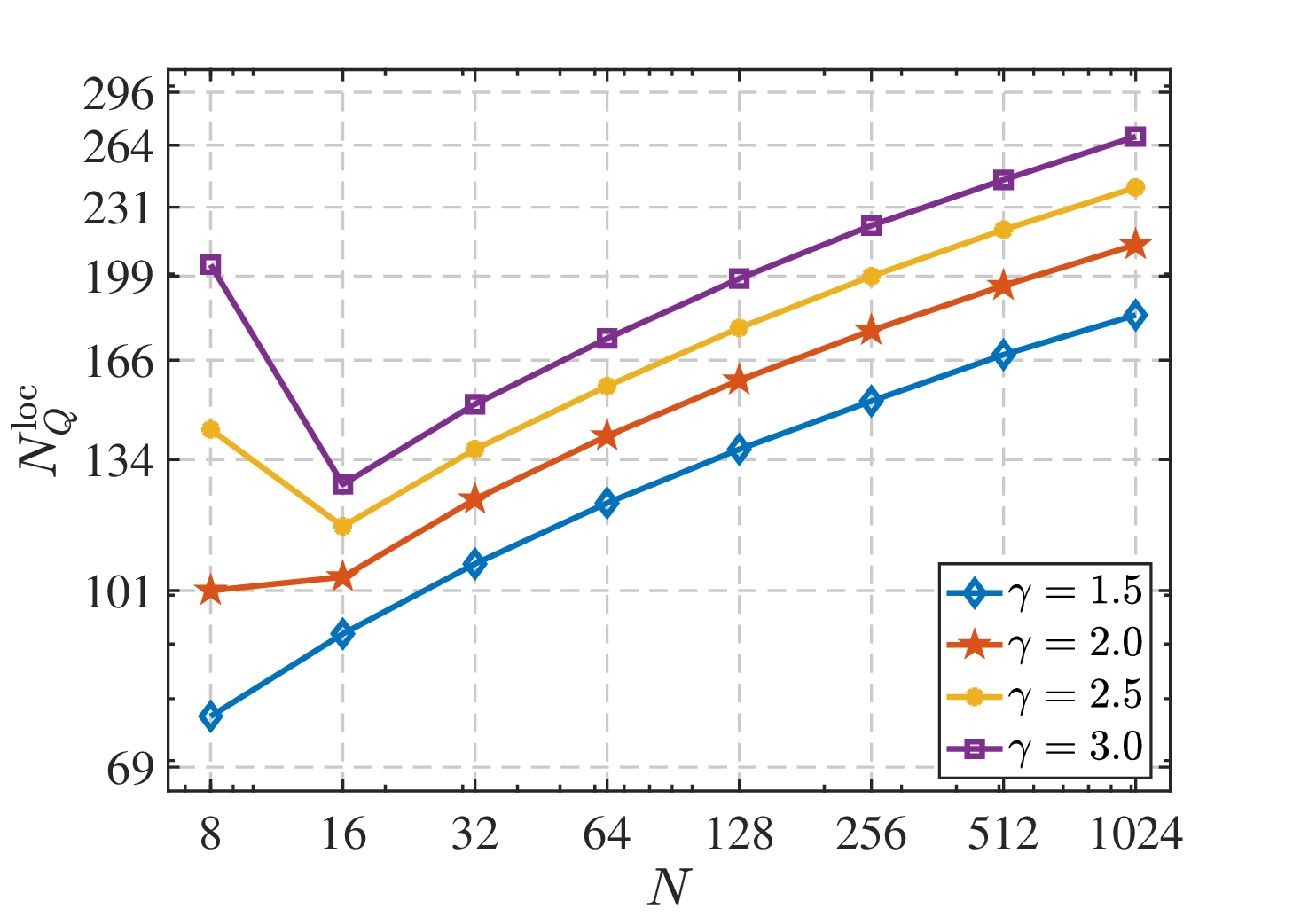}
		\caption{$N_Q^{\mathrm{loc}}$ in \eqref{NQloc} for the interval $[t_1,t_5]$.}
		\label{fig:Nqloc_diffgam_tol}
	\end{subfigure}
	\caption{Numerical results for Example~\ref{ex_fracint} using the two-stage Radau IIA scheme  with $\alpha=0.6$, $\beta=0.4$, on graded meshes \eqref{gmesh} with $\gamma = 1.5,2,2.5,3$ and adaptive tolerance \eqref{tol_adap}.}
	\label{fig:Adaptol_results_diffgam}
\end{figure}

In Figure~\ref{fig:num_quadcomp} we compare the total number of quadrature nodes required to approximate the fractional integral in Example~\ref{ex_fracint} on graded meshes, by using the algorithm in \cite{LoSau16} and by using the fast and oblivious algorithm in this section. As shown, the algorithm in this section significantly improves both efficiency and memory usage compared to the original gCQ  implementation in \cite{LoSau15apnum}. Specifically, the total number of quadrature nodes, this is $N_Q^{\mathrm{loc}}+N_Q^{\mathrm{his}}$, required by our fast algorithm grows like \(O(\log N)\) on the graded mesh \eqref{gmesh}, implying an overall computational complexity of \(O(N \log N)\) and a memory usage of \(O(\log N)\). This growth turns out to be due only to the evaluation of the history term $\Sv^{\mathrm{his}}_n$.  For a fixed high tolerance requirement $\tol = 10^{-14}$ in \eqref{NQloc}, which is independent of $N$, we need to take a large value $N_Q^{\mathrm{loc}} = 273$ to produce the results displayed in Figure~\ref{fig:num_quadcomp}. A better complexity can be achieved if the tolerance is chosen adaptively based on $N$, the total number of time steps, and our error bounds. A careful perturbation analysis for the effect of this quadrature in the global accuracy is beyond the scope of this paper but it is to be expected that this error grows proportionally to some power of $N$ -- see the result in \cite[Theorem 5.3]{HipLoPa2014} for a different application and for standard CQ. With this in mind we have heuristically found  the choice
	{\begin{equation}\label{tol_adap}
			\mathrm{tol} = \frac{1}{10} N^{-\left(\min\left\{ \gamma(\alpha + \beta),\, p,\, q + 1 + \alpha \right\}+\min\left\{\frac{4}{5}\gamma,\frac{5}{4}\right\}\right)} \, \max\limits_{t = t_i, \, 1 \le i \le N} |f(t)|,
	\end{equation}}
	to be good. We have used this choice to produce Figure~\ref{fig:Adaptol_results_diffgam}, which shows that this adaptive tolerance leads to  fewer quadrature points compared to the fixed high tolerance requirement $\tol = 10^{-14}$.

We emphasize that thanks to the improvements in this section we are able to choose very strongly graded meshes, with grading exponent $\gamma$ up to $20$. For the original implementation of the gCQ method in \cite{LoSau15apnum,LoSau16}, the complexity starts to become unaffordable in practice for $\gamma>2$.



\subsection{Algorithm for the linear subdiffusion equation}\label{subsec:subdiff}
This subsection presents the adaptation of the algorithm developed in Section~\ref{sec:algconv} to discretize the linear subdiffusion equation.
\begin{equation}\label{homo_fracdiffusion}
\begin{array}{rcl}
	\partial_t^\alpha v(t) +\mathcal{A}v(t) & = & f(t)-\mathcal{A}u_0,	\quad t\in (0,T].\\
	v(0)&=&0,
\end{array}
\end{equation}
where $0<\alpha<1$, $\mathcal{A}$ denotes an elliptic operator and $u_0$ is given. Define\[\Vv(t_n)=\left(v(t_{n-1}+c_i\tau_n)\right)_{i=1}^s,\quad \Vv_n=\left(V_{ni}\right)_{i=1}^s,\]
with $V_{ni}$ approximating $v(t_{n-1} + c_i \tau_n)$. To apply the fast gCQ  algorithm for the fractional derivative, we observe that, by definition,
\[
\partial_t^{\alpha} =  \partial_t^{\alpha-1} \partial_t,\quad 0<\alpha<1.
\]
Then, according to the discrete composition rule, it also holds that
\[
\left(\partial^{\Delta}_t \right)^{\alpha} =  \left(\partial^{\Delta}_t \right)^{\alpha-1} \partial^{\Delta}_t,
\]
where $\partial^{\Delta}_{t}$ denotes the approximation of the first order derivative provided by the gCQ  method, this is, \cite[(29)]{LoSau16}
\begin{equation}\label{gCQ_1stderv}
\left[\partial^{\Delta}_{t} \Vv\right]_n = (\tau_n \Av)^{-1}\left(\Vv_n-v_{n-1}\bone\right), \quad 1\le n\le N.
\end{equation}
Thus, we can write
\[
\sum_{j=1}^n \Wv_{n,j}  \left[\partial^{\Delta}_{t} \Vv \right]_{j}  + \mathcal{A}\Vv_n  =  \fv(t_n)-\mathcal{A}u_0\bone, \qquad n\ge 1,
\]
where \( \Wv_{n,j} \) are the gCQ weights \eqref{gCQ_w} for the fractional integral operator \( \partial_t^{\alpha-1} \), corresponding to
\[G(x)= \frac{\sin(\pi (1-\alpha))}{\pi} x^{\alpha-1}.\]
Subsequently, we have
\[
\Wv_{n,n}  \left[\partial^{\Delta}_{t} \Vv \right]_{n}  + \mathcal{A}\Vv_n  = \fv(t_n)-\mathcal{A}\Uv_0-\sum_{j=1}^{n-1} \Wv_{n,j}  \left[\partial^{\Delta}_{t} \Vv\right]_{j}, \qquad n\ge 1.
\]
Noting from \eqref{gCQ_wContour} that \(\Wv_{n,n} = (\tau_n \Av)^{1 - \alpha}\), and combining this with \eqref{gCQ_1stderv}, the numerical scheme can be reformulated as
\begin{equation}\label{scheme_subd}
\left(( \tau_n\Av)^{-\alpha}I + \mathcal{A} \right) \Vv_n=  \fv(t_n)-\mathcal{A}\Uv_0+ (\tau_n\Av)^{-\alpha}re v_{n-1}\bone -\sum_{j=1}^{n-1} \Wv_{n,j}  \left[\partial^{\Delta}_{t} \Vv\right]_{j}, \quad n\ge 1.
\end{equation}
The efficient evaluation of the right-hand side summation term employs the algorithm from Section~\ref{sec:algconv} by splitting it into two parts
\begin{equation}
\sum_{j=1}^{n-1} \Wv_{n,j}  \left[\partial^{\Delta}_{t} \Vv \right]_{j}=\sum_{j=1}^{n-n_0} \Wv_{n,j}  \left[\partial^{\Delta}_{t} \Vv \right]_{j} + \sum_{j=\max(1,n-n_0+1)}^{n-1} \Wv_{n,j}  \left[\partial^{\Delta}_{t} \Vv \right]_{j}.
\end{equation}
The history part is approximated via
\[
\sum_{j=1}^{n-n_0} \Wv_{n,j} \left[\partial^{\Delta}_{t} \Vv \right]_{j} \approx \sum_{l=1}^{N_{Q}^{\mathrm{his}}} \varpi_l  G(x_l)\left(\prod_{j=n-n_0+1}^{n} \left(\Rv(-\tau_j x_l)\ev_s^\top\right)\right) \Qv^{\mathrm{his}}_{n-n_0}(x_l),
\]
where \( \Qv^{\mathrm{his}}_{n-n_0} \) is defined as in \eqref{recursion_his}, with \( \fv_j \) taken to be \( \left[\partial^{\Delta}_{t} \Vv \right]_j \).  For the local part, the approximation is given by
\[\sum_{j=\max(1,n-n_0+1)}^{n-1} \Wv_{n,j} \left[\partial^{\Delta}_{t} \Vv \right]_{j} \approx \sum_{l=1}^{N_{Q}^{\mathrm{loc}}} w_l z_l^{\alpha-1} \Rv(\tau_n z_l) \ev_s^\top \Qv^{\mathrm{loc}}_{n-1}(z_l),
\]
where $\Qv^{\mathrm{loc}}_{n-1}$  satisfies \eqref{recursion_loc}, with $\fv_{n-1}$ replaced by $ \left[\partial^{\Delta}_{t} \Vv \right]_{n-1}$.

\section{Numerical experiments}\label{sec:num_test}
In this section, we consider the application of Runge--Kutta based gCQ to compute convolution integrals, solve fractional diffusion equations, and solve nonlinear wave equations with convolution-type damping, in order to test the algorithm in Section~\ref{sec:algconv} and verify the theoretical findings presented in Theorem~\ref{thm:conv_gCQrkGrad}. The solutions of the problems under consideration exhibit nonsmooth behavior  {toward} the origin and in some cases also at times away from the origin. In all cases, the approximation on uniform time meshes suffers from some sort of order reduction. To achieve full-order convergence, we have implemented  the gCQ scheme  {with the two-stage Radau IIA method of order $p=3$ and stage order $q=2$, and the four-stage Lobatto IIIC method of order $p=6$ and stage order $q=3$,} on graded meshes with a quadrature tolerance $\tol=10^{-14}$ for the fast and oblivious algorithm. When possible, the grading parameter is selected according to Theorem~\ref{thm:conv_gCQrkGrad} and, in the case of Example 5, which is not fully covered by our theory, heuristically. Notice that in this last example we deal with a nonlinear wave problem and error analysis techniques based uniquely on the preservation of the composition rule by the gCQ  method are not enough to justify the observed high order of convergence, cf. \cite{BaBanPtas24}. 
\subsection{Computation of convolution with different kernels}
\begin{figure}[H]
\centering
\includegraphics[width=0.48\textwidth]{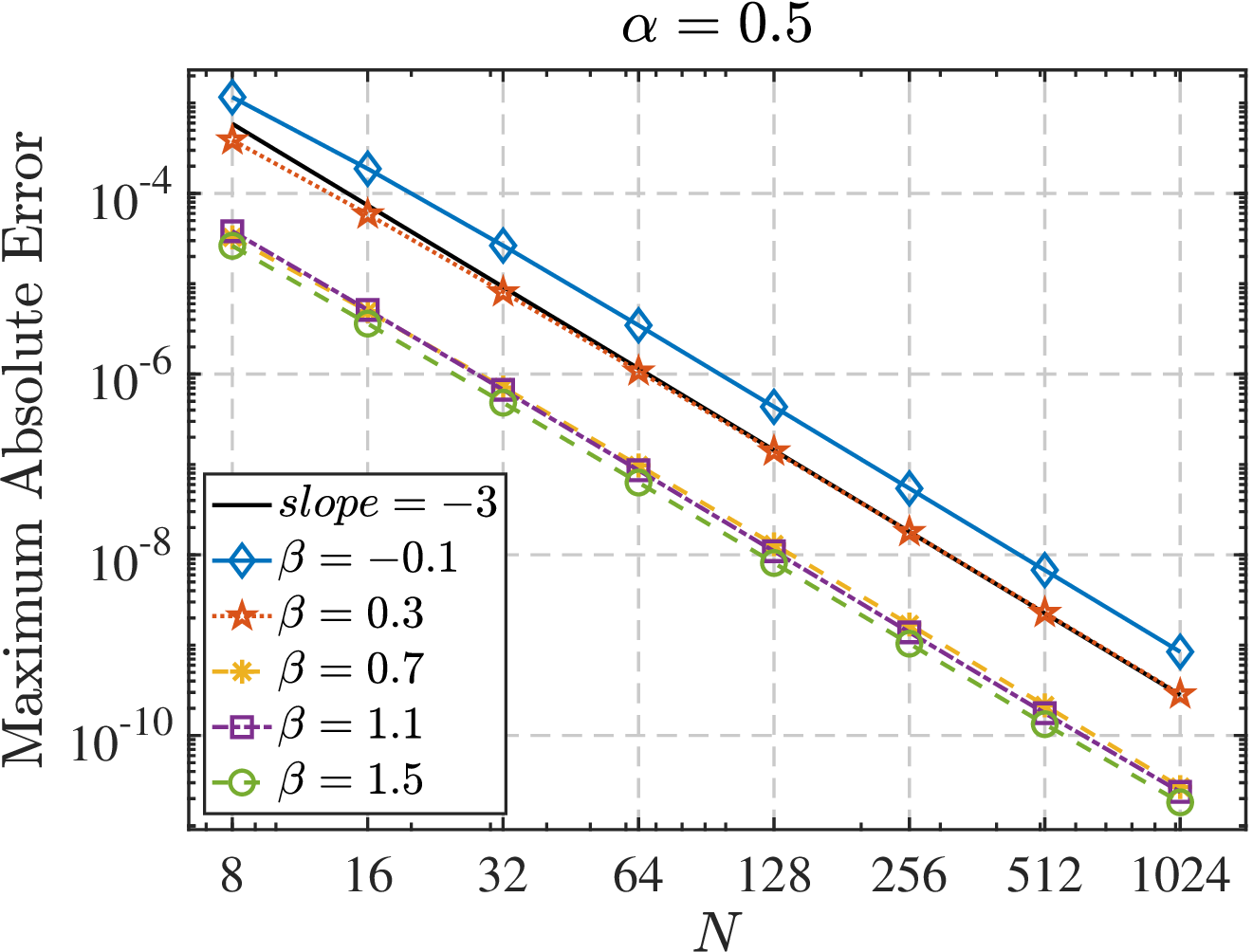}
\hfill
\includegraphics[width=0.48\textwidth]{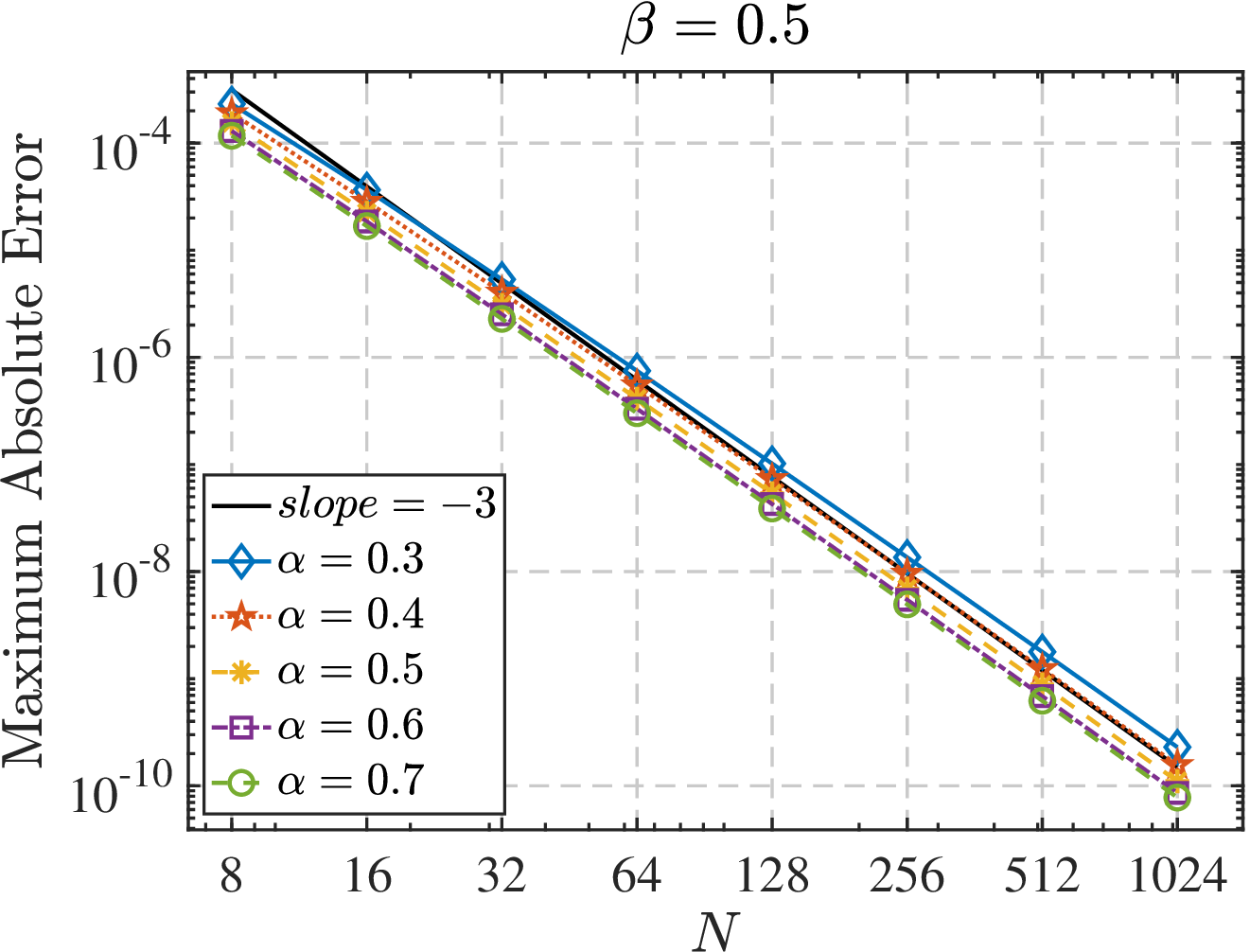}
\caption{Maximum absolute error of the two-stage Radau IIA scheme for Example~\ref{ex_fracint}  on \eqref{gmesh} with $\gamma=\frac{3}{\alpha+\beta}$,  { (Left: $\alpha = 0.5$ fixed and $\beta = -0.1,\,0.4,\,0.7,\,1.1,\,1.5$.
		Right: $\beta = 0.5$ fixed and $\alpha = 0.3,\,0.4,\,0.5,\,0.6,\,0.7$.)}}
\label{fig:fracint_orderOpt}
\end{figure}

\begin{example}\label{ex_fracint}
Consider the  Runge--Kutta  based gCQ approximation \eqref{gCQ_appInt} for  the fractional integral  
\begin{equation}\label{fracint}
	\partial_t^{-\alpha}f(t)=\int_0^t\frac{(t-s)^{\alpha-1}}{\Gamma(\alpha)} f(s)\, ds,
\end{equation}
Let $f(t)=t^\beta$. The exact solution is given by
\begin{equation}\label{fracint_exact}
	\partial_t^{-\alpha}t^\beta=\frac{\Gamma(\beta+1)}{\Gamma(\alpha+\beta+1)}t^{\alpha+\beta}.
\end{equation}
\end{example}
Set $T=1$.  The maximum absolute error in the numerical results is measured by
\begin{equation}\label{MaxAbsErr}
\max_{1\le n\le N}\left| \left[\partial_t^{-\alpha}f\right](t_n) - \left[\left( \partial^{\Delta}_t \right) ^{-\alpha}f \right]_n \right|.
\end{equation}
By applying the two-stage Radau IIA based gCQ  scheme and selecting different values for $\alpha$ and $\beta$, we illustrate the maximum absolute error on the graded mesh \eqref{gmesh} with $\gamma=3/(\alpha+\beta)$ in Figure~\ref{fig:fracint_orderOpt}. The results demonstrate third-order convergence, which aligns with the theoretical findings in Theorem~\ref{thm:conv_gCQrkGrad}.

\begin{example}\label{ex_genker}
Consider the Runge--Kutta based gCQ  approximation \eqref{gCQ_appInt} for convolutions \eqref{conv} with the following two kernels:
\begin{equation}\label{k_genker}
	k_a(t) = -\frac{d}{dt} E_{\alpha,1}(-t^\alpha), \quad k_b(t) = \frac{t^{\alpha-1}}{\Gamma(\alpha)} e^{-t},
\end{equation}
acting on the function
\begin{equation*}
	f(t) = t^\beta.
\end{equation*}
\end{example}
The Laplace transforms of the kernels in \eqref{k_genker} are given by
\[
K_a(z) = \frac{1}{z^\alpha + 1}, \quad K_b(z) = \frac{1}{(z + 1)^\alpha},
\]
both satisfy Assumption~\ref{assumptionK}. The gCQ  approximations take the form
\begin{align*}
\left[ K_a(\partial_t^\Delta) f \right]_n = \int_{0}^{\infty} G_a(x) y_n^a(x)\, dx, \quad
\left[ K_b(\partial_t^\Delta) f \right]_n = \int_{0}^{\infty} G_b(x) y_n^b(x)\, dx,
\end{align*}
where
\begin{equation*}
G_a(x) = \frac{\sin(\pi \alpha)}{\pi}  \frac{x^\alpha}{x^{2\alpha} + 2x^\alpha \cos(\pi \alpha) + 1}, \quad G_b(x) = \frac{\sin(\pi \alpha)}{\pi}  x^{-\alpha}, \quad 0 < \alpha < 1,
\end{equation*}
and
\[
y_n^a(x) = y_n(x), \quad y_n^b(x) = y_n(x + 1),
\]
with \( y_n(x) \) given in \eqref{RK_sol}. 

Notice that the function $G_a$ does not fit the behavior for which the quadrature in \cite{BanLo19} was designed. On the one hand, $G_a$ has poles, with zeros of its denominator at $x^\alpha=\e^{\pm \opi \pi(1-\alpha)}$. On the other hand, it is regular at $x=0$, where it is zero. Thus, we use trapezoidal quadrature to evaluate the  history part of the convolution \( \left[ K_a(\partial_t^\Delta) f \right]_n \). Notice that the substitution \( x = e^{\mu/\alpha} \) transforms the integral \eqref{S_hisint}, with \( G(x) = G_a(x) \), into
\begin{align}
	\label{SvnGa_his}
\Sv_n^{\mathrm{his}}
= \int_{-\infty}^{\infty}g(\mu)\,d\mu,
\end{align}
where
	\begin{equation}\label{intgrand_g}
g(\mu):=\frac{1}{\alpha}\e^{\mu/\alpha}G_a(\e^{\mu/\alpha})
\left(\prod_{j=n-n_0+1}^{n}
\Rv(-\tau_j\e^{\mu/\alpha})\ev_s^\top\right)
\Qv^{\mathrm{his}}_{n-n_0}(\e^{\mu/\alpha}).
	\end{equation}
Applying the trapezoidal rule yields the approximation  
\begin{align}
\label{SvnGa_Trap}
\Sv_n^{\mathrm{his}}
\approx h\sum_{l=-\widetilde M}^{\widetilde M}g(lh).
\end{align}
Let \(z=\mu+\opi\eta\). Under \(x=\e^{z/\alpha}\), the zeros of the denominator
of \(G_a(x)\) identified above correspond to
\(\e^z=\e^{\pm\opi\pi(1-\alpha)}\). Hence, the nearest poles of
\(G_a(\e^{z/\alpha})\) are \(z=\pm\opi\pi(1-\alpha)\), so this factor is
holomorphic for \(|\eta|<\pi(1-\alpha)\). The Runge--Kutta factors in \eqref{intgrand_g}
involve resolvents of the form
\[
(\Iv-\rho_j\Av)^{-1},\qquad \rho_j=-\tau_j\e^{z/\alpha}.
\]
Assumption~\ref{ass_RK} guarantees that these resolvents are nonsingular
whenever \(\Re\rho_j\le0\). For \(|\eta|<\pi\alpha/2\),
\[
\Re\rho_j=-\tau_j\e^{\mu/\alpha}\cos(\eta/\alpha)<0.
\]
Thus, the resolvents are holomorphic in this strip, and
\eqref{recursion_his} implies the same for
\(\Qv_{n-n_0}^{\mathrm{his}}(\e^{z/\alpha})\). Combining this with
the kernel analyticity, the integrand \(g\) in \eqref{intgrand_g} is
continuous on the closed strip
\begin{equation}
	\label{Dd}
D_d:=\{z\in\mathbb C:|\Im z|\le d\},\qquad
0<d<\min\left\{\pi(1-\alpha),\frac{\pi\alpha}{2}\right\},
\end{equation}
and holomorphic in its interior.
Moreover, the estimates in Lemmas~\ref{lem:Rz_bnd} and \ref{lem:gCQ w_bnd} and the asymptotic behavior of $G_a(\e^{\mu/\alpha})$, as $\mu \to \pm\infty$,  imply the bounds, uniformly for \(|\eta|\le d\),
\begin{equation}
	\label{g_etabnd}
\|g(\mu+\opi\eta)\|\le C_0\e^{(1+1/\alpha)\mu}
\quad \mu\to-\infty,
\qquad
\|g(\mu+\opi\eta)\|\le C_0\e^{-\mu}
\quad\mu\to+\infty.
\end{equation}
Since \(1+1/\alpha>1\) for $0<\alpha<1$, on the real axis,  we have
\begin{equation}
	\label{g_asympbnd}
\|g(\mu)\|\le C_0\e^{-|\mu|},
\qquad \mu\in\mathbb R.
\end{equation}

By \cite[Theorem~3.4]{Steng}, conditions
\eqref{Dd}--\eqref{g_asympbnd} imply that the error of the
trapezoidal rule in \eqref{SvnGa_Trap} satisfies
\[
\left\|\int_{-\infty}^{\infty}g(\mu)\,d\mu
-h\sum_{l=-\widetilde M}^{\widetilde M}g(lh)\right\|
\le C\e^{-\sqrt{2\pi d\widetilde M}},
\qquad
h=\sqrt{\frac{2\pi d}{\widetilde M}}.
\]
Based on this bound and the graded-mesh gCQ error estimate \(O\!\left(N^{-\min\{p,q+1+\alpha,\gamma(\alpha+\beta)\}}\right)\) in
Theorem~\ref{thm:conv_gCQrkGrad}, we set
\begin{equation*}
h=\sqrt{\frac{2\pi d}{\widetilde M}},\quad \widetilde M=
\left\lceil
\frac{\left(\min\{p,q+1+\alpha,\gamma(\alpha+\beta)\}\log N+2\right)^2}
{2\pi d}
\right\rceil,
\end{equation*}
with \(d=0.9\min\left\{\pi(1-\alpha),\frac{\pi\alpha}{2}\right\}\) in the numerical test for evaluating \eqref{SvnGa_Trap}.

\begin{figure}[H]
\centering
\begin{subfigure}[t]{0.49\textwidth}
	\centering
	\includegraphics[width=1\textwidth]{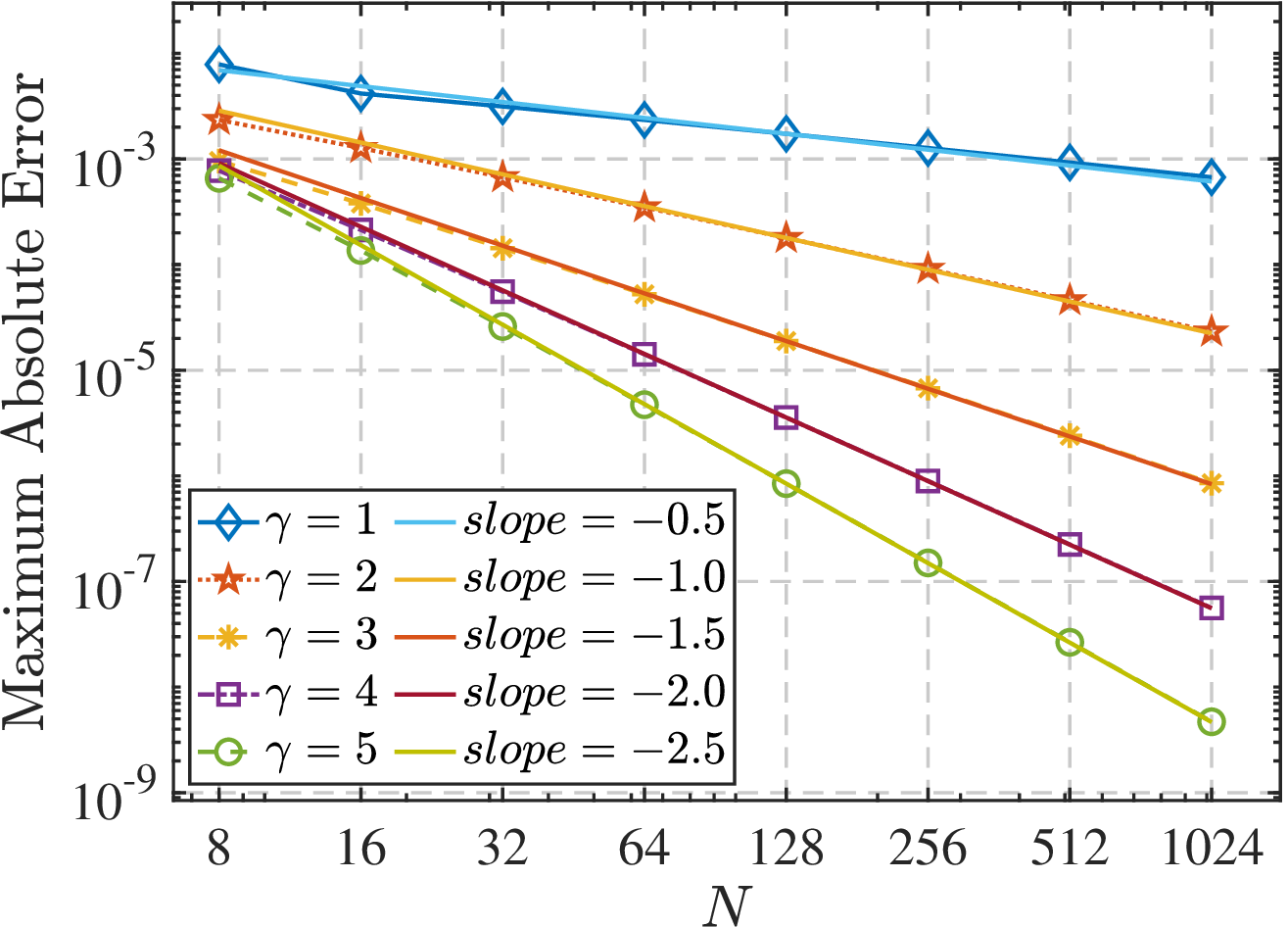}
	\caption{Maximum absolute error.}\label{fig:conv_kadiffGam}
\end{subfigure}
\hfill
\begin{subfigure}[t]{0.49\textwidth}
	\centering
	\includegraphics[width=1\textwidth]{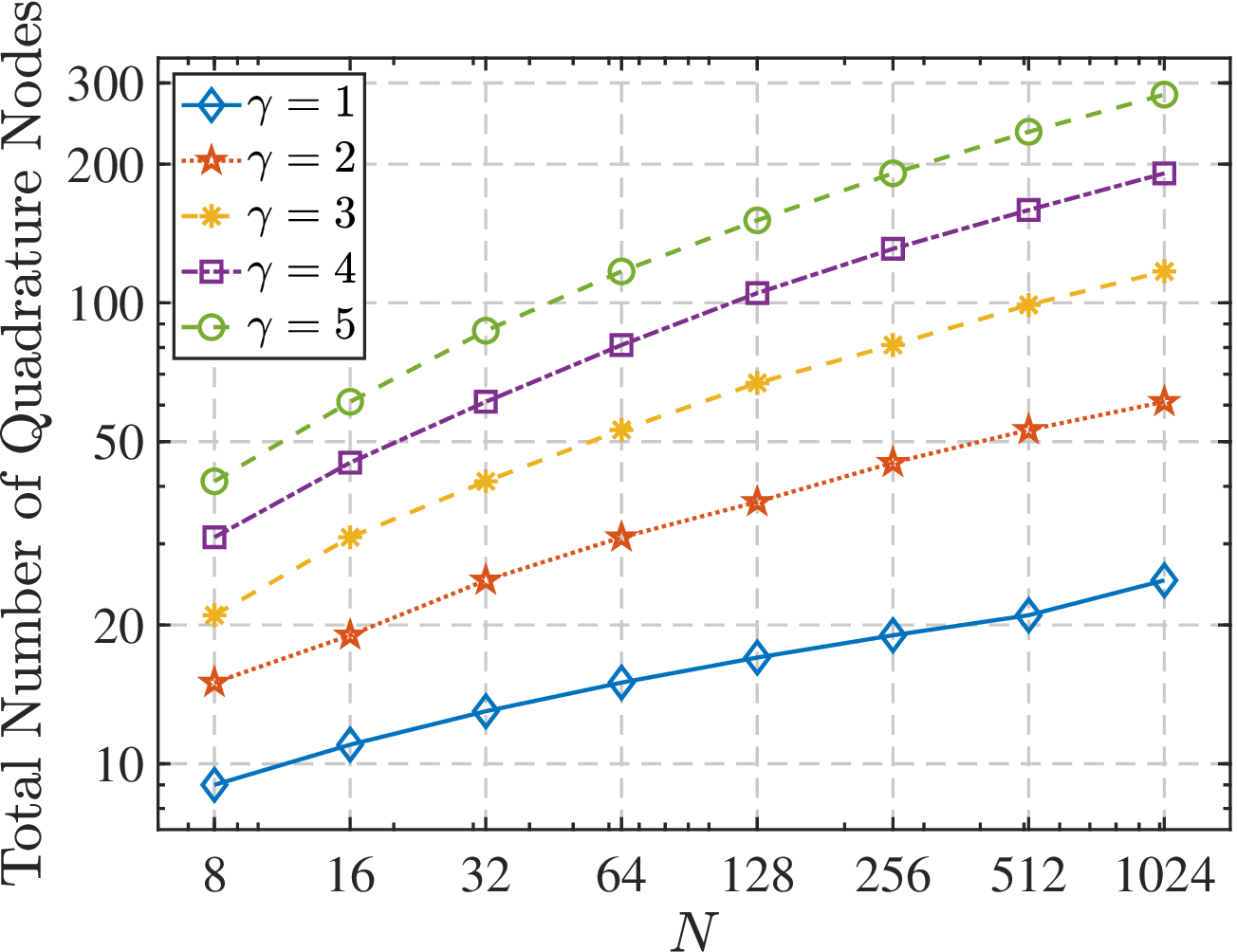}
	\caption{Total number of quadrature nodes.}\label{fig:conv_kadiffGamWidtM}
\end{subfigure}
\caption{Numerical results for the two-stage Radau IIA scheme for Example~\ref{ex_genker} with $k=k_a$, $\alpha = 0.3$, $\beta = 0.2$  and varying $\gamma = 1,\,2,\,3,\,4,\,5$.}\label{fig:conv_ka}
\end{figure}

\begin{figure}[H]
\centering
\begin{subfigure}[t]{0.49\textwidth}
	\centering
	\includegraphics[width=1\textwidth]{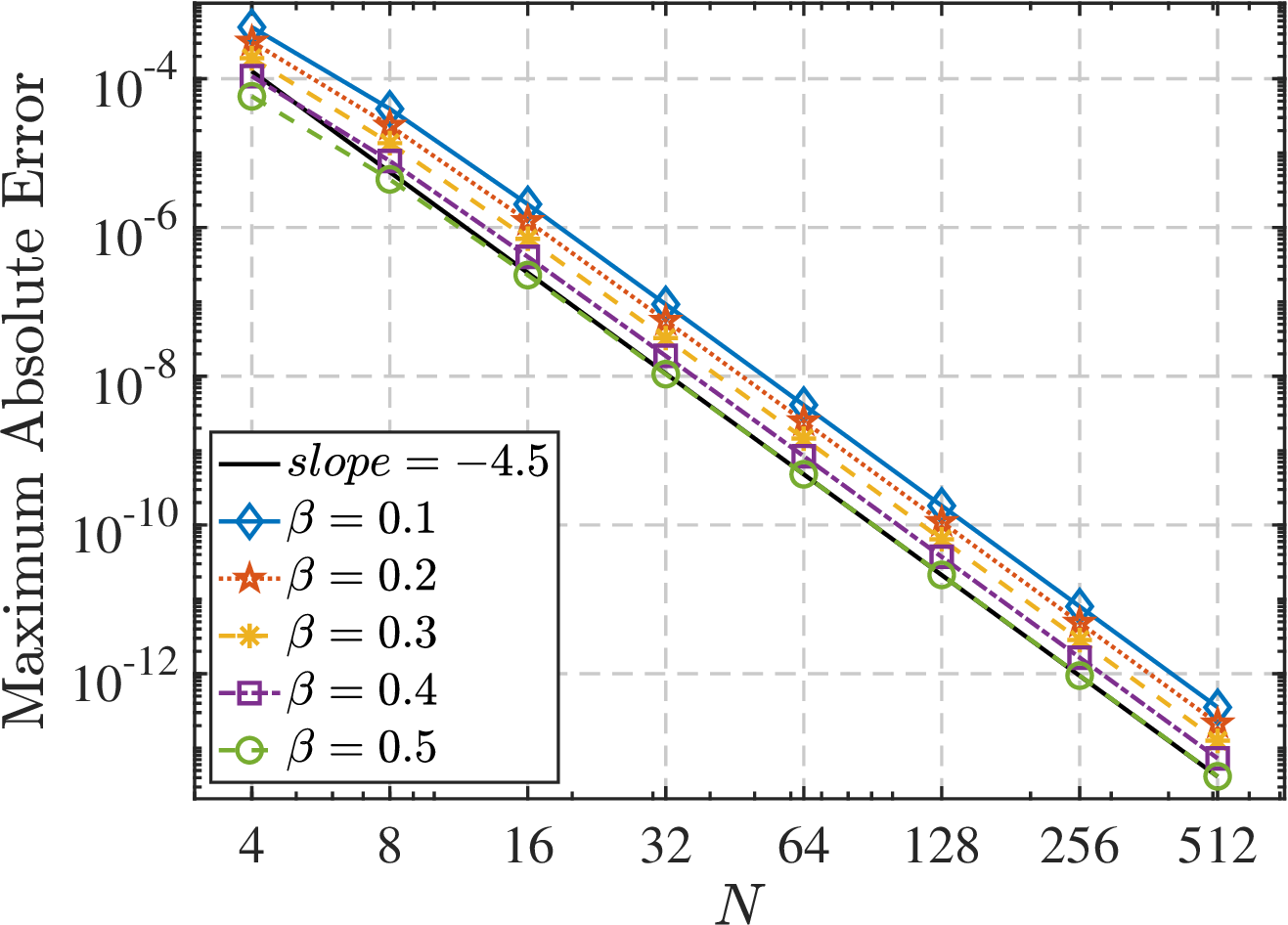}
	\caption{$\alpha=0.5$  { fixed and $\beta = 0.1,\,0.2,\,0.3,\,0.4,\,0.5$,} using $\gamma=\frac{4.5}{\alpha+\beta}$.}\label{fig:conv_kadiffBetLobalp05}
\end{subfigure}
\hfill
\begin{subfigure}[t]{0.49\textwidth}
	\centering
	\includegraphics[width=1\textwidth]{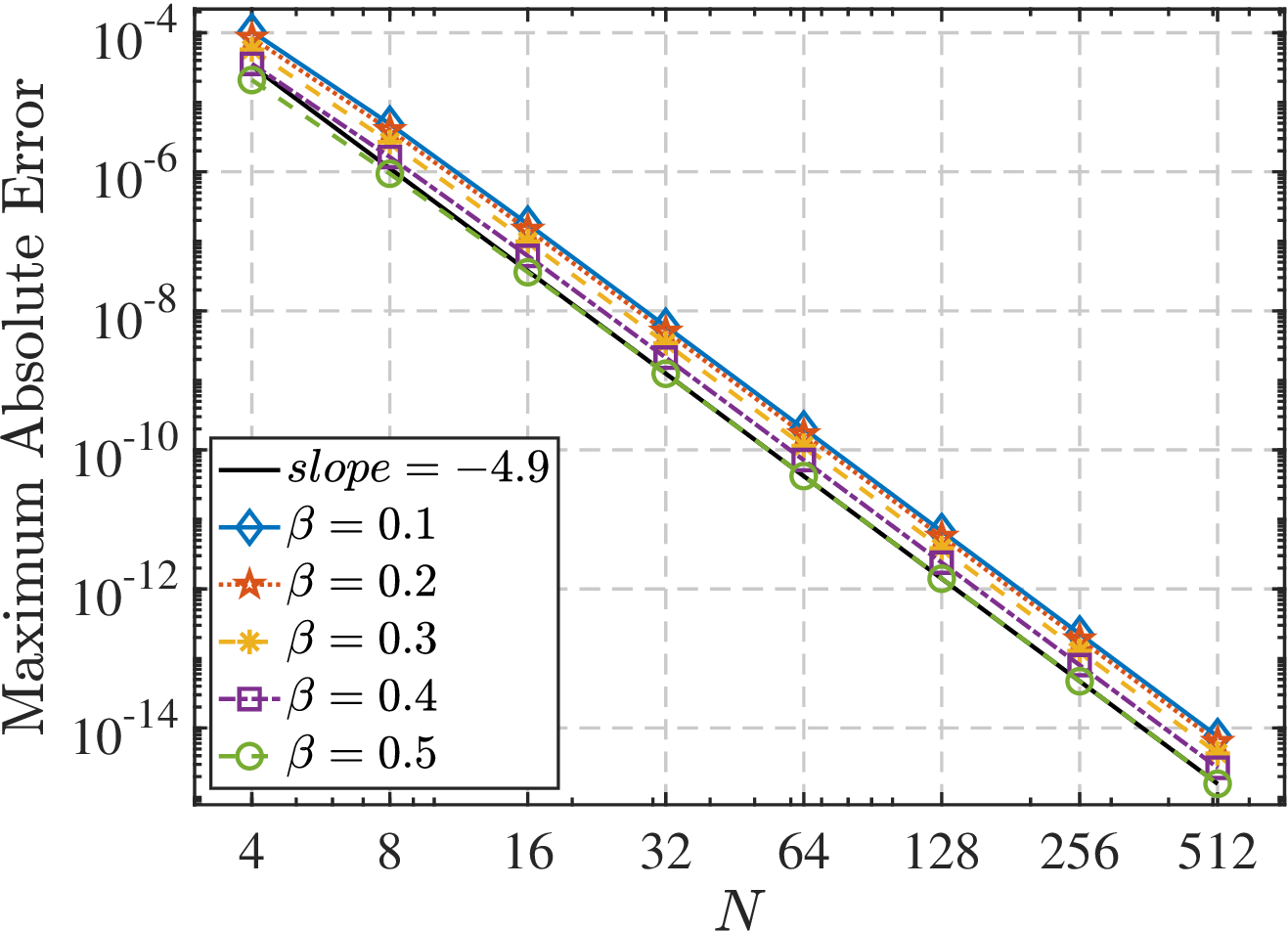}
	\caption{$\alpha=0.9$  { fixed and $\beta = 0.1,\,0.2,\,0.3,\,0.4,\,0.5$,} using $\gamma=\frac{4.9}{\alpha+\beta}$.}\label{fig:conv_kadiffBetLobalp09}
\end{subfigure}
\caption{Maximum absolute error of the four-stage Lobatto IIIC scheme for Example~\ref{ex_genker} with $k=k_a$.}\label{fig:conv_kadiffBetLob}
\end{figure}

The behavior of $G_a$ is actually more favorable from the implementation point of view than the case considered in \cite{BanLo19}. It actually allows to reduce the splitting into a local and a history part in the discrete convolution to the case $n_0=1$ and use a unique quadrature rule from the beginning of the time integration procedure, to deal with the history. Indeed, setting $n_0=1$ in Section~\ref{sec:algfi}, the local contribution \eqref{Snloc} involves only the current time step and the evaluation of \(K((\tau_n \Av)^{-1})\fv_n\), that can be directly obtained by diagonalizing $\Av$.  In this case, no quadrature nodes are needed to deal with the local part of the convolution, i.e., we set \(N_Q^{\mathrm{loc}}=0\). The results of these choices are reported in Figures~\ref{fig:conv_ka} and~\ref{fig:conv_kadiffBetLob}. The total number of quadrature nodes is \(N_Q^{\mathrm{loc}}+N_Q^{\mathrm{his}}=2\widetilde M+1\).

The maximum absolute error for \( K(z) = K_a(z) \) is computed in the same manner as in \eqref{MaxAbsErr}. The exact solution is given by  
\[
(K_a(\partial_t)f)(t) = \Gamma(\beta + 1) t^{\alpha + \beta} E_{\alpha, \alpha + \beta + 1}(-t^\alpha),
\]
with the Mittag-Leffler function evaluated using the routine from \cite{Gar15}.

The results in Figure~\ref{fig:conv_kadiffGam} show that the two-stage
Radau IIA based gCQ method converges with order
\(\min\{3,\gamma(\alpha+\beta)\}\) on the graded mesh \eqref{gmesh}.
Figure~\ref{fig:conv_kadiffGamWidtM} demonstrates that the required total
number of quadrature points increases moderately as the gCQ approximation becomes more accurate. Additionally, Figure~\ref{fig:conv_kadiffBetLob} presents the case of \(p > q + 1 + \alpha\), verifying that the  Runge--Kutta  based gCQ method can attain a convergence order of \(\min\{p, q + 1 + \alpha\}\) and maintain this rate for nonsmooth \(f(t)\) on a graded mesh \eqref{gmesh} with \(\gamma = \min\{p, q + 1 + \alpha\}/(\alpha + \beta)\). These results are consistent with  Theorems~\ref{thm:conv_gCQrkGenM} and ~\ref{thm:conv_gCQrkGrad}.

\begin{figure}[H]
\centering
\begin{subfigure}[t]{0.49\textwidth}
	\centering
	\includegraphics[width=1\textwidth]{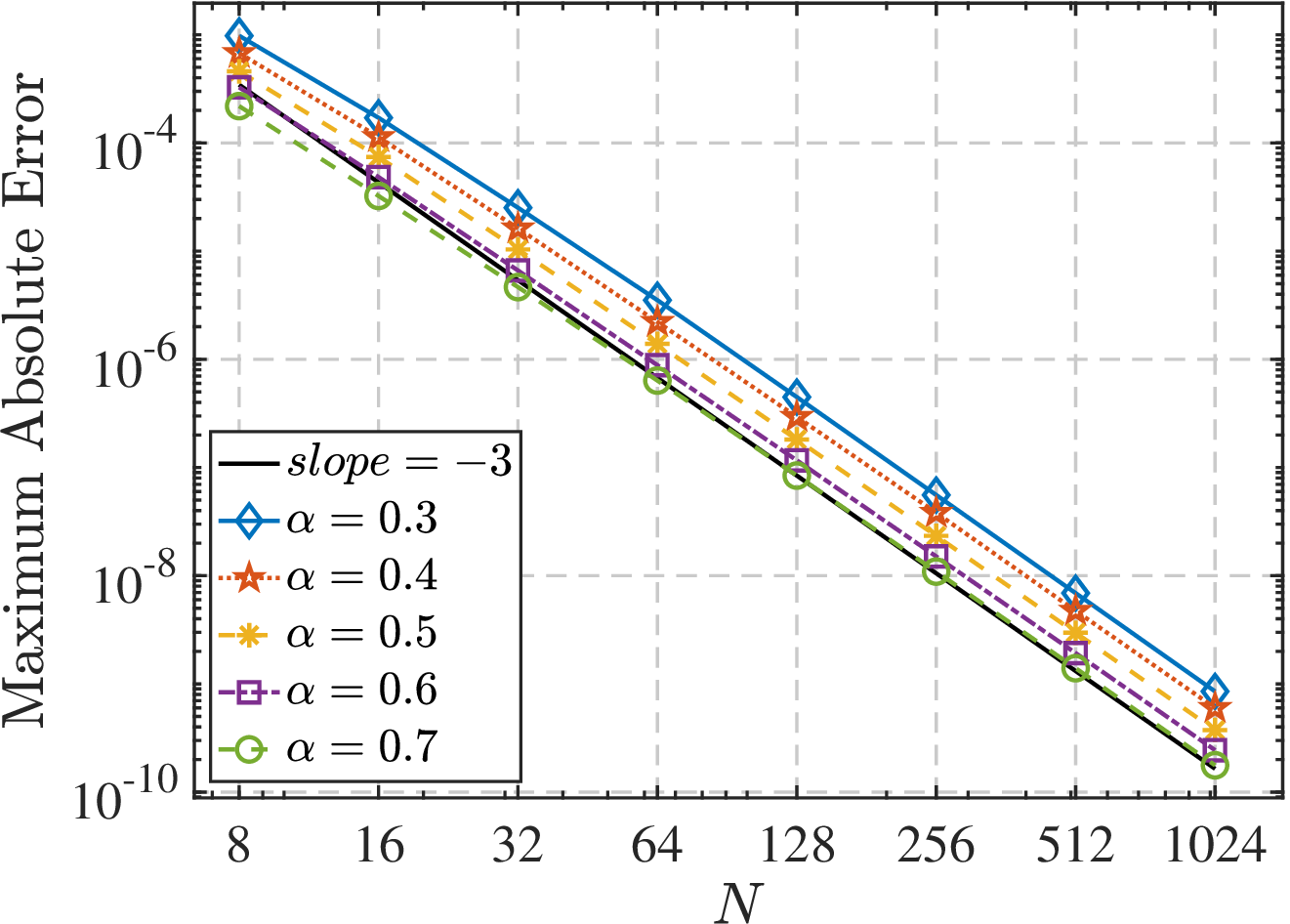}
	\caption{$\beta=0.2$  { fixed and $\alpha = 0.3,\,0.4,\,0.5,\,0.6,\,0.7$,} using  $\gamma=\frac{3}{\alpha+\beta}$.}\label{fig:conv_kboptGam}
\end{subfigure}
\hfill
\begin{subfigure}[t]{0.49\textwidth}
	\centering
	\includegraphics[width=1\textwidth]{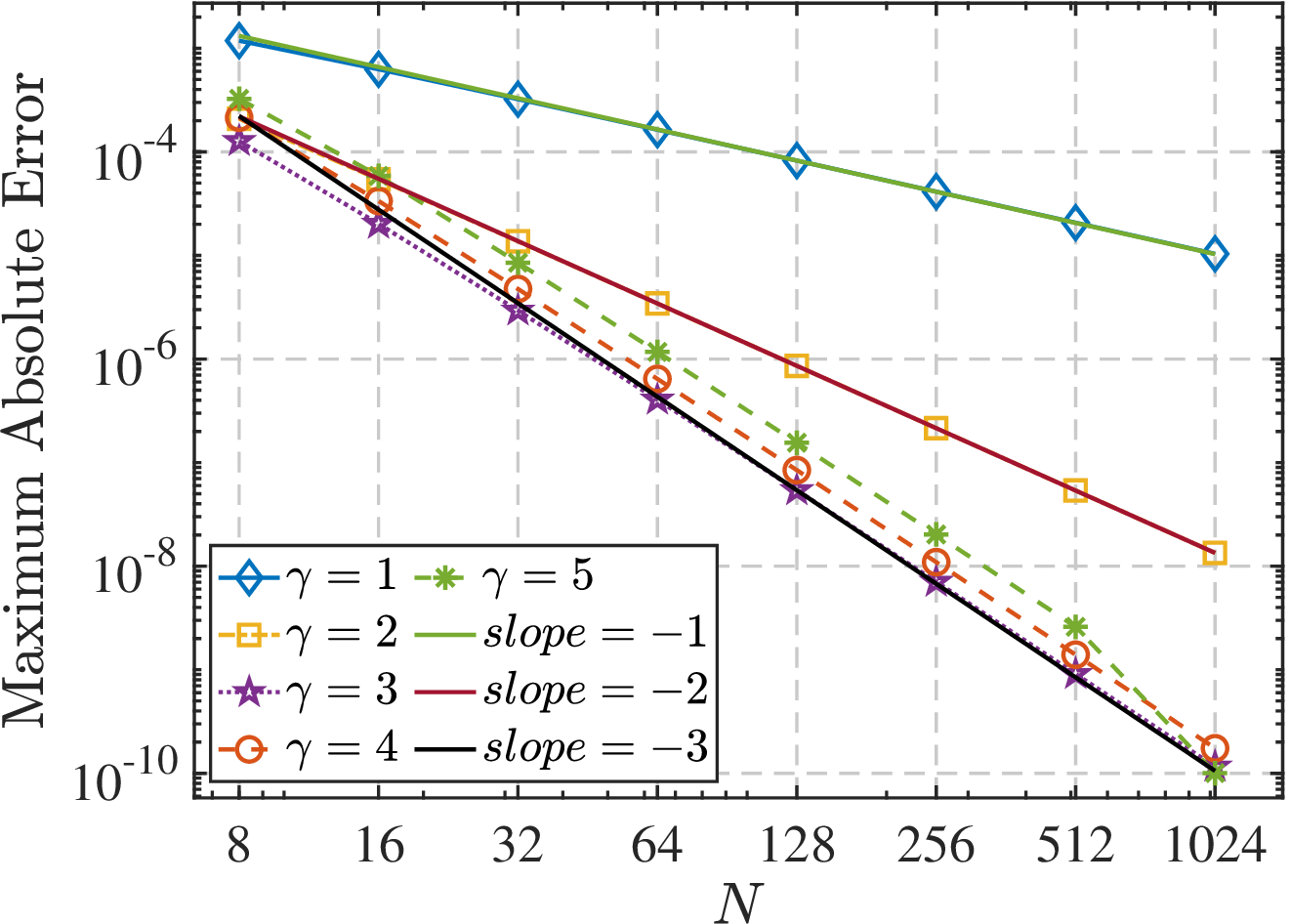}
	\caption{$\alpha=0.4$, $\beta=0.6$  {and varying $\gamma = 1,\,2,\,3,\,4$}.}\label{fig:conv_kbdiffGam}
\end{subfigure}
\caption{Maximum absolute error of the two-stage Radau IIA scheme for Example~\ref{ex_genker} with $k=k_b$.}\label{fig:conv_kb}
\end{figure}

Let $u_n^{\mathrm{ref}}$ denote the numerical solution \(\left[ K_b \left( \partial^{\Delta}_t \right) f \right]_n\) computed on the finer graded time mesh
\begin{equation}\label{gmesh_finer}
t_n = \left(\frac{n}{2N}\right)^\gamma, \quad 1 \leq n \leq 2N,  
\end{equation}
the maximum absolute error in Figure~\ref{fig:conv_kb} is computed by
\begin{equation*}
\max_{1 \leq n \leq N} \left\vert  u_{2n}^{\mathrm{ref}} - \left[ K_b \left( \partial^{\Delta}_t \right) f \right]_n \right\vert.
\end{equation*}

The results in Figure~\ref{fig:conv_kb} show that the fast algorithm for computing $\left[ K_b \left( \partial^{\Delta}_t \right) f \right]_n$ produces the predicted convergence rates on strongly graded meshes, in agreement with Theorem~\ref{thm:conv_gCQrkGrad}.

\subsection{Differential equations}\label{subsec:expFDE}

\begin{figure}[H]
\centering
\begin{subfigure}[t]{0.49\textwidth}
	\centering
	\includegraphics[width=1\textwidth]{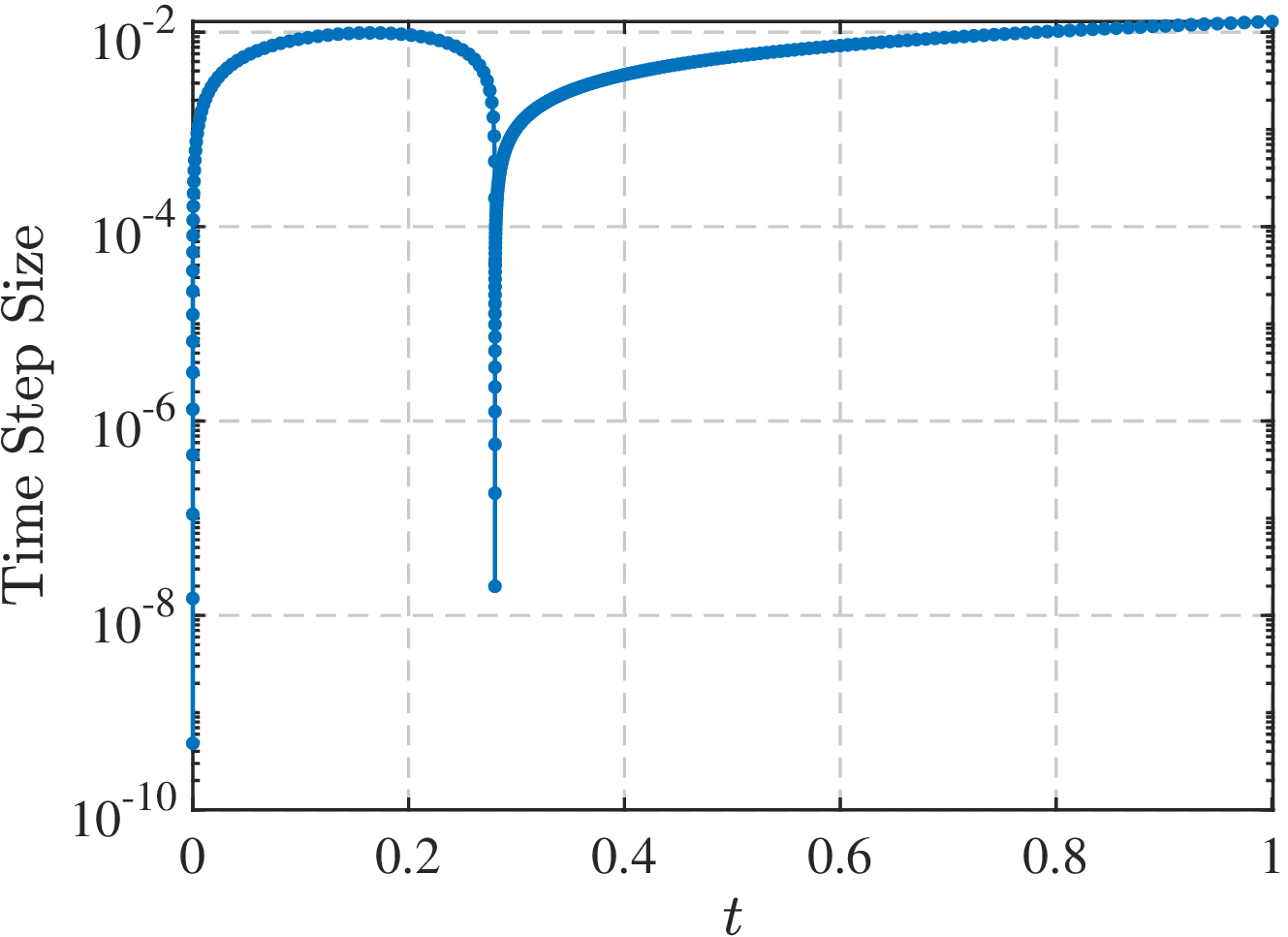}
	\caption{Mesh.}\label{fig:ODE_Mesh}
\end{subfigure}
\hfill
\begin{subfigure}[t]{0.49\textwidth}
	\centering
	\includegraphics[width=1\textwidth]{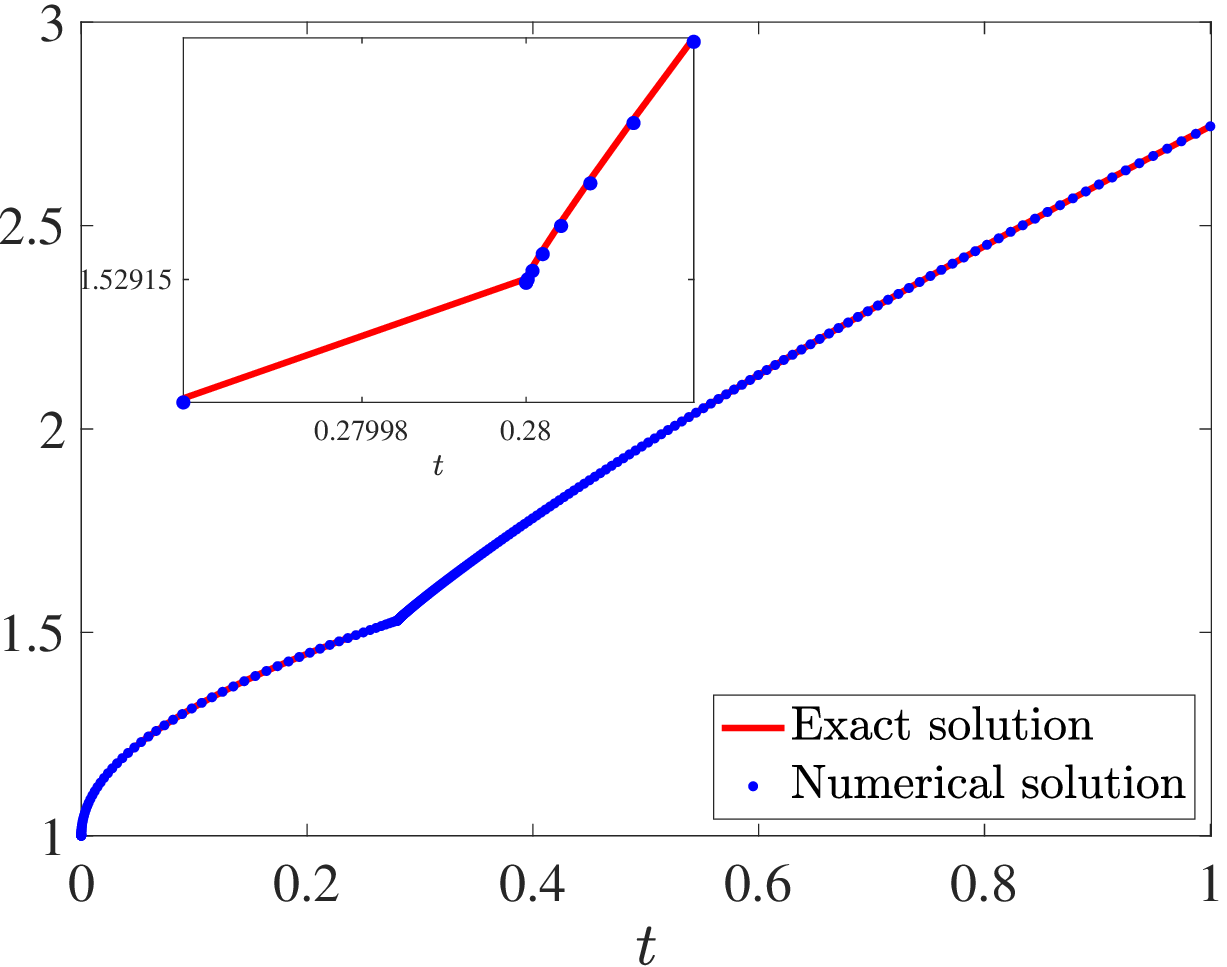}
	\caption{Exact and numerical solution.}\label{fig:ODE_Sol}
\end{subfigure}
\caption{The structure of mesh \eqref{mesh2sig} and two-stage Radau IIA based gCQ solution with $N=256$, $\sigma=0.28$, {$\alpha=0.5$}, $\beta_1=0.5$, $\beta_2=0.9$, $\gamma_1=3/\beta_1$, and $\gamma_2=3/\beta_2$. In the right plot, we zoom in on the exact and numerical solutions around the nonsmooth point $\sigma=0.28$.}\label{fig:ODE_MeshSol}
\end{figure}

\begin{figure}[H]
\centering
\begin{subfigure}[t]{0.49\textwidth}
	\centering
	\includegraphics[width=1\textwidth]{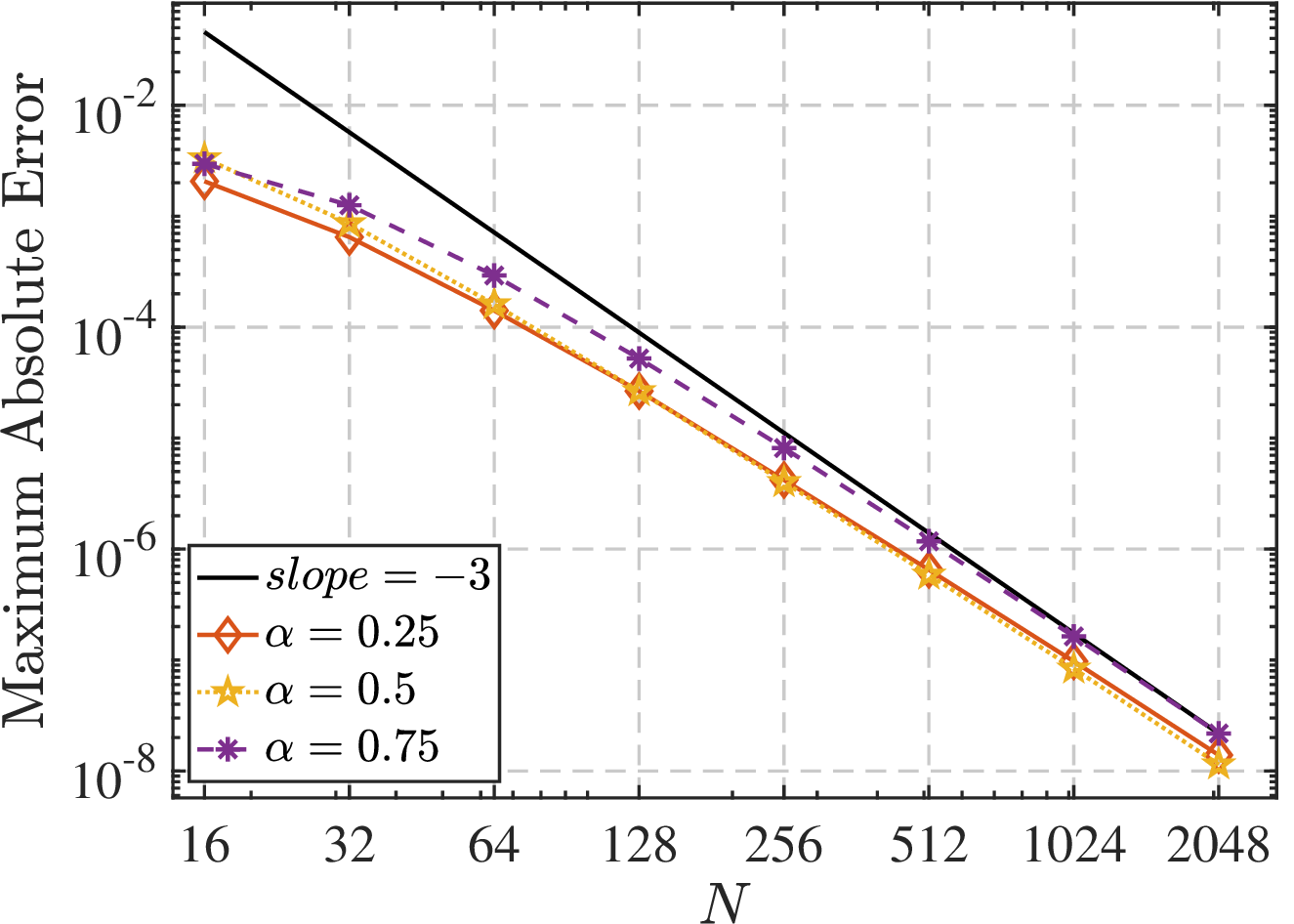}
	\caption{$\sigma=0.28$.}
\end{subfigure}
\hfill
\begin{subfigure}[t]{0.49\textwidth}
	\centering
	\includegraphics[width=1\textwidth]{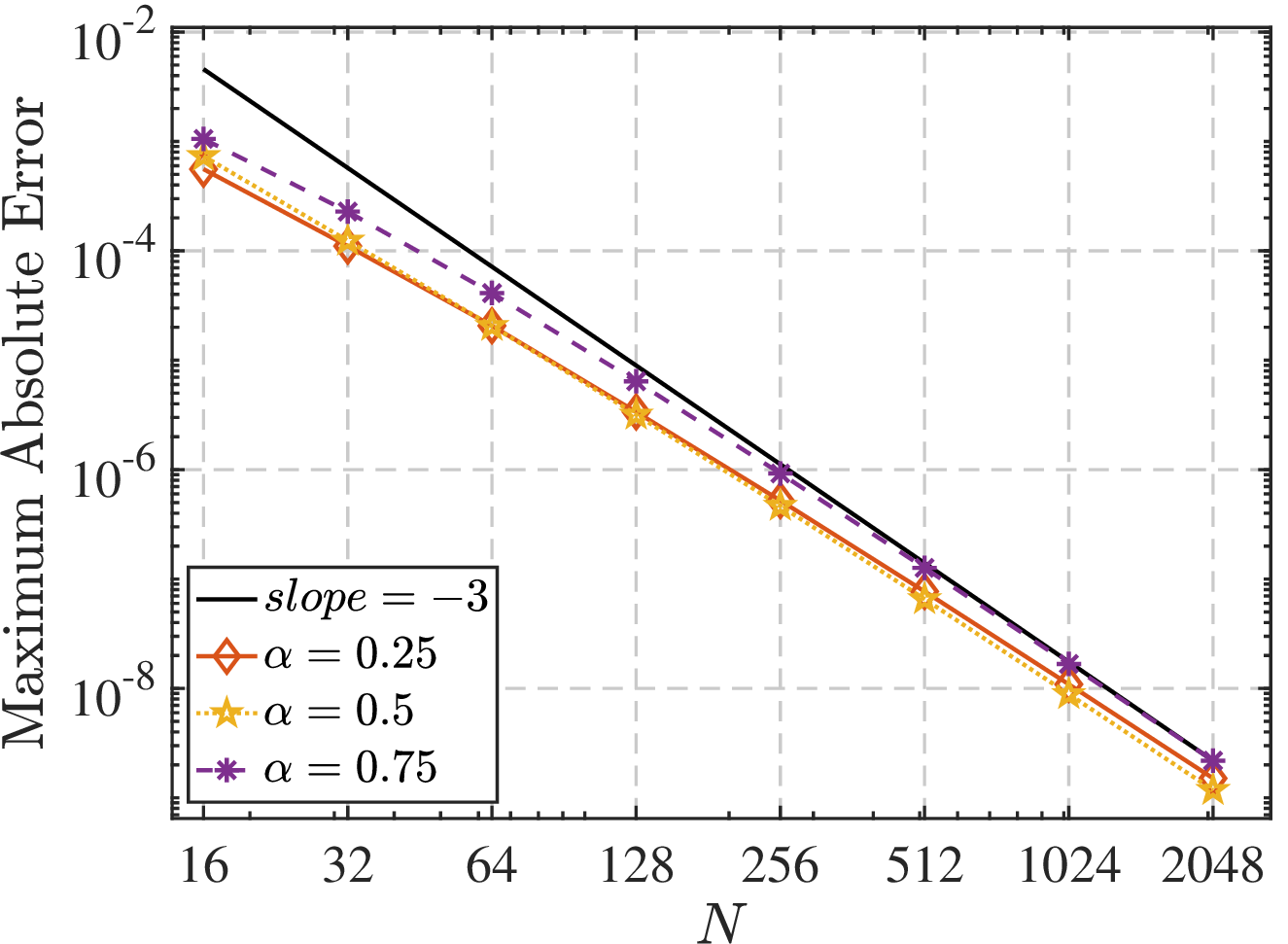}
	\caption{$\sigma=0.72$.}
\end{subfigure}
\caption{Maximum absolute error of the two-stage Radau IIA scheme for Example~\ref{ex_fracODE} with {varying $\alpha = 0.25,\,0.5,\,0.75$}, $\beta_1 = 0.5$, $\beta_2 = 0.9$, $\gamma_1 = 3/\beta_1$, and $\gamma_2 = 3/\beta_2$.}\label{fig:ODE_heavi}
\end{figure}

\begin{example}\label{ex_fracODE}
We consider the fractional ordinary diffusion equation
\begin{equation*}\label{sca_pro}
	D_t^\alpha u(t) + u(t) = f(t), \quad u(0) = 1, \quad 0 < t \leq 1,
\end{equation*}
with the exact solution
\begin{equation*}
	u(t) = 1 + t^{\beta_1} + H(t - \sigma)(t - \sigma)^{\beta_2}, \quad t \ge 0,
\end{equation*}
where \( \beta_1 > -1 \), \( \beta_2 > \alpha \), and \( H(\cdot) \) denotes the Heaviside function. The corresponding source term takes the form
\begin{equation*}
	f(t) = u(t) + \frac{\Gamma(\beta_1 + 1)}{\Gamma(\beta_1 - \alpha + 1)} t^{\beta_1 - \alpha} + \frac{\Gamma(\beta_2 + 1)}{\Gamma(\beta_2 - \alpha + 1)} H(t - \sigma)(t - \sigma)^{\beta_2 - \alpha},
\end{equation*}
exhibiting piecewise smooth behavior.
\end{example}
For this example, we employ a graded mesh with higher density near the nonsmooth points \( t = 0 \) and \( t = \sigma \), defined as
\begin{equation}\label{mesh2sig}
\begin{split}
	&N_1 = \lfloor N \sigma / T \rfloor, \quad 
	\lv(n) = n^{\gamma_1 - 1} (N_1 - n + 1)^{\gamma_2 - 1} N_1^{-\gamma_1 - \gamma_2 + 1}, \quad 1 \le n \le N_1,\\[.5em]
	&t(0) = 0, \quad 
	t(n) = t(n - 1) + \tau(n), \quad 
	\tau(n) = \frac{\sigma}{\| \lv \|_{l^1}} \lv(n), \quad 1 \le n \le N_1,\\[.5em]
	&t(k) = t(N_1) + (T - \sigma) \left( \frac{k - N_1}{N - N_1} \right)^{\gamma_2}, \quad 
	N_1 + 1 \le k \le N-1, \quad 
	t(N) = T.
\end{split}
\end{equation}
In Figure~\ref{fig:ODE_Mesh}, the structure of the mesh defined in~\eqref{mesh2sig} is illustrated, showing that the time step size near \( t = 0 \) and \( t = \sigma \) is of order \( O(N^{-\gamma_1}) \) and \( O(N^{-\gamma_2}) \), respectively, while the maximum step size is of order \( O(N^{-1}) \). The numerical solution, displayed in Figure~\ref{fig:ODE_Sol}, accurately captures the exact solution behavior, even at the nonsmooth point \( t = \sigma \). Moreover, as confirmed in Figure~\ref{fig:ODE_heavi}, the optimal convergence order \( \min\{p,\, q + 1 + \alpha\} \) is asymptotically attained by choosing \( \gamma_1 = \min\{p,\, q + 1 + \alpha\} / \beta_1 \) and \( \gamma_2 = \min\{p,\, q + 1 + \alpha\} / \beta_2 \).

\begin{example}\label{ex_fracPDE}
In this example, we consider the approximation of the two dimensional partial differential equations
\begin{equation*}
	\partial_t^\alpha u-u_{xx}-u_{yy}=f,
\end{equation*}
with initial conditions
\[u(x,y,0)=0,  \quad (x,y)\in  \Omega,\]
and boundary conditions
\[u(x,y,t)=0,\quad (x,y)\in\partial \Omega,\quad 0<t\le 1.\]
The exact solution is given by
\begin{equation*}
	u(x,y,t)=t^\alpha\cos \left(\frac{\pi}{2}x \right) \cos\left(\frac{\pi}{2}y \right),
\end{equation*}
and the source term is
\[f(x,y,t)=\left(\Gamma(\alpha+1)+\frac{\pi^2}{2} t^\alpha\right)\cos \left(\frac{\pi}{2}x \right) \cos\left(\frac{\pi}{2}y \right).\]

\end{example}
In our numerical experiments, we consider the domain \( \Omega = (-1,1)^2 \) and employ a fourth-order compact finite difference scheme \cite{ZhuJuZhao} with uniform grid spacing \( \Delta_x = \Delta_y = \frac{1}{128} \) for spatial discretization. {We apply the algorithm described in Subsection~\ref{subsec:subdiff}, with $\mathcal{A}=A_h$ the discrete Laplacian obtained by applying this spatial discretization scheme. For a given spatial grid, the analysis in \cite[Section 4]{GuoLo} applies to the semi-discrete problem in space, since $A_h$ is a self-adjoint operator and all its eigenvalues are negative. We do not address here the derivation of full error estimates in space and time and only test  the sharpness of our error estimates in time.}

\begin{figure}[H]
\centering
\begin{subfigure}[t]{0.49\textwidth}
	\centering
	\includegraphics[width=1\textwidth]{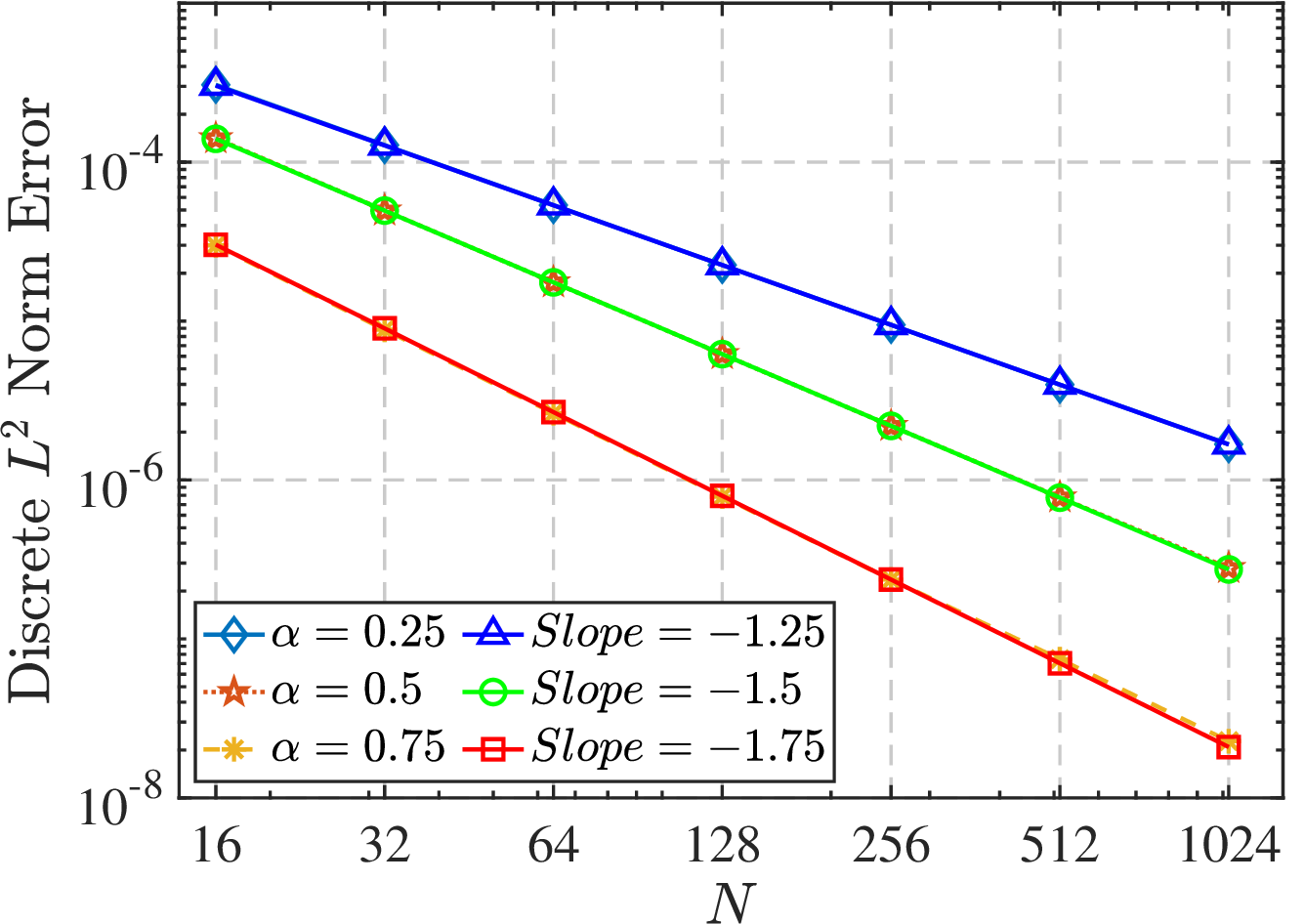}
	\caption{Corrected BDF3 CQ \cite{JinLiZh}.}
\end{subfigure}
\hfill
\begin{subfigure}[t]{0.49\textwidth}
	\centering
	\includegraphics[width=1\textwidth]{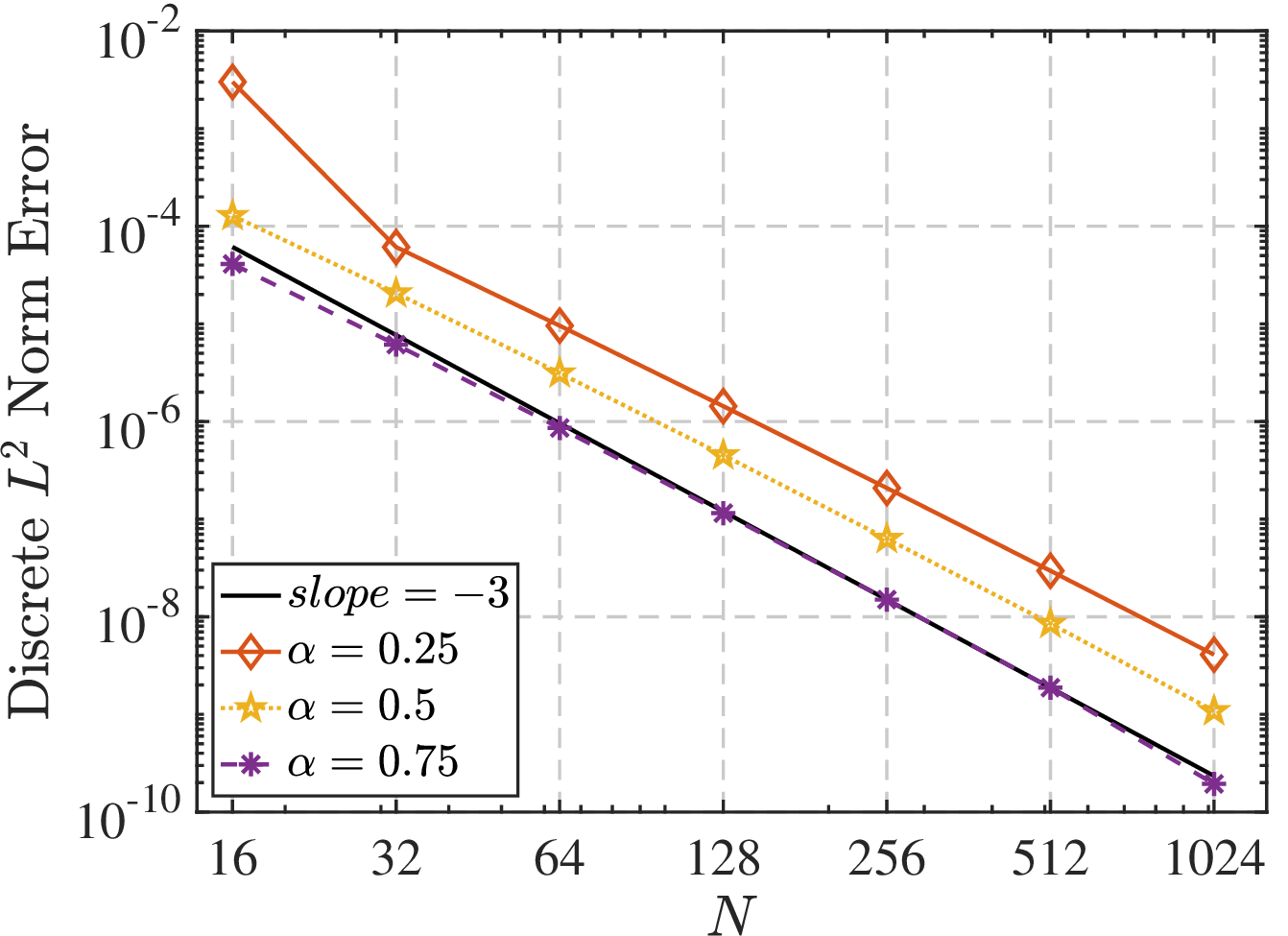}
	\caption{Two-stage Radau IIA gCQ   on \eqref{gmesh} with  $\gamma=\frac{3}{\alpha}$.}
\end{subfigure}
\caption{Comparison of discrete $L^2$-norm errors for Example~\ref{ex_fracPDE}  {with varying $\alpha = 0.25,\,0.5,\,0.75$} at the final time $t_N = 1$.}
\label{fig:fracdiff_errtN}
\end{figure}

\begin{figure}[H]
\centering
\begin{subfigure}[t]{0.49\textwidth}
	\centering
	\includegraphics[width=1\textwidth]{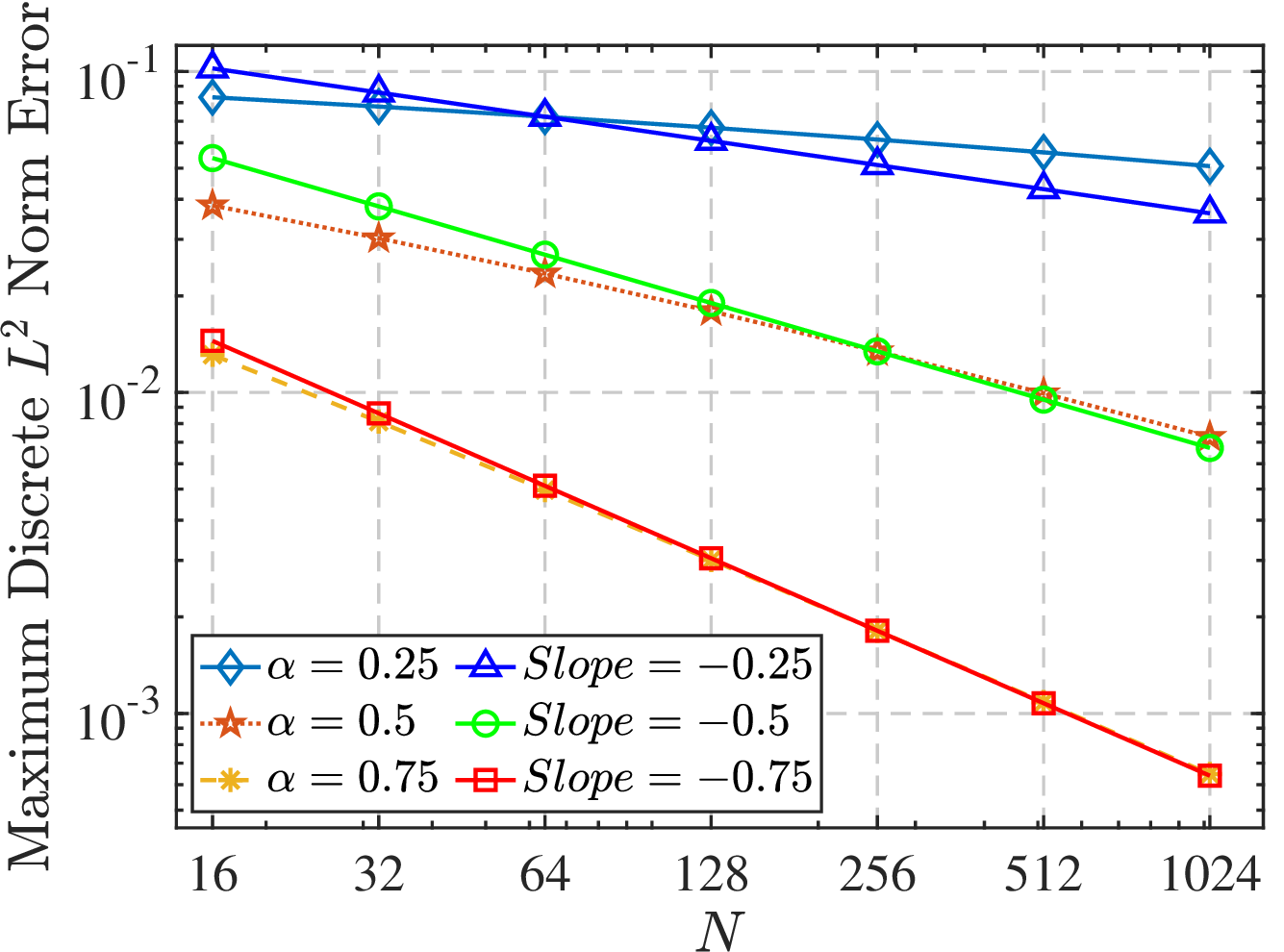}
	\caption{Corrected BDF3 CQ \cite{JinLiZh}.}
\end{subfigure}
\hfill
\begin{subfigure}[t]{0.49\textwidth}
	\centering
	\includegraphics[width=1\textwidth]{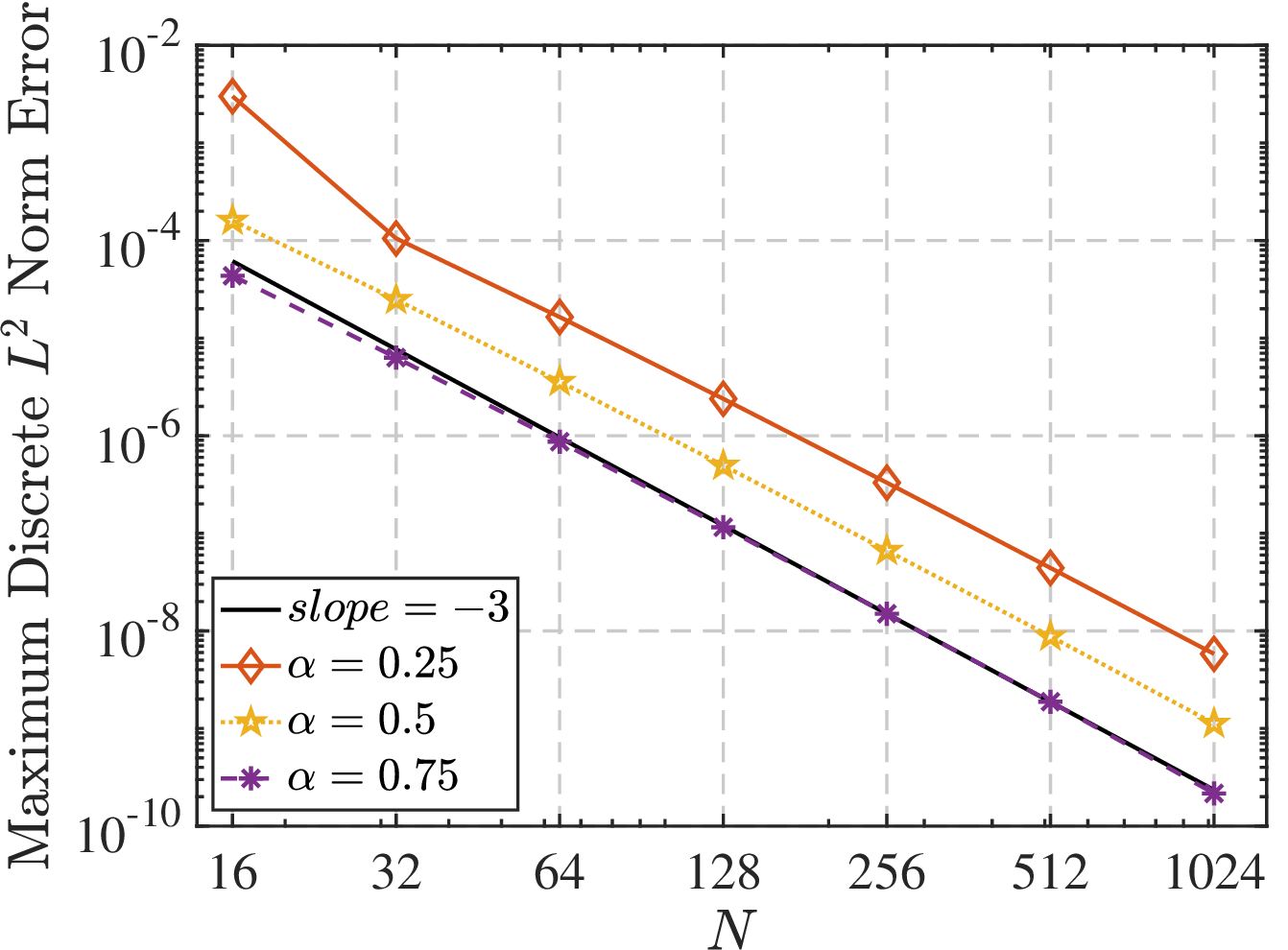}
	\caption{Two-stage Radau IIA gCQ   on \eqref{gmesh} with  $\gamma=\frac{3}{\alpha}$.}
\end{subfigure}
\caption{Comparisons of maximum discrete $L^2$ norm errors for Example~\ref{ex_fracPDE}  {with varying $\alpha = 0.25,\,0.5,\,0.75$}.}\label{fig:fracdiff_errMax}
\end{figure}

In Figures~\ref{fig:fracdiff_errtN} and \ref{fig:fracdiff_errMax}, we compare the two-stage Radau IIA based gCQ method with the corrected third-order backward differentiation method from \cite{JinLiZh}. These results show that for solutions with limited temporal regularity, the full convergence order is not preserved when using the corrected scheme in \cite{JinLiZh}, whereas the Runge--Kutta based gCQ method on graded meshes recovers optimal convergence both at the final time point and across all evaluated time steps. Notice that the optimal grading parameter in this case is $\gamma=\frac{3}{\alpha}$, which corresponds to the case $\beta=0$, $p=3$, $q=2$, in Theorem~\ref{thm:conv_gCQrkGrad}. According to the analysis in \cite[Section 4]{GuoLo},  the optimal time mesh is determined by the behaviour of the source term $f$ close to 0, rather than the behaviour of the exact solution $u$, and the numerical experiments confirm the theory.  
\begin{example}\label{ex_Westervelt}
Consider the following  Westervelt equation \cite{BaBanPtas24}
\begin{equation}\label{eq:West}
	\begin{cases}
		(1-2\kappa u)\partial_t^2 u - u_{xx} - k \ast \partial_t u_{xx} = 2\kappa(\partial_t u)^2+f, & x \in (-8,8), \quad t \in (0,2], \\[.5em]
		u(-8,t) = u(8,t) = 0, & t \in (0,2], \\[.5em]
		u(x,0) = \frac{\e^{-x^2}}{2}, \quad \partial_t u(x,0) = 0,& x \in (-8,8),
	\end{cases}
\end{equation}
where
\[k(t)=\frac{t^{\alpha-1}}{\Gamma(\alpha)}\e^{-t},\quad f(x,t)=\left(1+\log(t)\right)\e^{-x^2},\quad t\in(0,2], \quad x\in (-8,8).\]
\end{example}
Applying the Runge--Kutta based gCQ method for temporal discretization and the fourth-order compact difference scheme \cite{ZhuJuZhao} for spatial discretization with $J$ uniform grid subdivisions, the fully discrete scheme for solving \eqref{eq:West} is given by
\begin{equation}\label{Weseq_sch}
\begin{cases}
	\left(1 - 2\kappa \Uv_n^h\right)\left[(\partial_t^\Delta)^2 \Uv^h\right]_n 
	- \mathbf{L} \Uv_n^h 
	- \left[K(\partial_t^\Delta)\left(\partial_t^\Delta \mathbf{L} \Uv^h\right)\right]_n 
	= 2\kappa \left[\partial_t^\Delta \Uv^h\right]_n^2+\fv_n^h,\\[1em]
		\Uv_0^h = \left(\frac{1}{2}\e^{-x_j^2}\right)_{j=1}^{J-1},\quad  
		\partial_t^\Delta \Uv_0^h = \mathbf{0},\quad  
	x_j=-8+jh,\quad h=16/J,\quad 1 \le j \le J - 1.
\end{cases}
\end{equation}
Here, \(\fv_n^h=\bigl(f(x_j,t_{n-1}+c_i\tau_n)\bigr)_{i=1,\ldots,s;\,j=1,\ldots,J-1}\), and the matrix \( \mathbf{L} \) is derived by applying the compact difference method \cite{ZhuJuZhao} to approximate \( u_{xx} \). The operator \( \partial_t^\Delta \) is defined in \eqref{gCQ_1stderv}. The term \( \left[(\partial_t^\Delta)^2 \Uv^h\right]_n \) is obtained by applying \( \partial_t^\Delta \) to \( \left[\partial_t^\Delta \Uv^h\right]_n \), which is given by
\[\left[(\partial_t^\Delta)^2\Uv^h \right]_n= (\tau_n \Av)^{-1}\left(\left[\partial_t^\Delta\Uv^h\right]_n-\ev_s^\top\left[\partial_t^\Delta\Uv^h\right]_{n-1}\bone\right).\]
In the implementation, the fast algorithm shown in Section~\ref{sec:algfi} is employed to approximate  $\left[K(\partial_t^\Delta)\left(\partial_t^\Delta \mathbf{L} \Uv^h\right)\right]_n $,  and  the nonlinear system \eqref{Weseq_sch} is solved by the fixed point iteration with a tolerance of $ 10^{-8}$.  With spatial mesh size \( h = 1/16\), the maximum discrete \( L^2 \) norm error is measured by  
\[
\max_{1 \leq n \leq N} \Vert \widetilde{u}_{2n}^h - u_n^h \Vert_{h},
\]  
where  
\[
\Vert v \Vert_h^2 = h \sum_{j=1}^{J-1} (v^j)^2
\]  
denotes the discrete \( L^2 \) norm, and \( \widetilde{u}_{2n}^h \) represents the numerical solution computed on the refined time mesh \eqref{gmesh_finer}.

\begin{figure}[H]
\centering
\begin{subfigure}[t]{0.49\textwidth}
	\centering
	\includegraphics[width=1\textwidth]{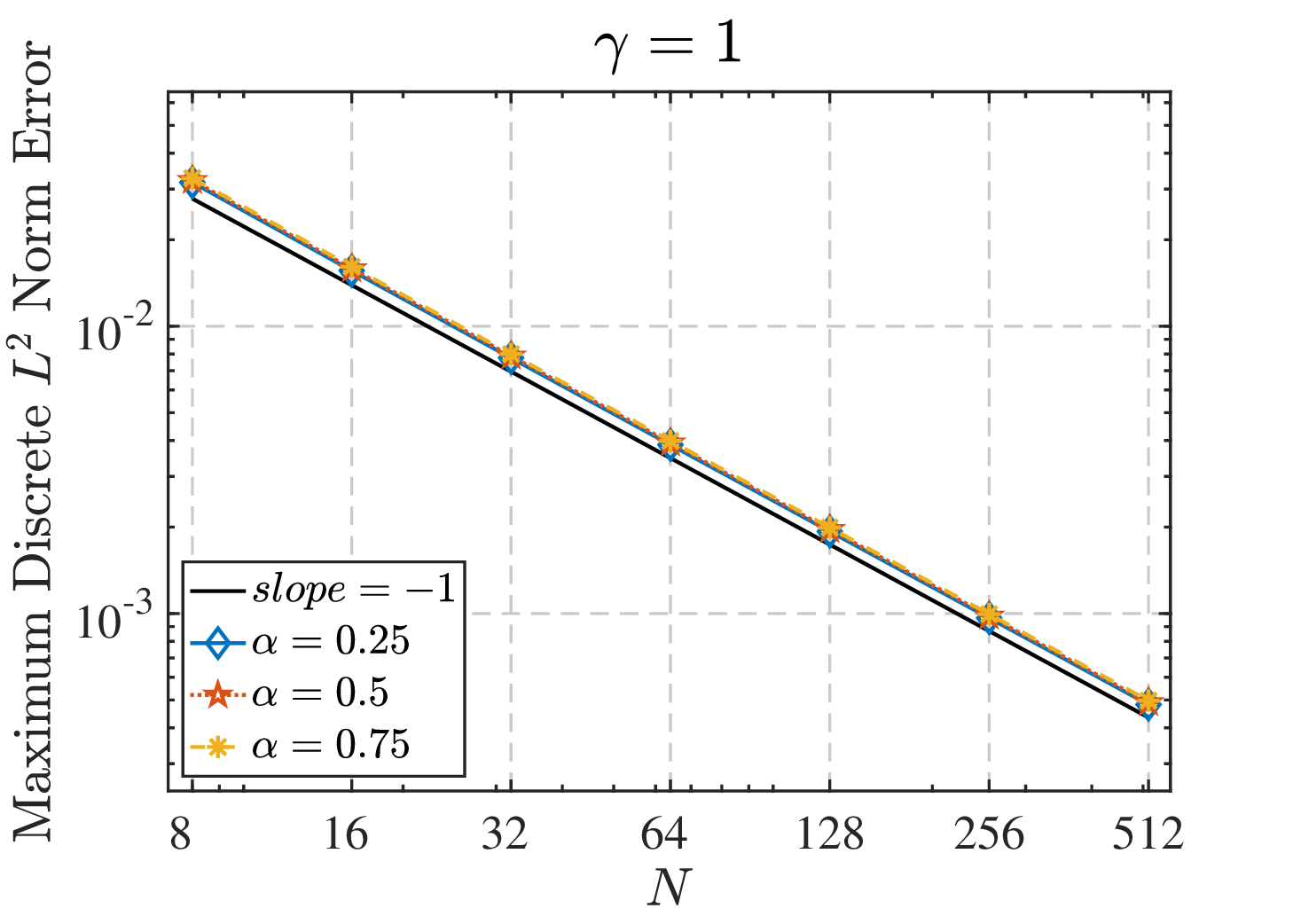}
\end{subfigure}
\hfill
\begin{subfigure}[t]{0.49\textwidth}
	\centering
	\includegraphics[width=1\textwidth]{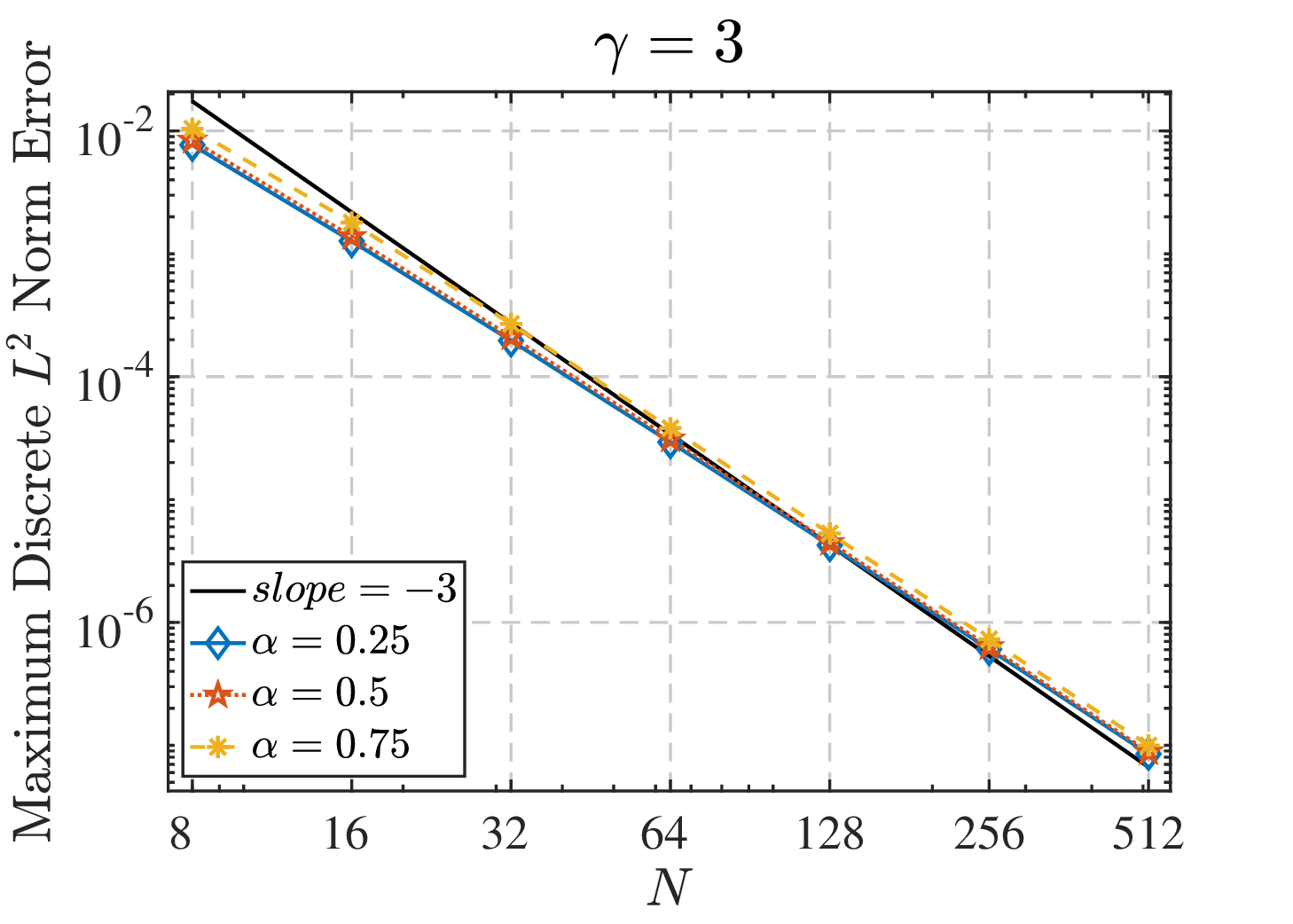}
\end{subfigure}
\caption{Maximum discrete $L^2$ norm errors of the  two-stage Radau IIA Runge--Kutta based gCQ  method  for Example~\ref{ex_Westervelt} with  {varying $\alpha = 0.25,\,0.5,\,0.75$} and $\kappa=0$.}\label{fig:WestV_kap0}
\end{figure}

\begin{figure}[H]
\centering
\begin{subfigure}[t]{0.49\textwidth}
	\centering
	\includegraphics[width=1\textwidth]{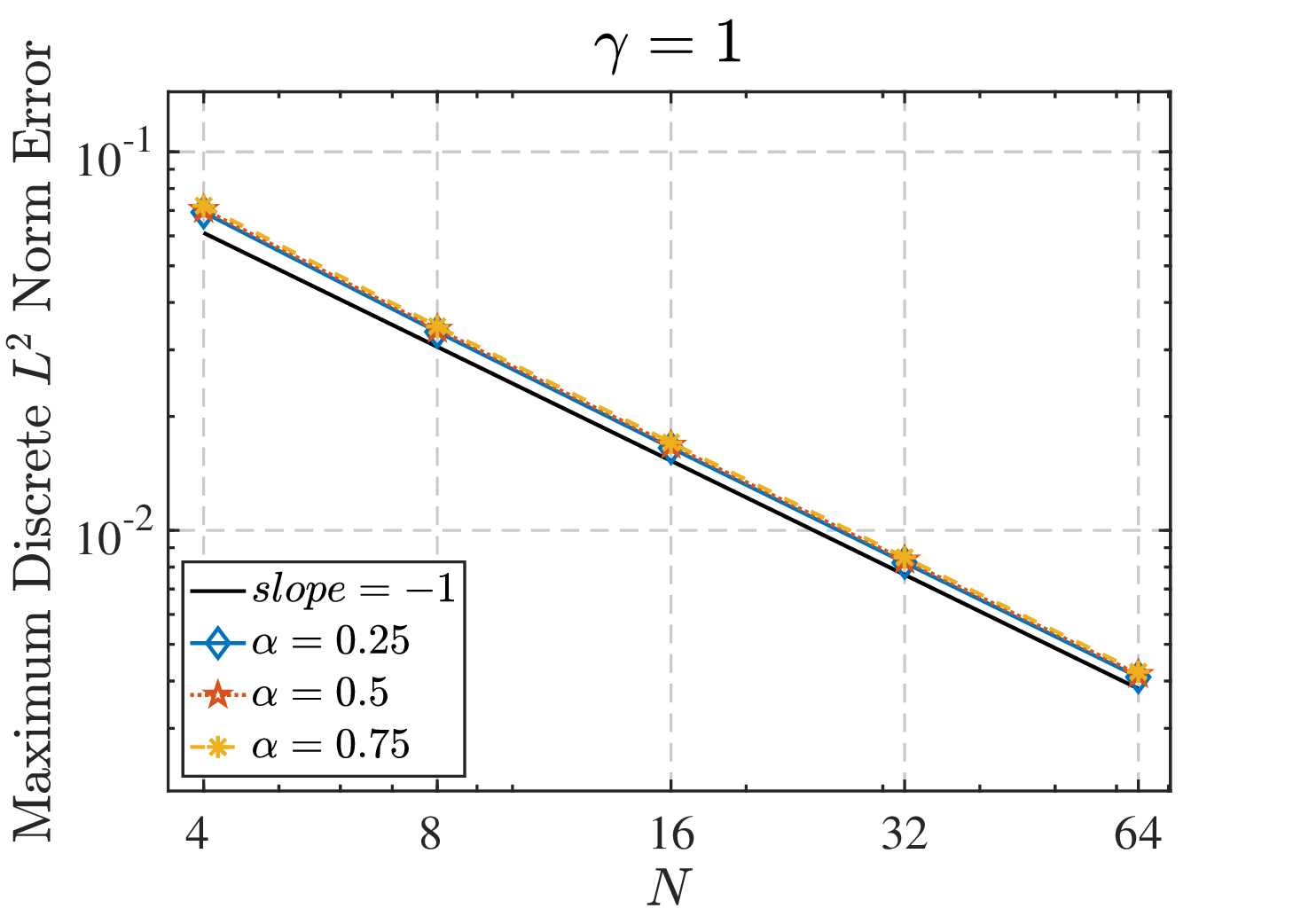}
\end{subfigure}
\hfill
\begin{subfigure}[t]{0.49\textwidth}
	\centering
	\includegraphics[width=1\textwidth]{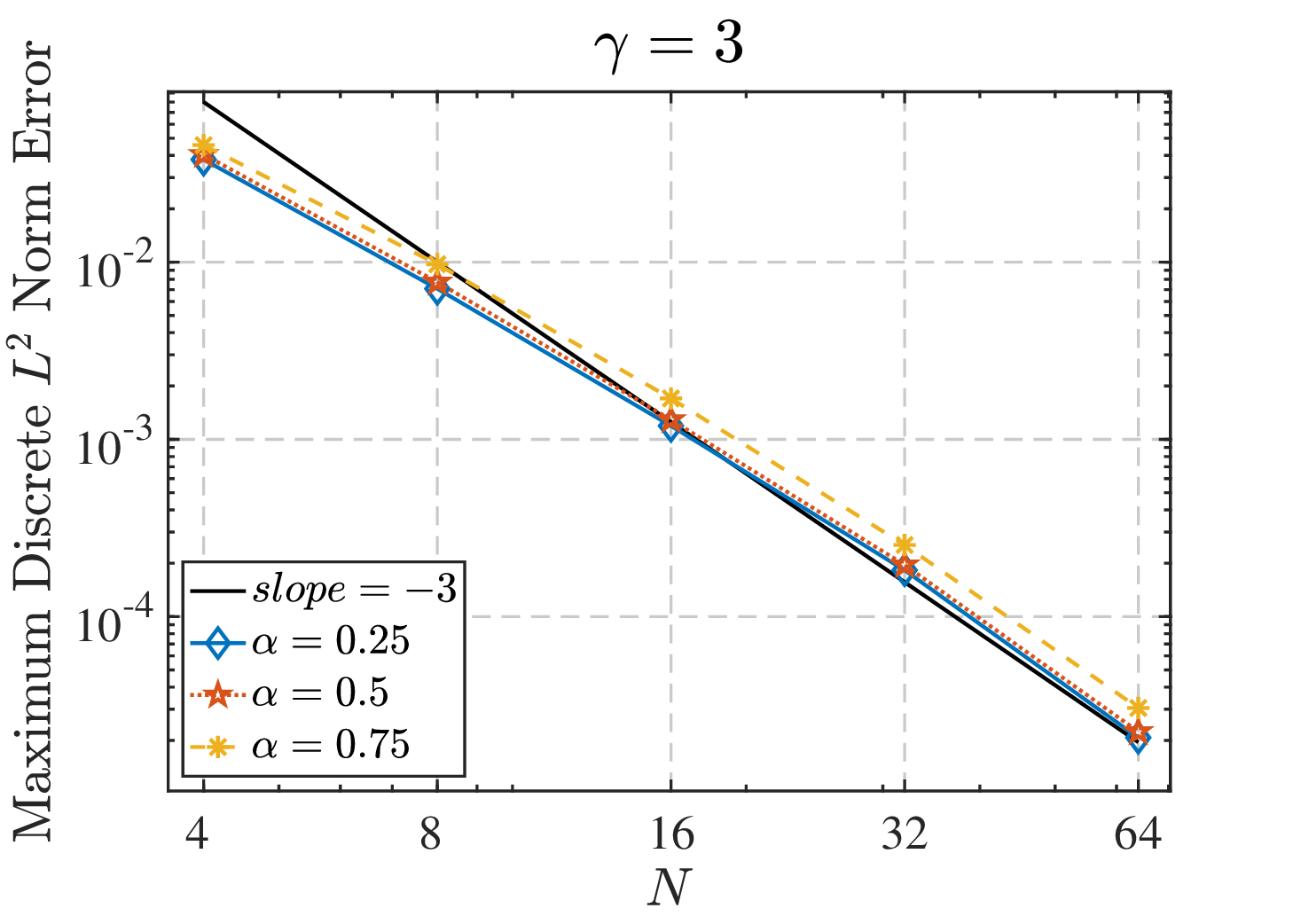}
\end{subfigure}
\caption{Maximum discrete $L^2$ norm errors of the  two-stage Radau IIA Runge--Kutta based gCQ  method  for Example~\ref{ex_Westervelt}  with   {varying $\alpha = 0.25,\,0.5,\,0.75$} and $\kappa=0.09$.}\label{fig:WestV_kap009}
\end{figure}

\begin{figure}[H]
\centering
\begin{subfigure}[t]{0.49\textwidth}
	\centering
	\includegraphics[width=1\textwidth]{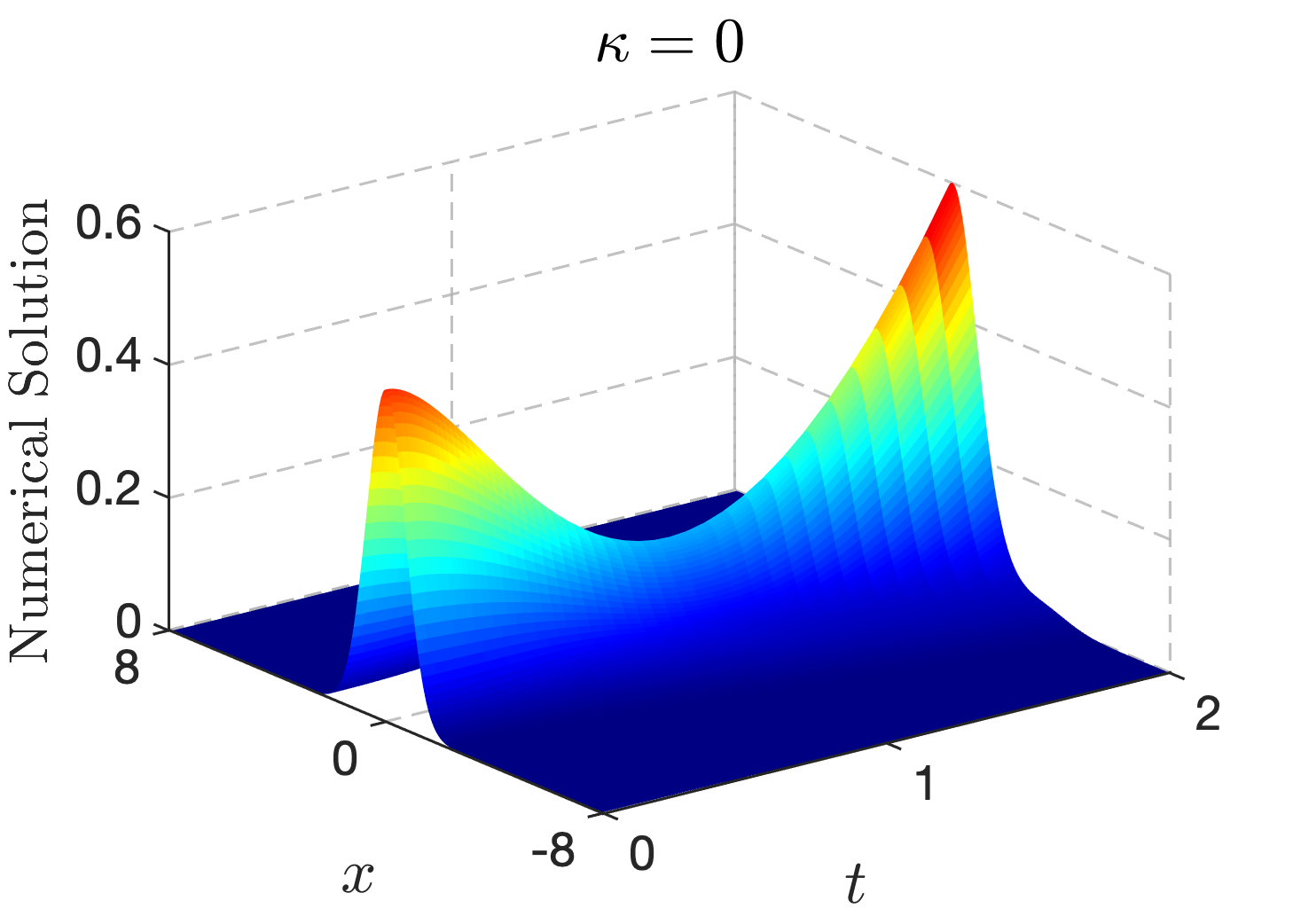}
\end{subfigure}
\hfill
\begin{subfigure}[t]{0.49\textwidth}
	\centering
	\includegraphics[width=1\textwidth]{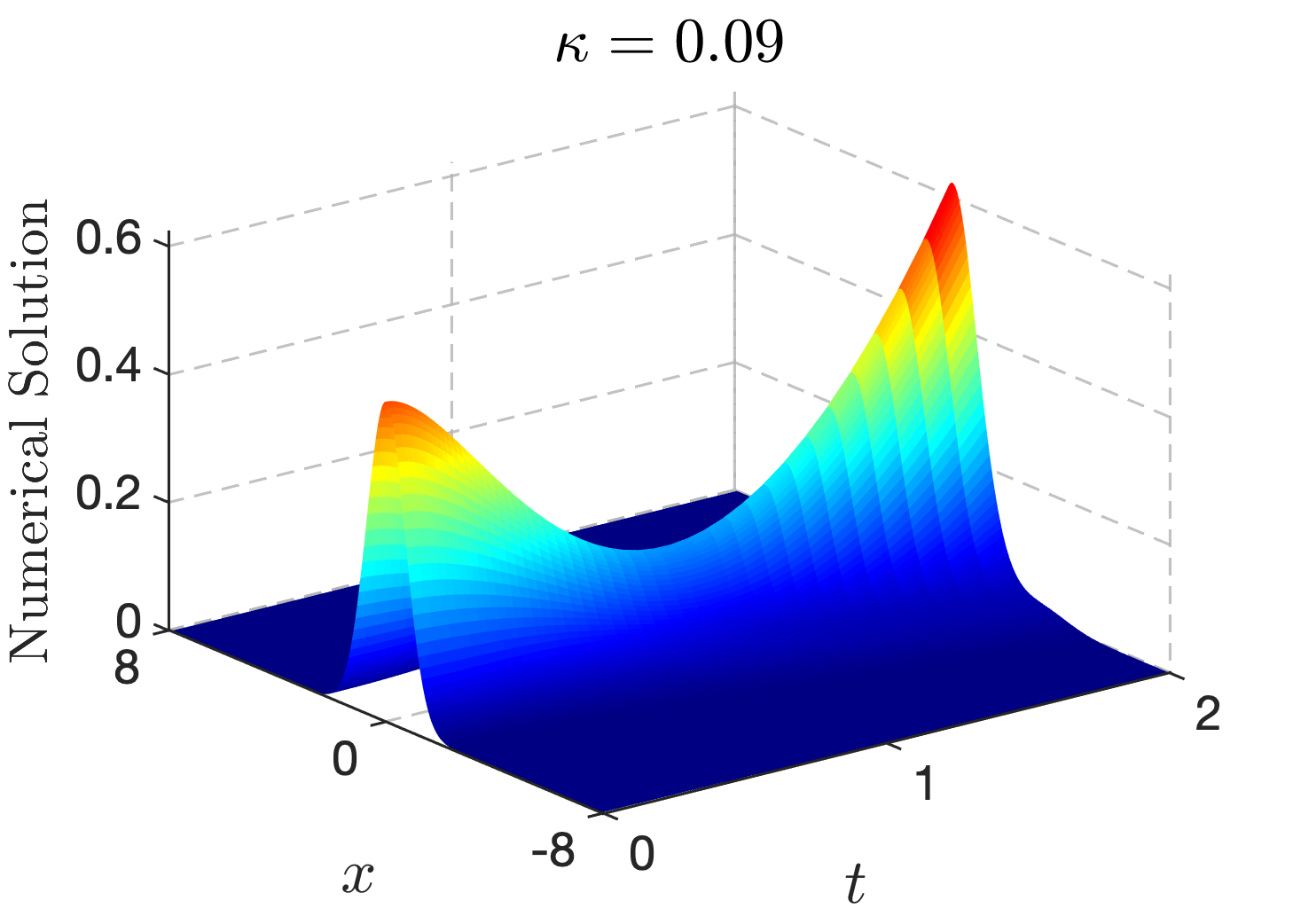}
\end{subfigure}
\caption{Numerical solution of the two-stage Radau IIA Runge--Kutta based gCQ method for Example~\ref{ex_Westervelt} with  $\alpha=0.5$, $\gamma=3$ and $N=64$.}
\label{fig:WestV_kap009NumSol}
\end{figure}

\begin{figure}[H]
\centering
\begin{subfigure}[t]{0.49\textwidth}
	\centering
	\includegraphics[width=1\textwidth]{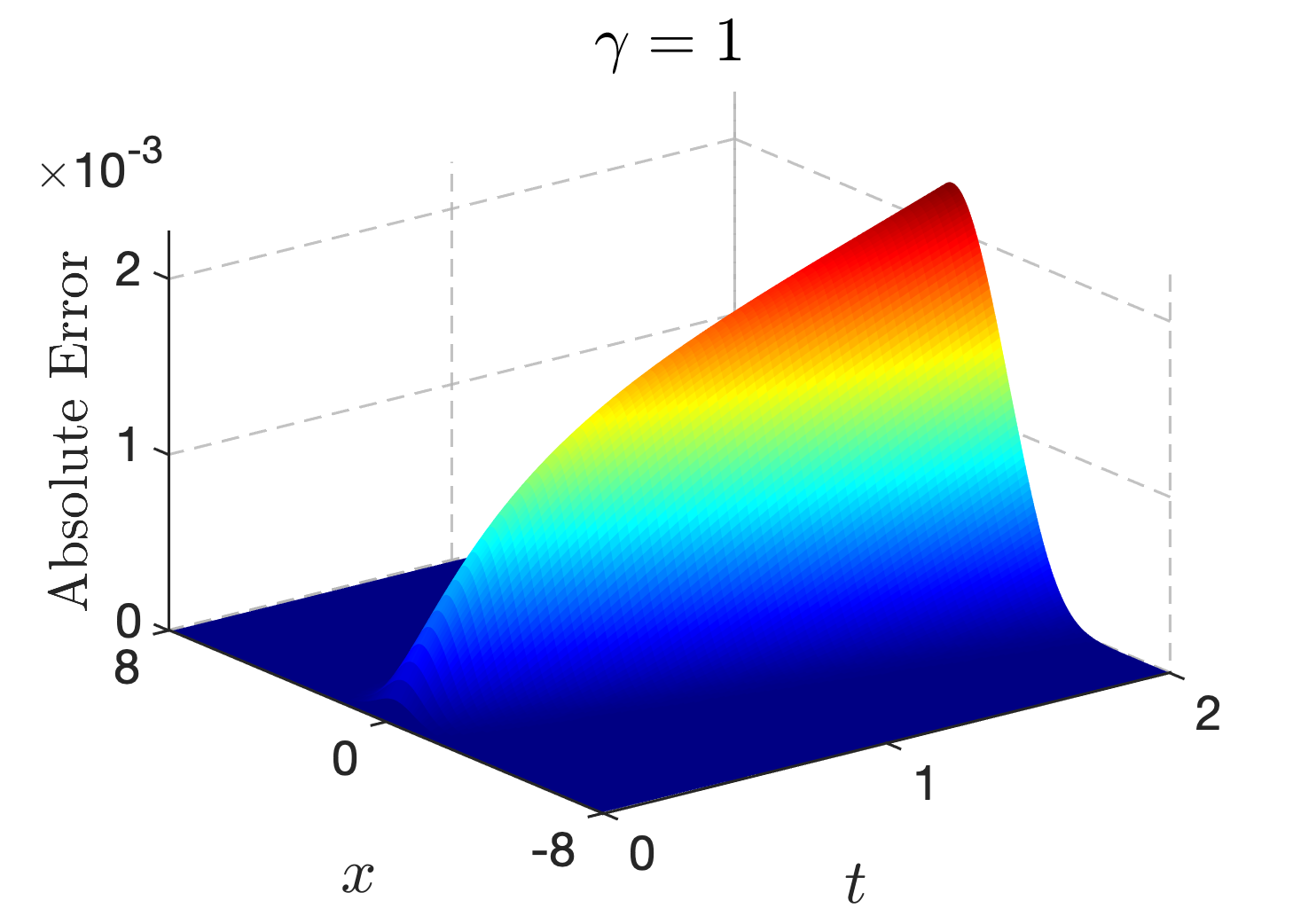}
\end{subfigure}
\hfill
\begin{subfigure}[t]{0.49\textwidth}
	\centering
	\includegraphics[width=1\textwidth]{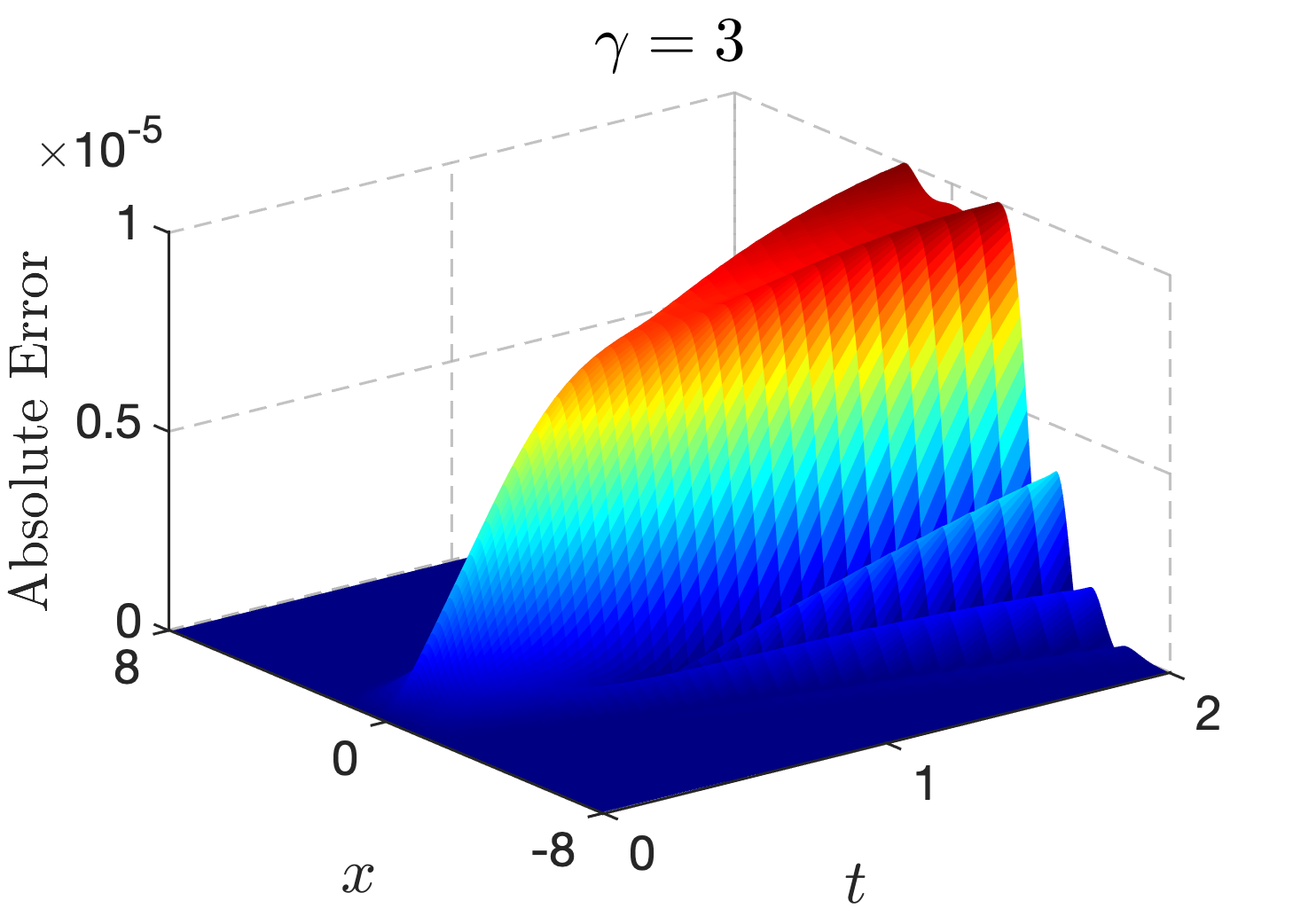}
\end{subfigure}
\caption{Absolute error of the two-stage Radau IIA Runge--Kutta based gCQ method for Example~\ref{ex_Westervelt} with  $\alpha=0.5$, $\kappa=0.09$ and $N=64$.}
\label{fig:WestV_kap009AbsErr}
\end{figure}

From the numerical results in Figures~\ref{fig:WestV_kap0} and~\ref{fig:WestV_kap009}, we observe that employing a graded mesh \eqref{gmesh} with \(\gamma = 3\) yields full third-order convergence asymptotically. In contrast, when using the uniform mesh, that is, \eqref{gmesh} with \(\gamma = 1\), the convergence rate is limited to first order.
The evolution of the numerical solution is shown in Figure~\ref{fig:WestV_kap009NumSol}. It can be observed that the case with \( \kappa = 0.09 \) attains slightly larger solution values at the final time \( T = 2 \) on the domain \( (-8, 8) \). Moreover, the temporal evolution of the numerical solution indicates that the maximum values in the \( x \)-direction for both \( \kappa = 0 \) and \( \kappa = 0.09 \) exhibit the same qualitative behavior: decreasing over the time interval \( (0,1] \) and increasing thereafter.
Figure~\ref{fig:WestV_kap009AbsErr} further demonstrates that the Runge--Kutta based gCQ method, when implemented on the graded mesh \eqref{gmesh} with $\gamma = 3$, yields smaller maximum errors near the final time than those obtained using a uniform mesh. Increasing the accuracy near the origin reduces the global error, even though its maximum is attained at the final time $T=2$; see Figure~\ref{fig:WestV_kap009AbsErr}.
These results demonstrate the effectiveness of the Runge--Kutta based gCQ method in solving equations involving convolution integrals.



\section*{Acknowledgements}
The second author has been supported by ``Proyecto 16 - proyecto G Plan Propio'' of the University of Malaga, by grant PID2022-137637NB-C21 funded by MICIU/AEI/10.13039/501100011033 and ERDF/EU and by grant PPRO-FQM216-G-2023 funded by Consejería de Universidad, Investigación e Innovación and by ERDF Andalusia Program 2021-2027. 

\section*{Data availability statement}
The codes implementing the algorithms discussed in this article are publicly available at: \url{https://doi.org/10.5281/zenodo.19482892}. No other data are associated with the manuscript.

\appendix
\section{Identities}\label{sec:regularity}
This appendix collects essential mathematical results used in our analysis. We begin with the definition of the Mittag-Leffler function, followed by key identities and estimates for solutions to \eqref{ode}. The proofs of these results follow the techniques in \cite[Lemma 6]{GuoLo}.
\begin{definition}\label{def:ml}
[Mittag-Leffler function] For $\alpha>0$, $\beta\in\mathbb{R}$, the two parameter Mittag-Leffler function $E_{\alpha,\beta}(z)$ is defined by
$$E_{\alpha,\beta}(z)=\sum\limits_{\ell=0}^{\infty}\frac{z^\ell}{\Gamma(\alpha \ell+\beta)},\quad z\in\mathbb{C}.$$
\end{definition}
\begin{lemma}\label{lem:identities}
For  \( \beta > -1 \), \( m \in \mathbb{N} \), it holds
\begin{enumerate}
	\item For $0<\alpha<1$, $\zeta>0$, we have
	\begin{equation}\label{int_Betafun}
		\int_0^\infty x^{-\alpha} \frac{1}{1+\zeta x}\,dx = B(\alpha,1-\alpha)\, \zeta^{\alpha-1},
	\end{equation}
	where \( B\) is the Beta function.
	\item \cite[Theorem 1.6]{Podlubny} For \(  {0< \alpha < 2 }\), \( \beta \in \mathbb{R} \), there exists a constant \( C > 0 \) such that
	\begin{equation*}\label{bound-absml-neg}
		|E_{\alpha,\beta}(-x)| \le \frac{C}{1+x}, \qquad x\ge 0.
	\end{equation*}
	\item 
	The solution to \eqref{ode} with \( f(t) = t^\beta\vv \) is
	{
		\begin{align*}
			\label{y}
			y(x,t) &= \Gamma(\beta+1) \sum_{\ell=0}^{\infty} \frac{(-x)^{\ell}}{\Gamma( \ell+\beta +2)}t^{\ell+\beta+1} \vv = \Gamma(\beta+1) t^{\beta+1} E_{1,\beta+2}(-xt) \vv, \quad \forall\beta>-1.
	\end{align*}}
	In particular, for $\beta\in\bN$, it can be rewritten as
	\begin{equation*}\label{y_int}
		y(x,t) =	\beta!(-x)^{-(\beta+1)}\left(
		\e^{-xt}-\sum_{\ell=0}^{\beta}  \frac{(-xt)^{\ell}}{\ell!}\right) \vv.
	\end{equation*}
	\item The derivatives of the solution to \eqref{ode} with \( f(t) = t^\beta\vv \) with $\beta>-1$ can be expressed explicitly. For \( \beta \notin \mathbb{N} \), the expression is given by
	\begin{eqnarray*}\label{ydiff_nonint}
		\frac{\partial^m y}{\partial t^m}(x, t) = \Gamma(\beta + 1)\, t^{\beta - m + 1} E_{1, \beta - m + 2}(-x t) \vv,
	\end{eqnarray*}
	while for \( \beta \in \mathbb{N} \), it takes the form
	\begin{equation*}\label{ydiff_intg}
		\frac{\partial^m y}{\partial t^m}(x, t) = 
		\begin{cases}
			\beta!\, t^{\beta - m + 1} E_{1, \beta - m + 2}(-x t) \vv, & \text{if } m \le \beta, \\[.5em]
			\beta!\, (-x)^{-\beta + m - 1} e^{-x t} \vv, & \text{if } m \ge \beta + 1.
		\end{cases}
	\end{equation*}
	In both cases, the derivatives satisfy
	\begin{align}\label{ydiff_bnd}
		\left\Vert  \frac{\partial^m y}{\partial t^m}(x, t) \right\Vert_{\mathcal{X}}  \le \Gamma(\beta + 1)\, t^{\beta - m + 1} \frac{C}{1 + x t},
	\end{align}
	where the constant \( C \) depends on $\Vert \vv\Vert_{\mathcal{X}}$, \( \beta \) and \( m \). In particular, if \( \beta \in \mathbb{N} \) and \( m \ge \beta + 1 \), then \( C = (m - \beta)! \Vert \vv\Vert_{\mathcal{X}}\).
	\item The $m$-th derivative of the solution to \eqref{ode} with \( f(t) = H(t-\eta)(t-\eta)^\beta \vv \), where \( \beta \ge 0 \), \( m \le \beta + 2 \), and \( \eta \ge 0 \), can be expressed explicitly. For \( \beta \notin \mathbb{N} \), the expression is given by
	\begin{align*}
		\frac{\partial^m y}{\partial t^m}(x, t) &=\Gamma(\beta + 1) H(t-\eta) (t-\eta)^{\beta - m + 1} E_{1, \beta - m + 2}(-x (t-\eta))\vv\\[.5em]
		&\hfill+\Gamma(\beta + 1)\delta(t-\eta) (t-\eta)^{\beta - m + 2} E_{1, \beta - m + 3}(-x(t-\eta)) \vv,
	\end{align*}
	while for \( \beta \in \mathbb{N} \), it takes the form
	\begin{equation*}
		\frac{\partial^m y}{\partial t^m}(x, t) = 
		\begin{cases}
			\beta!\, H(t-\eta)(t-\eta)^{\beta - m + 1} E_{1, \beta - m + 2}(-x (t-\eta)) \vv, & \text{if } m \le \beta, \\[.5em]
			\beta!\, H(t-\eta) e^{-x (t-\eta)}\vv, & \text{if } m = \beta + 1\\[.5em]
			\beta!\, \left(-H(t-\eta)x e^{-x (t-\eta)}+\delta(t-\eta) e^{-x (t-\eta)} \right) \vv, & \text{if } m = \beta + 2.
		\end{cases}
	\end{equation*}
	In both cases, the derivatives satisfy, for \(t\ge\eta\),
	\begin{align}\label{ydiff_bndpeano}
		\left\Vert  \frac{\partial^m y}{\partial t^m}(x, t) \right\Vert_{\mathcal{X}}  \le C\Gamma(\beta+1) \frac{H(t-\eta) (t-\eta)^{\beta - m + 1}+\delta(t-\eta) (t-\eta)^{\beta - m + 2}}{1 + x(t-\eta)}.
	\end{align}
\end{enumerate}
\end{lemma}
\section{Estimates on the graded mesh}
Some estimates on the graded mesh defined in \eqref{gmesh} are provided below.
\begin{lemma}\label{lem:grad_par}
Let $m > 0$, and consider the time mesh defined in \eqref{gmesh} with step sizes \(\tau_j\). Then for any \(\sigma_j \in [t_{j-1}, t_j]\), \(2 \le j \le n\), the following estimates hold
\begin{align}
	\label{grad_prodbnd}
	\tau_j^{\nu_1+m} \sigma_j^{\nu_2- m} &\le C  {T^{ \nu_1+\nu_2}}
	\begin{cases}
		N^{-\gamma (\nu_1+\nu_2)}, & \gamma (\nu_1+\nu_2) < m, \\[0.5em]
		N^{-m}, & \gamma(\nu_1+\nu_2)\ge m,
	\end{cases}
\end{align}
and moreover,
\begin{align}
	\label{grad_sumprodbnd}
	\sum_{j=2}^n \tau_j^{m} \sigma_j^{\nu - m} &\le C  {T^{ \nu}}
	\begin{cases}
		N^{-\gamma \nu}, & \gamma \nu < m -1, \\[0.5em]
		N^{-m+1} \log(n), & \gamma \nu = m - 1, \\[0.5em]
		N^{-m+1}, & \gamma \nu > m - 1,
	\end{cases}
\end{align}
where \( C > 0 \) is a constant independent of the mesh.
\end{lemma}

\begin{proof}
From the graded mesh definition \eqref{gmesh}, it follows that
\[
\tau_j \leq \gamma j^{\gamma - 1} \tau^\gamma,
\]
and
\[
\sigma_j^{\nu_2 - m} \leq \max\{(j-1)^{\gamma(\nu_2 - m)}, j^{\gamma(\nu_2 - m)}\} \tau^{\gamma(\nu_2 - m)}.
\]
Hence,
\[
\tau_j^{\nu_1+m} \sigma_j^{\nu_2 - m} \leq \gamma^{\nu_1 + m} j^{(\nu_1 + m)(\gamma - 1)} \max\{(j-1)^{\gamma(\nu_2 - m)}, j^{\gamma(\nu_2 - m)}\} \tau^{\gamma(\nu_1 + \nu_2)}.
\]
Using the inequality \( j/2 \leq j - 1 \leq j \) for \( j \geq 2 \), we obtain
\begin{equation*}
	\tau_j^{\nu_1+m} \sigma_j^{\nu_2- m}  \leq \max\{1, 2^{\gamma(m - \nu_2)}\} \gamma^{\nu_1+m}j^{\gamma (\nu_1+\nu_2) - m} \tau^{\gamma (\nu_1+\nu_2)} .
\end{equation*}
This leads to
\begin{align*}
	\tau_j^{\nu_1+m} \sigma_j^{\nu_2- m} 
	&\leq  \max\{1, 2^{\gamma(m - \nu_2)}\} \gamma^{\nu_1+m} {T^{ \nu_1+\nu_2}}
	\begin{cases}
		N^{-\gamma(\nu_1+\nu_2) }, & \gamma (\nu_1+\nu_2) < m, \\[0.5em]
		N^{-m}, & \gamma  (\nu_1+\nu_2) \geq m,
	\end{cases}
\end{align*}
and
\begin{align*}
	\sum_{j=2}^n \tau_j^m \sigma_j^{\nu - m} 
	&\leq \max\{1, 2^{\gamma(m - \nu)}\} \gamma^m \tau^{\gamma \nu}  \sum_{j=2}^n j^{\gamma \nu - m} \\
	&\le C \max\{1, 2^{\gamma(m - \nu)}\} \gamma^m T^\nu
	\begin{cases}
		N^{-\gamma \nu}, & \gamma \nu < m - 1, \\[0.5em]
		N^{-m + 1} \log (n), & \gamma \nu = m - 1, \\[0.5em]
		N^{-m + 1}, & \gamma \nu > m - 1.
	\end{cases}
\end{align*}
This completes the proof.
\end{proof}
\bibliographystyle{abbrv}
\bibliography{gCQRK_sectorial_R2}

@article {McSloTho2006,
	AUTHOR = {McLean, William and Sloan, Ian H. and Thom\'ee, Vidar},
	TITLE = {Time discretization via {L}aplace transformation of an
	integro-differential equation of parabolic type},
	JOURNAL = {Numer. Math.},
	FJOURNAL = {Numerische Mathematik},
	VOLUME = {102},
	YEAR = {2006},
	NUMBER = {3},
	PAGES = {497--522},
	ISSN = {0029-599X,0945-3245},
	MRCLASS = {65R20 (44A10 45K05 65D32)},
	MRNUMBER = {2207270},
	MRREVIEWER = {Ezio\ Venturino},
	DOI = {10.1007/s00211-005-0657-7},
	URL = {https://doi.org/10.1007/s00211-005-0657-7},
}

@article {LuScha2002,
	AUTHOR = {Lubich, Christian and Sch\"adle, Achim},
	TITLE = {Fast convolution for nonreflecting boundary conditions},
	JOURNAL = {SIAM J. Sci. Comput.},
	FJOURNAL = {SIAM Journal on Scientific Computing},
	VOLUME = {24},
	YEAR = {2002},
	NUMBER = {1},
	PAGES = {161--182},
	ISSN = {1064-8275,1095-7197},
	MRCLASS = {44A35 (35L05 35Q40 44A10 65R10)},
	MRNUMBER = {1924419},
	DOI = {10.1137/S1064827501388741},
	URL = {https://doi.org/10.1137/S1064827501388741},
}

@article {JingRebecca_fi,
	AUTHOR = {Li, Jing-Rebecca},
	TITLE = {A fast time stepping method for evaluating fractional
	integrals},
	JOURNAL = {SIAM J. Sci. Comput.},
	FJOURNAL = {SIAM Journal on Scientific Computing},
	VOLUME = {31},
	YEAR = {2009/10},
	NUMBER = {6},
	PAGES = {4696--4714},
	ISSN = {1064-8275,1095-7197},
	MRCLASS = {65D30 (26A33 41A55)},
	MRNUMBER = {2594999},
	MRREVIEWER = {Kai\ Diethelm},
	DOI = {10.1137/080736533},
	URL = {https://doi.org/10.1137/080736533},
}

@article {BanMak2022,
	AUTHOR = {Banjai, Lehel and Makridakis, Charalambos G.},
	TITLE = {A posteriori error analysis for approximations of
	time-fractional subdiffusion problems},
	JOURNAL = {Math. Comp.},
	FJOURNAL = {Mathematics of Computation},
	VOLUME = {91},
	YEAR = {2022},
	NUMBER = {336},
	PAGES = {1711--1737},
	ISSN = {0025-5718},
	MRCLASS = {65M06 (35R11 65M12 65M15)},
	MRNUMBER = {4435945},
	MRREVIEWER = {Luis Vazquez},
	DOI = {10.1090/mcom/3723},
	URL = {https://doi.org/10.1090/mcom/3723},
}

@article{LoSau13,
	Author = {Lopez-Fernandez, Maria and Sauter, Stefan},
	Doi = {10.1093/imanum/drs034},
	Fjournal = {IMA Journal of Numerical Analysis},
	Issn = {0272-4979},
	Journal = {IMA J. Numer. Anal.},
	Mrclass = {65R20 (44A35)},
	Mrnumber = {3119712},
	Number = {4},
	Pages = {1156--1175},
	Title = {Generalized convolution quadrature with variable time stepping},
	Url = {https://doi.org/10.1093/imanum/drs034},
	Volume = {33},
	Year = {2013}}

@article{LoSau16,
	Author = {Lopez-Fernandez, M. and Sauter, S.},
	Doi = {10.1007/s00211-015-0761-2},
	Fjournal = {Numerische Mathematik},
	Issn = {0029-599X},
	Journal = {Numer. Math.},
	Mrclass = {65R20},
	Mrnumber = {3520008},
	Mrreviewer = {Kai Diethelm},
	Number = {4},
	Pages = {743--779},
	Title = {Generalized convolution quadrature based on {R}unge-{K}utta methods},
	Url = {https://doi.org/10.1007/s00211-015-0761-2},
	Volume = {133},
	Year = {2016}}

@article{Lu86,
	Author = {Lubich, Ch},
	Journal = {SIAM J. Math. Anal.},
	Number = {3},
	Pages = {704--719},
	Publisher = {SIAM},
	Title = {Discretized fractional calculus},
	Volume = {17},
	Year = {1986}}

@article{Gar15,
	Author = {Garrappa, Roberto},
JOURNAL = {SIAM J. Numer. Anal.},
FJOURNAL = {SIAM Journal on Numerical Analysis},
	Number = {3},
	Pages = {1350--1369},
	Publisher = {SIAM},
	Title = {{Numerical evaluation of two and three parameter Mittag-Leffler functions}},
	Volume = {53},
	Year = {2015}}

@article{BanLo20,
	Author = {Banjai, L. and L\'{o}pez-Fern\'{a}ndez, M.},
	Doi = {10.1016/j.jcp.2019.109155},
	Fjournal = {Journal of Computational Physics},
	Issn = {0021-9991},
	Journal = {J. Comput. Phys.},
	Mrclass = {65M99},
	Mrnumber = {4046171},
	Pages = {109155, 21},
	Title = {Numerical approximation of the {S}chr\"{o}dinger equation with concentrated potential},
	Url = {https://doi.org/10.1016/j.jcp.2019.109155},
	Volume = {405},
	Year = {2020}}

@article{BanLo19,
	Author = {Banjai, Lehel and L{\'o}pez-Fern{\'a}ndez, Mar{\'\i}a},
	Doi = {10.1007/s00211-018-1004-0},
	Fjournal = {Numerische Mathematik},
	Issn = {0029-599X},
	Journal = {Numer. Math.},
	Mrclass = {65M60 (26A33 35R11 65D30 65L06 65M15)},
	Mrnumber = {3905428},
	Mrreviewer = {Eihab B. M. Bashier},
	Number = {2},
	Pages = {289--317},
	Title = {Efficient high order algorithms for fractional integrals and fractional differential equations},
	Url = {https://doi.org/10.1007/s00211-018-1004-0},
	Volume = {141},
	Year = {2019}}

@article{LoSau15apnum,
	Author = {Lopez-Fernandez, Maria and Sauter, Stefan},
	Doi = {10.1016/j.apnum.2015.03.004},
	Fjournal = {Applied Numerical Mathematics. An IMACS Journal},
	Issn = {0168-9274},
	Journal = {Appl. Numer. Math.},
	Mrclass = {65R20 (44A35)},
	Mrnumber = {3342585},
	Mrreviewer = {Eleonora Messina},
	Pages = {88--105},
	Title = {Generalized convolution quadrature with variable time stepping. {P}art {II}: {A}lgorithm and numerical results},
	Url = {https://doi.org/10.1016/j.apnum.2015.03.004},
	Volume = {94},
	Year = {2015}}

@article{BanLuMe,
	Author = {Banjai, Lehel and Lubich, Christian and Melenk, Jens Markus},
	Doi = {10.1007/s00211-011-0378-z},
	Fjournal = {Numerische Mathematik},
	Issn = {0029-599X},
	Journal = {Numer. Math.},
	Mrclass = {65D30 (65E05 65R10)},
	Mrnumber = {2824853},
	Mrreviewer = {Michael J. Carley},
	Number = {1},
	Pages = {1--20},
	Title = {{R}unge-{K}utta convolution quadrature for operators arising in wave propagation},
	Url = {https://doi.org/10.1007/s00211-011-0378-z},
	Volume = {119},
	Year = {2011}}

@article{CuLuPa,
	Author = {Cuesta, Eduardo and Lubich, Christian and Palencia, Cesar},
	Doi = {10.1090/S0025-5718-06-01788-1},
	Fjournal = {Mathematics of Computation},
	Issn = {0025-5718},
	Journal = {Math. Comp.},
	Mrclass = {65R20 (26A33 35R10 45K05)},
	Mrnumber = {2196986},
	Mrreviewer = {Axel Maas},
	Number = {254},
	Pages = {673--696},
	Title = {Convolution quadrature time discretization of fractional diffusion-wave equations},
	Url = {https://doi.org/10.1090/S0025-5718-06-01788-1},
	Volume = {75},
	Year = {2006}}

@article{Lu94,
	Author = {Lubich, Ch.},
	Doi = {10.1007/s002110050033},
	Fjournal = {Numerische Mathematik},
	Issn = {0029-599X},
	Journal = {Numer. Math.},
	Mrclass = {65M15 (65D30 65M60 65R20)},
	Mrnumber = {1269502},
	Mrreviewer = {Rudolf Gorenflo},
	Number = {3},
	Pages = {365--389},
	Title = {On the multistep time discretization of linear initial-boundary value problems and their boundary integral equations},
	Url = {https://doi.org/10.1007/s002110050033},
	Volume = {67},
	Year = {1994}}

@article{Lu2004,
	Author = {Lubich, Christian},
	Doi = {10.1023/B:BITN.0000046813.23911.2d},
	Fjournal = {BIT. Numerical Mathematics},
	Issn = {0006-3835},
	Journal = {BIT},
	Mrclass = {65R10 (41A55)},
	Mrnumber = {2106013},
	Number = {3},
	Pages = {503--514},
	Title = {Convolution quadrature revisited},
	Url = {https://doi.org/10.1023/B:BITN.0000046813.23911.2d},
	Volume = {44},
	Year = {2004}}

@article{LuOs,
	Author = {Lubich, Ch. and Ostermann, A.},
	Doi = {10.2307/2153158},
	Fjournal = {Mathematics of Computation},
	Issn = {0025-5718},
	Journal = {Math. Comp.},
	Mrclass = {65M12 (65D30 65M15 65R20 76D05 76M25)},
	Mrnumber = {1153166},
	Mrreviewer = {John C. Strikwerda},
	Number = {201},
	Pages = {105--131},
	Title = {Runge-{K}utta methods for parabolic equations and convolution quadrature},
	Url = {https://doi.org/10.2307/2153158},
	Volume = {60},
	Year = {1993}}

@article{Lu88II,
	Author = {Lubich, C.},
	Doi = {10.1007/BF01462237},
	Fjournal = {Numerische Mathematik},
	Issn = {0029-599X},
	Journal = {Numer. Math.},
	Mrclass = {65D30 (65R20)},
	Mrnumber = {932708},
	Mrreviewer = {Syvert P. N\o rsett},
	Number = {4},
	Pages = {413--425},
	Title = {Convolution quadrature and discretized operational calculus. {II}},
	Url = {https://doi.org/10.1007/BF01462237},
	Volume = {52},
	Year = {1988}}

@article{Lu88I,
	Author = {Lubich, C.},
	Doi = {10.1007/BF01398686},
	Fjournal = {Numerische Mathematik},
	Issn = {0029-599X},
	Journal = {Numer. Math.},
	Mrclass = {65D30 (44A55 65L05 65R20)},
	Mrnumber = {923707},
	Mrreviewer = {Syvert P. N\o rsett},
	Number = {2},
	Pages = {129--145},
	Title = {Convolution quadrature and discretized operational calculus. {I}},
	Url = {https://doi.org/10.1007/BF01398686},
	Volume = {52},
	Year = {1988}}

@book{Henrici_II,
	Address = {New York},
	Author = {Henrici, Peter},
	Isbn = {0-471-01525-3},
	Mrclass = {30-02 (33-02 44A10)},
	Mrnumber = {0453984 (56 \#12235)},
	Mrreviewer = {A. E. Heins},
	Note = {Special functions---integral transforms---asymptotics---continued fractions},
	Pages = {ix+662},
	Publisher = {Wiley Interscience [John Wiley \& Sons]},
	Title = {Applied and computational complex analysis. {V}ol. 2},
	Year = {1977}}

@book{Gorenflo2020,
	Author = {Gorenflo, Rudolf and Kilbas, Anatoly A. and Mainardi, Francesco and Rogosin, Sergei},
	Doi = {10.1007/978-3-662-61550-8},
	Isbn = {978-3-662-61550-8; 978-3-662-61549-2},
	Mrclass = {33-02 (26A33 33C60 33E12 34A08 44A20 45K05 60G22)},
	Mrnumber = {4179587},
	Note = {Second edition [of 3244285]},
	Pages = {xvi+540},
	Publisher = {Springer, Berlin},
	Series = {Springer Monographs in Mathematics},
	Title = {Mittag-{L}effler functions, related topics and applications},
	Url = {https://doi.org/10.1007/978-3-662-61550-8},
	Year = {2020}}

@book{Podlubny,
	Author = {Podlubny, Igor},
	Isbn = {0-12-558840-2},
	Mrclass = {26A33 (34K05)},
	Mrnumber = {1658022},
	Mrreviewer = {Anatoly Kilbas},
	Note = {An introduction to fractional derivatives, fractional differential equations, to methods of their solution and some of their applications},
	Pages = {xxiv+340},
	Publisher = {Academic Press, Inc., San Diego, CA},
	Series = {Mathematics in Science and Engineering},
	Title = {Fractional differential equations},
	Volume = {198},
	Year = {1999}}

@article{HipLoPa2014,
	Author = {Hiptmair, Ralf and Lopez-Fernandez, Maria and Paganini, Alberto},
	Doi = {10.1016/j.cam.2013.12.025},
	Fjournal = {Journal of Computational and Applied Mathematics},
	Issn = {0377-0427},
	Journal = {J. Comput. Appl. Math.},
	Mrclass = {65M60 (65R20 78M10)},
	Mrnumber = {3162369},
	Pages = {500--517},
	Title = {Fast convolution quadrature based impedance boundary conditions},
	Url = {https://doi.org/10.1016/j.cam.2013.12.025},
	Volume = {263},
	Year = {2014}}

@book{BanSay22,
	title={Integral Equation Methods for Evolutionary PDE: A Convolution Quadrature Approach},
	author={Banjai, Lehel and Sayas, Francisco-Javier},
	volume={59},
	year={2022},
	publisher={Springer Nature}
}

@book{SamkoKilbasMarichev,
	title={Fractional integrals and derivatives: theory and applications},
	author={Samko, S. G. and Kilbas, A. A. and Marichev, O. I.},
	year={1993},
	publisher={CRC Press},
	address={Basel}
}

@article {BaBanPtas24,
		AUTHOR = {Baker, Katherine and Banjai, Lehel and Ptashnyk, Mariya},
		TITLE = {{Numerical analysis of a time-stepping method for the Westervelt equation with time-fractional damping}},
		JOURNAL = {Math. Comp.},
		FJOURNAL = {Mathematics of Computation},
		VOLUME = {},
		YEAR = {2024},
		PAGES = {},
		ISSN = {},
		MRCLASS = {},
		MRNUMBER = {},
		DOI = {10.1090/mcom/3945},
		URL = {https://doi.org/10.1090/mcom/3945},
	}

@book{HairWan10,
		author    = {Ernst Hairer and Gerhard Wanner},
		title     = {Solving Ordinary Differential Equations {II}, Second Revised Edition},
		year      = {2010},
		publisher = {Springer-Verlag},
		address   = {Berlin, Heidelberg},
		doi       = {10.1007/978-3-642-05199-9},
		isbn      = {978-3-642-05198-2},
	}

@article{GuoLo,
  title={Generalized convolution quadrature for non smooth sectorial problems},
author={Guo, Jing and Lopez-Fernandez, Maria},
journal={Calcolo},
volume={62},
number={1},
pages={1--37},
year={2025},
publisher={Springer}
}

@article {LoLuPaScha,
	AUTHOR = {L\'{o}pez-Fern\'{a}ndez, Mar\'{\i}a and Lubich, Christian and Palencia,
	Cesar and Sch\"{a}dle, Achim},
	TITLE = {Fast {R}unge-{K}utta approximation of inhomogeneous parabolic
	equations},
	JOURNAL = {Numer. Math.},
	FJOURNAL = {Numerische Mathematik},
	VOLUME = {102},
	YEAR = {2005},
	NUMBER = {2},
	PAGES = {277--291},
	ISSN = {0029-599X},
	MRCLASS = {65M20},
	MRNUMBER = {2206466},
	MRREVIEWER = {Othmar Koch},
	DOI = {10.1007/s00211-005-0624-3},
	URL = {https://doi.org/10.1007/s00211-005-0624-3},
}

@phdthesis{Maria_thesis,
	title={Discretizaciones de orden espectral de integrales de contorno sectoriales y aplicaciones a problemas de evoluci{\'o}n},
	author={L\'{o}pez-Fern\'{a}ndez, M.},
	year={2005},
	school={Universidad de Valladolid}
}

@article{JinLiZh,
	title={{Correction of high-order BDF convolution quadrature for fractional evolution equations}},
	author={Jin, Bangti and Li, Buyang and Zhou, Zhi},
	JOURNAL = {SIAM J. Sci. Comput.},
	FJOURNAL = {SIAM Journal on Scientific Computing},
	volume={39},
	number={6},
	pages={A3129--A3152},
	year={2017},
	publisher={SIAM}
}

@article{BanFer,
	title={{Runge--Kutta convolution quadrature based on Gauss methods}},
	author={Banjai, Lehel and Ferrari, Matteo},
        JOURNAL = {Numer. Math.},
	FJOURNAL = {Numerische Mathematik},
	volume={156},
	number={5},
	pages={1719--1750},
	year={2024},
	publisher={Springer}
}

@article{BanLu11,
	title={{An error analysis of Runge--Kutta convolution quadrature}},
	author={Banjai, Lehel and Lubich, Christian},
	Fjournal={BIT Numerical Mathematics},
Journal={BIT},
	volume={51},
	pages={483--496},
	year={2011},
	publisher={Springer}
}

@incollection{SautSch,
	title={Boundary element methods},
	author={Sauter, Stefan A and Schwab, Christoph},
	booktitle={Boundary Element Methods},
	pages={183--287},
	year={2010},
	publisher={Springer}
}

@article{BanFer23,
title={Generalized convolution quadrature based on the trapezoidal rule},
author={Banjai, Lehel and Ferrari, Matteo},
journal={IMA Journal of Numerical Analysis},
pages={draf141},
year={2026},
publisher={Oxford University Press}
}

@article{ZhuJuZhao,
  title={{Fast high-order compact exponential time differencing Runge--Kutta methods for second-order semilinear parabolic equations}},
  author={Zhu, Liyong and Ju, Lili and Zhao, Weidong},
FJOURNAL = {Journal of Scientific Computing},
JOURNAL = {J. Sci. Comput.},
  volume={67},
  pages={1043--1065},
  year={2016},
  publisher={Springer}
}

@article{LiMa,
  title={Exponential convolution quadrature for nonlinear subdiffusion equations with nonsmooth initial data},
  author={Li, Buyang and Ma, Shu},
  Fjournal = {SIAM Journal on Numerical Analysis},
  Journal = {SIAM J. Numer. Anal.},
  volume={60},
  number={2},
  pages={503--528},
  year={2022},
  publisher={SIAM}
}

@article{LiMaRa,
author = {Li, Buyang and Lin, Yanping and Ma, Shu and Rao, Qiqi},
title = {An Exponential Spectral Method Using VP Means for Semilinear Subdiffusion Equations with Rough Data},
journal = {SIAM Journal on Numerical Analysis},
volume = {61},
number = {5},
pages = {2305-2326},
year = {2023},
doi = {10.1137/22M1512041}
}

@article{Steng,
	title={Numerical methods based on Whittaker cardinal, or sinc functions},
	author={Stenger, Frank},
	journal={Siam Review},
	volume={23},
	number={2},
	pages={165--224},
	year={1981},
	publisher={SIAM}
}

@misc{CodesGuoLo25Zenodo,
	author       = {Guo, Jing and Lopez-Fernandez, Maria},
	title        = {{Code accompanying the paper ``Runge--Kutta generalized Convolution Quadrature for sectorial problems''}},
	year         = {2026},
	version      = {v2025.06},
	publisher    = {Zenodo},
	doi          = {10.5281/zenodo.19482892},
	url          = {https://doi.org/10.5281/zenodo.19482892},
	howpublished = {\url{https://doi.org/10.5281/zenodo.19482892}}
}

\end{document}